\documentclass[12pt]{book}

\usepackage[dvips]{graphics}

\begin{document}

\raggedbottom

\pagenumbering{roman}
\pagestyle{plain}

\begin{center}
\vspace*{.25in}
Stochastic Spatial Models of Plant Diseases \\
\vspace{.5in}
By \\
David Herbert Brown \\
B.A. (St. John's College, Santa Fe, NM) 1992 \\
\vspace{.3in}
DISSERTATION \\
\vspace{.3in}
Submitted in partial satisfaction of the requirements for the degree of \\
DOCTOR OF PHILOSOPHY \\
in \\
Applied Mathematics \\
in the \\
OFFICE OF GRADUATE STUDIES \\
of the \\
UNIVERSITY OF CALIFORNIA \\
DAVIS \\
\vspace{.2in}
Approved: \\
\vspace{.2in}
\rule{3in}{.01in} \\
\vspace{.2in}
\rule{3in}{.01in} \\
\vspace{.2in}
\rule{3in}{.01in} \\
Committee in Charge \\
2001

\end{center}

\newpage
\mbox{}
\thispagestyle{empty}
\newpage

\section*{Acknowledgements}

I could not have completed my Ph.D.\ without the support, guidance,
and friendship of many people; it is a pleasure to acknowledge a few
of them here.  First, I thank my parents for their patient support and
encouragement through these many years of schooling.  I cannot express
my gratitude to them for enabling me to pursue a career based solely on 
intellectual interest.  They have been my first and best teachers.
I thank my sister Heather for blazing the Ph.D.\ trail and dropping
breadcrumbs to help me find my way.

Six years ago, I did not know that the field of theoretical ecology
existed (it does).  My love of it arose largely because of the guidance
of my advisor, Alan Hastings.  His enthusiasm and expertise helped
me gain a foothold in a field that I expect will provide me with a
lifetime of challenges and rewards.  He is also a Zen master of thesis
advising; his subtle hints and mild suggestions kept me on course even
as I pursued questions of my own devising.  I thank Janko Gravner for
teaching me probability, and Carole Hom and Ed Caswell--Chen for their
careful comments on my manuscripts.

Finally, and most importantly, I thank my wife Lisa for her kindness,
intelligence, humor, and love.  I took a fair amount of putting up 
with during this process; her encouragement kept me going when I 
doubted the outcome.  I am forever in her debt.
Also, I thank John Elway for the Superbowl wins.  And whoever
invented doughnuts. 

\newpage
\mbox{}
\thispagestyle{empty}
\newpage

\tableofcontents

\newpage
\mbox{}
\thispagestyle{empty}
\newpage

\pagestyle{myheadings}
\markright{}
\pagenumbering{arabic}

\chapter[Introduction]{Introduction}
\thispagestyle{myheadings}
\markright{}

	The importance of infectious diseases in natural communities has
been recieving increased attention from ecologists over the past two decades
(Grenfell and Dobson, 1995).
Prior to this surge in interest, our understanding of disease in natural
populations lagged behind our understanding of other biotic forces such
as predation and competition.  Pathogenic microorganisms are hard to
detect, their effects on hosts may be subtle, and in complex natural
communities disease incidence is often patchy in space and time (Burdon, 1987).
Several factors contributed to the new focus on infectious disease as an
ecological factor.  Improvements in pathogen detection techniques revealed
the ubiquitous presence of potential disease agents in every kind of 
ecosystem (Dinoor and Eshed, 1984).  
At the same time, introductions of novel pathogens had
spectacular consequences, sometimes completely changing the structure
of the communities they invaded.  The effects of \emph{Phytophthora
cinnamomi} on Australian eucalypt forests (Weste and Marks, 1987), of 
chestnut blight on North American hardwood forests (Stephenson, 1986), and 
of avian malaria
on Hawaiian birds (Atkinson et al., 1995) illustrated the vast potential 
of pathogens for
impacting natural populations and communities.  The reemergence of
infectious diseases of humans, especially HIV, contributed to a general
interest in disease dynamics (Mollison, 1995).  This helped fuel the 
effort to model
disease processes in human and natural population; mechanistic (as
opposed to statistical) epidemiological modeling, which had been largely 
dormant since the foundational papers by Kermack and McKendrick (1927), 
was reinvigorated by the seminal work of Anderson and May (1979; May and
Anderson, 1979).

	Meanwhile, the importance of spatial structure in ecological 
processes has been a growing concern of theoretical ecologists 
(Tilman and Kareiva, 1997).  Biotic
interactions are now understood to be scale dependent (Levin, 1992), 
with spatial
structure being an important factor at a number of scales, from that
of the individual to populations, metapopulations, landscapes, and 
continents.  Each of these spatial scales has been incorporated in 
the effort to develop the theory of disease ecology.  For example,
reaction--diffusion and integrodifferential models have been used to 
study the spread of 
rabies across Europe (Murray et al., 1986) and the spread of disease 
within agricultural fields (Zadoks, 2000).  Metapopulation and 
landscape models have been important in
explaining the maintenance of resistance polymorphisms, by showing that
host--pathogen coevolution can be decoupled by spatially induced
heterogeneity (Burdon and Thrall, 1999; Damgaard, 1999; Stahl et al., 1999).  
All of these modeling frameworks incorporate spatial
structure above the scale of the individual; they assume local mixing
of individuals in order to derive larger scale spatial patterns.  
However, Mollison (1977) recognized the importance of individual--scale 
spatial structure in disease ecology.  The progress of an epidemic depends on
transmission between individuals, and the effect of a disease on its
host may depend on the host's local environment.  Recent advances in fine 
scale spatial modeling have allowed us to begin to understand how
the spread and impact of a disease depends on spatial structure at the
scale of individuals.

	As Durrett and Levin (1994) pointed out, the importance of fine scale 
spatial structure is intertwined with the discreteness of individual
organisms and the stochasticity that characterizes biotic interactions.
Thus, a modeling framework has been developed that incorporates
discrete individuals and explicit spatial structure into a stochastic 
process.  When space is represented as a discrete lattice, these
models are known as probabilistic cellular automata or interacting
particle systems; when space is continuous, the models are called
point processes.  Epidemic models of this type have demonstrated the
importance of local spatial structure on the spread, impact, and
evolution of diseases.  Generally, the spread of an epidemic is much
slower in a stochastic spatial model than in an analogous nonspatial
model (Holmes, 1997; Bolker, 1999).  This occurs because local disease 
transmission tends to
deplete the local availability of susceptible individuals, so that the
realized transmission rate is much lower than what would be achieved
in a well--mixed population.  More interestingly, the spatial models
predict qualitative results that cannot occur in a nonspatial model.
For example, there may be an upper transmission threshold above which
the disease cannot spread (Rand et al., 1995; Holmes, 1997); if the disease 
depletes local susceptibles
too quickly, it may burn itself out.  Also, the disease may be able
to drive the host extinct (Sato et al., 1994).  This cannot happen in a 
simple nonspatial
model, because the disease always fades out before the host crashes 
(i.e. it cannot stabilize the $S = 0$ equilibrium).  However, spatial
clustering of individuals may allow transmission rates to remain
high even as the host population approaches zero, so that the disease
drives the population lower even at arbitrarily low density.
Finally, the evolution of virulence and resistance may be qualitatively
changed when host or pathogen dispersal is localized (Rand et al., 1995; 
Boots and Sasaki, 2000).  These results
suggest that disease ecology in natural populations, especially of
sessile organisms, may be highly dependent on spatial structure at
the scale of individuals.  Although these models usually do not incorporate
a high level of biological detail and realism, they indicate ways in
which the theoretical framework of disease ecology should be modified
from predictions made by earlier, nonspatial models.

	While stochastic spatial models reveal many important aspects
of local population interactions, the analysis of them is difficult.
The local interactions of discrete individuals that drive the models'
interesting behaviors prevent the type of continuum limit or averaging
methods that have been applied so successfully in other contexts.
For example, one may derive reaction--diffusion approximations to
interacting particle systems in the limit of rapid local stirring,
but the resulting equations do not always capture the features of
interest in the original model (Durrett and Levin, 1994).  A fair amount 
may be learned about
such models simply by simulating them on a computer.  However, exploration
of parameter space via simulation is tedious, and care must be taken to
ensure that threshold calculations are not an artifact of the size of 
the simulations or the timescales used.  Furthermore, simulations often
have little explanatory power; we would like to know not only what the model
does, but how its behavior can be explained in terms of the interaction
structure and parameter values.  For sufficiently simple lattice models,
a variety of probabilistic methods may yield exact results (Liggett, 1985;
Durrett, 1995).  These
techniques are especially useful in certain problems of 
statistical mechanics, where the physics of the system imposes a
structure (such as a conserved quantity or a minimization principle)
that simplifies analysis.  However, in ecological models this type
of structure is usually absent: quantities are not conserved, and
Hamiltonian functions are not specified \emph{a priori}.  As a result,
ecological modelers have begun to develop approximation techniques
that capture the local spatial structure of these models in a
simpler framework.  The methods of pair approximation and moment
closure allow the incorporation of local spatial structure for
interacting particle systems and point processes, respectively, 
into systems of differential equations.  Since their initial introduction
in the past decade, these techniques have been adapted to describe a
variety of ecological systems, culminating in a recent edited volume
largely dedicated to them (Dieckmann et al., 2000). 

\section*{Pair Approximation}

	The technique of pair approximation in lattice population models
was introduced by Matsuda \emph{et al.} (1992).  The authors studied
population growth on a regular lattice, with local reproduction and 
migration.  Competition and altruism were incorporated by allowing
death rates to depend on local population densities.  They studied the
model by writing down differential equations for $P_{\sigma}$, the
probability that a randomly chosen site is in state $\sigma$.  Because
of local interactions, changes in $P_{\sigma}$ are governed by the
``pair densities'' $P_{\sigma\sigma'}$, the probability that a randomly
chosen site is in state $\sigma$ and a randomly chosen neighbor of it
is in state $\sigma'$.  Under the mean field assumption, 
$P_{\sigma\sigma'} = P_{\sigma}P_{\sigma'}$.  However, this omits all
spatial structure from the system and does not approximate the behavior
of the model very well.  Thus, the authors proceeded by writing down
the differential equations for the pair densities themselves.  Now, 
the state of one site in a pair can change due to interaction with
the other site in the pair, or with a third (neighboring) site.  Thus,
changes in the pair densities depend on triplets of sites, 
$P_{\sigma\sigma'\sigma''}$.  In order to obtain a closed system, 
Matsuda \emph{et al.} (1992) used the approximation $P_{\sigma\sigma'\sigma''}
= \frac{P_{\sigma\sigma'}P_{\sigma'\sigma''}}{P_{\sigma'}}$.  In terms
of conditional probabilities, this can be written $P_{\sigma|\sigma'\sigma''}
= P_{\sigma|\sigma'}$.  This means that the state of a site depends on
its neighbors, but not on its neighbors' neighbors.  The resulting
differential equations incorporate an approximation to the local spatial
structure of the model.  The authors found that the pair approximation
equations performed better than the mean field equations, and in particular
were able to predict the success of altruist behaviors that cannot
persist in the absence of local spatial structure. 

	Since its introduction, the pair approximation technique has
undergone a number of modifications and has been used to study a wide
range of ecological phenomena (reviewed in Rand, 1999).  
Sato \emph{et al.} (1994) introduced an
``improved pair approximation'' to try to correct the fact that pair
approximations typically underestimate the amount of local clustering.
They argued that, when state $\sigma$ is rare, $P_{\sigma|\sigma'\sigma}
 >> P_{\sigma|\sigma'}$.  That is, knowing that another $\sigma$ is
nearby is significant, since it indicates that we are likely in a region
of higher $\sigma$ density.  To model this, they assumed that
$P_{\sigma|\sigma'\sigma''} = \epsilon P_{\sigma|\sigma'}$ for some
$\sigma'' \ne \sigma$ and $\epsilon < 1$.  From the conservation
rules for probabilities, it follows that $P_{\sigma|\sigma'\sigma} \ne 0$
even when  $P_{\sigma|\sigma'} = 0$.  This clustering of types at low
densities was crucial for explaining the ability of an infectious disease
to bring about host extinction in the lattice model.  Deviations from
the original pair approximation may also be expected when we are studying
irregular networks (such as may arise from social interactions) or
regular lattices in which interactions are not restricted to the four
cardinal directions (van Baalen, 2000).  In that case, triplets may 
form triangles, in
which they are all mutually neighbors.  When the three sites are connected
in a triangle, a more appropriate approximation is $P_{\sigma\sigma'\sigma''}
= \frac{P_{\sigma\sigma'}P_{\sigma'\sigma''}P_{\sigma\sigma''}}{
P_{\sigma}P_{\sigma'}P_{\sigma''}}$ (Rand, 1999).  

	Thus, the type of
pair approximation used depends on the expected influence of the triplet
members on each other; this in turn depends on the geometry of space or
the interaction network (reviewed in Rand, 1999; van Baalen, 2000).
Unfortunately, no definite prescription can be given
for the appropriate pair approximation in any given model.  Modifications
of the original approach are important in systems that do not have
a square lattice geometry; however, they may also be necessary when the
details of local clustering govern behaviors of interest (Sato et al., 1994).
Moreover, the justification for the pair approximation is entirely heuristic.
With any approximation, we would ideally like to have an \emph{a priori}
estimate of the magnitude of the error introduced (as in a Taylor
expansion).  Such estimates are not likely to be feasible with
pair approximations; they are not even possible for the mean field
approximation.  However, we at least know under what limiting behaviors
the mean field assumption becomes valid.  In the limit of long--distance
interactions or rapid local stirring, a lattice model will converge to
the solution of the mean field equations (Durrett, 1995).  
Are there analogous interaction
structures or types of stirring that will bring the lattice model into
agreement with pair approximations?  It is not entirely obvious how to
break down the triplet structure while maintaining pairwise correlations.
Without knowing how to bring models into agreement with pair approximations,
it is hard to predict the suitability of any pair approximation approach
for a given model.  For now, all we can do is employ heuristic arguments
and compare the predictions carefully with simulation results. 

	Nevertheless, pair approximations have proven useful
for investigating the effects of local structure on a variety of 
ecological processes.  They have yielded new insights into plant
competition (Harada and Iwasa, 1994; Takenaka et al., 1997), 
succession (Sato and Konno, 1995), 
and forest gap dynamics (Iwasa, 2000).  In disease
ecology, they have been used to explain periodic and chaotic fluctuations 
(Keeling et al., 1997; Rand, 1999), host
extinction (Sato et al., 1994), and the effect of space on the 
epidemic threshold (Keeling, 1999).
When used in an ESS (evolutionarily stable state) approach to studying
phenotypic evolution, pair approximations clarify the success of 
altruism (Harada et al., 1995; van Baalen, 1998) and the evolution 
of intermediate levels of pathogen
virulence and disease resistance (Rand et al., 1995; Boots and Sasaki, 2000).  
In an interesting non--equilibrium
study, Ellner et al. (1998) used pair approximations to determine the 
rate of spread of an invading organism.   
  
	A common goal of pair approximation studies is to determine
critical parameter values for the lattice models.  These critical values
are typically thresholds for the ability of a particular state to invade
the system from low density.  That is, we want to determine parameter 
values for which $\frac{\dot{P}_{\sigma}}{P_{\sigma}} = 0$ when 
$P_{\sigma} = 0$.  One way to do this would be to augment the equations for
the singleton densities with the pair equations, and study the dynamics
of the full system.  However, the resulting system is often cumbersome, 
and its behavior sheds little insight.  Matsuda \emph{et al.} (1992) introduced
a second assumption that simplifies the analysis while increasing
its explanatory power.  They observed that when $\sigma$ is rare, the 
dynamics of the quantities $P_{\sigma'|\sigma}$ are much faster than for
$P_{\sigma}$.  This occurs because the equations for the conditional
probabilities are governed by terms that need not be small when $\sigma$
is rare.  Thus, we can separate the timescale for the development of 
local spatial structure from that of the overall invasion.  By solving
$\dot{P}_{\sigma'|\sigma} = 0$ when $P_{\sigma} = 0$, we can obtain a
pseudoequilibrium description of the spatial structure early in the invasion.
We can then include these terms as fixed parameters in the singleton
equations.  The resulting system has the simple structure of the mean
field equations, but includes corrections due to local spatial structure.
Of course, this approximation breaks down as the invader achieves a
non--trivial density, but it provides a powerful way to study invasion
criteria.  Like the pair approximation itself, the separation of timescales
does not have a formal justification, but it is easily verified by
simulations of the full model.  

	The development of the pair approximation approach appears to
have been self--contained within the field of theoretical ecology.  This
is somewhat surprising, given the superficial similarity of lattice models
in ecology to those in statistical mechanics.  Indeed, the idea of pair
approximations is related to the renormalization group approach used widely
in physics (Wilson, 1975; Zhou et al., 1994; Newman and Watts, 1999).  
In this approach, one recursively ``coarse--grains'' a
model; i.e. one studies the model in terms of larger and larger ``blocks''.
If the transformation map that converts the model to the next--coarsest
grain converges to a fixed point when iterated, one obtains a scaling
relationship that can be used, for example, to detect the existence of
a phase transition, or to model processes at scales smaller than can
conveniently be simulated.  The pair approximation approach of studying
the dynamics of pairs of sites is similar, in moving the focus from 
individual sites to blocks.  However, there are fundamental differences
between the two approaches.  The renormalization group is an iterative
scheme for detecting scaling laws of the system.  Pair approximations
only carry out the coarse--graining process one step, after which a
heuristic closure is invoked.  Moreover, the pair approximation is not
motivated by a scaling argument, but by the pairwise interactions that
determine the dynamics.  Pair approximations are more closely related to
the local structure theory of deterministic cellular automata (Gutowitz
et al., 1987), in which the dynamics of blocks of sites are used to
characterize model behavior.  Describing the dynamics of ecological models
in terms of blocks larger than two, however, can be tedious.  The effort
may be deemed worthwhile when we are trying to determine the precise
algorithmic properties of deterministic cellular automata, but not when
we are studying the qualitative effects of local clustering in models
that are understood to be simplistic representations of nature.  Thus,
it is unclear to what extent the pair approximation approach will continue
to evolve in isolation or be informed by the use of block dynamics in
statistical mechanics and cellular automata.

\section*{Moment Closure}

	The technique of moment closure was introduced by Bolker and
Pacala (1997) to approximate the dynamics of point processes.  As with
lattice models, dynamics are governed by interactions between pairs of
individuals; the dynamics of the pairs are in turn dependent on triplets.
Again, the goal is to close the system at the level of pairs by 
approximating the triplet densities in terms of lower order terms.  Since
the kernels used to describe spatial interactions are generally nonzero
over infinite distances, each triplet now constitutes a mutually interacting
triangle.  As a result, there are now several \emph{a priori} plausible
ways to express triplet densities as additive or multiplicative combinations
of pair densities (Dieckmann and Law, 2000).  

	Each of the closure schemes
leads to a set of integrodifferential equations for the dynamics of the
singleton and pair densities.  However, different closure assumptions
lead to systems that differ in their structure and performance.  The original
technique of Bolker and Pacala (1997, 1999) uses an additive (linear) 
combination
of pair densities.  The resulting equations are linear in the pair terms,
and include convolution terms arising from interactions between the two
neighbors of the focal individual.  The linearity of the system allows
techniques such as Fourier analysis to be used, at least for special types
of kernels.   Multiplicative closures result in systems that are nonlinear
in the spatial (pair) terms, and thus cannot be solved analytically.  
The simplest multiplicative closure ignores interactions between the
neighbors of the focal site, and thus does not result in convolution
terms in the pair density equations.  In fact, the original pair approximation
can be seen as an implementation of this closure when the interaction
kernels are uniform across nearest neighbors in a lattice.  In general,
however, the moment closure approach differs from pair approximations
because it includes pairwise structure at all distances.  
 
	Moment closure analysis of point processes has not been applied
to as wide a variety of systems as pair approximation.  It has proven
useful in studying how plant communities are structured by competition
and dispersal phenomena occurring at different spatial scales
(Bolker and Pacala, 1997, 1999; Law and Dieckmann, 2000).  In addition,
Bolker (1999) used the additive moment closure to study the dynamics
of a simple epidemic in randomly distributed and clumped host populations. 
He was able to explain how the local depletion of susceptibles typically slows
the epidemic relative to mean field predictions; however, he found that
clustering of the hosts could accelerate the epidemic's early progress.
As with pair approximations, a separation of timescales may be invoked
to determine invasion criteria, since the conditional pairwise spatial
structure of the system typically develops much faster than global densities.
Convergence of the spatial structure to a low density pseudoequilibrium
is not guaranteed, however; the lack of such convergence prevented Bolker
(1999) from determining the effects of spatial structure on the epidemic
threshold.  Thus, the different closure assumptions yield systems that
differ in their accuracy, tractability, and convergence properties.  Since
there is no formal and little heuristic justification for choosing one
method over another, we are for now forced to proceed by trial and error,
choosing the method that best suits our needs for a given question
(Dieckmann and Law, 2000).

\section*{Thesis Organization}

	In this thesis I present three studies of plant--pathogen interactions.
They are motivated by a desire to understand how the kinds of subtle
ecological interactions and spatial structure found in natural systems
determine the effects of pathogens.  The models are also linked by a 
common approach; each uses a version of pair approximation or moment
closure along with a separation of timescales argument to determine the
effects of spatial clustering on threshold structure.  By computing the
spatial structure early in an invasion, I find explicit corrections to
mean field theory.  In each case, the resulting pair or moment equations
are too complex to be solved analytically.  Hence, all of my results
come in the form of numerical calculations of the structure of parameter
space.  Nevertheless, the approximation schemes provide a powerful
pseudoanalytic tool: they allow rapid threshold calculations and often
provide insight into the qualitative effects of local spatial structure.

	In the first chapter, I model a disease that is not directly
lethal to its host, but rather affects its ability to compete with
neighbors.  This was motivated by the observation that pathogens
in natural systems often do not cause catastrophic damage to their hosts
(Burdon, 1987; Dobson and Crawley, 1994).
As I show, such diseases can still have a major impact if they place their
hosts at a competitive disadvantage.  Since competition and disease 
transmission are both spatially localized processes, their interaction
may reasonably be assumed to depend on spatial structure.  Since my goal
was to study the generic properies of this kind of interaction, I 
used a simple lattice model.  I used a version of the improved pair
approximation (Sato et al., 1994) to capture the type of clustering that
might be important for events such as host extinction.  

	In the second chapter, I study the basic SIR epidemic point
process introduced by Bolker (1999).  My goal was to address the fundamental
topic of the epidemic threshold: how does the ability of a pathogen to
invade depend on its dispersal and on the spatial structure of the host
population?  As Bolker found, the additive moment closure approach does
not allow one to address this question because it fails to converge to
a pseudoequilibrium spatial structure early in the invasion.  Thus, I
implemented a multiplicative moment closure.  I found that the resulting
equations did have the desired convergence properties, which allowed me
to compute the threshold transmission rate as a function of spatial
parameters.  

	In the final chapter, I study a problem in the evolution of
pathogen resistance suggested by Rice and Westoby (1982).  
I present a model of two plant species that share
a pathogen.  When the two hosts are also competitors, evolution may
lead to non--resistance by a host that finds the disease to be a useful
weapon.  Again, since I was interested in the generic properties of
such a system, I used a lattice model.  I used the ordinary pair approximation
assumption to study phenotypic evolution via repeated invasions by
novel strains.  Here again, local spatial structure proved critical; 
the evolution of non--resistance is possible in a spatial model but not
in the mean field equations. 

\pagebreak

\chapter[Nonlethal Diseases and Competition]{Nonlethal Diseases
	and Competition}
\thispagestyle{myheadings}
\markright{}

\section*{Abstract}

	Plant diseases that weaken but do not kill their hosts can have a 
major impact on natural communities by changing competitive hierarchies.
This paper presents a simple, spatially explicit model of interspecific
competition in which the superior competitor loses its advantage when
infected.  Pair approximations are used to determine conditions for
invasion and coexistence. 
Several aspects of the model's behavior are driven by the local outbreak
nature of infectious diseases.  Thus, the impact of diseases may be
qualitatively different from other sources of competitive heterogeneity.

\pagebreak

\section*{Introduction}

	Infectious diseases and competition have each been recognized as
important biotic forces shaping natural plant communities (Dinoor and
Eshed 1984, Grace and Tilman 1990, Dobson and Crawley 1994).  Ecologists
are increasingly aware of important interactions between these processes.
Diseases
can affect competitive interactions in three distinct ways.  First, 
apparent competition can occur among host plant species that share
a pathogen (Holt and Pickering 1985, Hudson
and Greenman 1998, Yan 1996).  Second, the increased mortality or 
reduced fecundity
of infected plants can afford new opportunities for growth of other
genotypes or species, shifting the competitive relationships between
the species (Chilvers and Brittain 1972, Alexander and Holt 1998).  
In this case, the combined effects of the 
disease and
competition on an individual plant are additive.  Finally, infection and
competition may interact nonadditively to determine the fate of each
individual.  This occurs when the physiological effects of infection
alter the host's ability to compete for a particular resource (Ayres and
Paul 1990, Clay 1990, Alexander and Holt 1998).  In this case, the combined
effect of disease and competition on an individual can be greater than the
sum of the two factors considered separately.

Diseases
have many physiological impacts on their plants, such as stunting growth,
causing leaf drop, interfering with nutrient uptake, or reducing seed
viability (Burdon 1987).  In some cases, the primary result of these 
physiological
effects is to change the host's competitive relationships with its
neighbors (for example, Paul and Ayres 1990, Lively et al. 1995).  Thus, 
a pathogen can have an important effect on the
host population not because it directly kills or sterilizes the host,
but because it lowers its competitive strength.  Since competition is
ubiquitous in plant communities, and many diseases in natural systems are
nonlethal ``debilitators'' (Burdon 1987, Dobson and Crawley 1994), this 
phenomenon may be very common.
This raises a basic question: How are natural plant populations and
communities impacted by sublethal pathogens whose sole or primary
effect is to lower their hosts' competitive abilities?  Experimental 
studies have demonstrated the existence of the phenomenon in several
simple (usually agricultural) systems (Ayres and Paul 1990).  However, 
we have little
theoretical or empirical basis for understanding how sublethal diseases
help structure natural plant communities. 

	A number of studies have documented the influence of pathogens on
competitive interactions in plants (reviewed in Ayres and
Paul 1990, Clay 1990, Alexander and Holt 1998).  The abiliy of sublethal
diseases to affect host populations primarily through changes in competitive
strength has been established in several systems.  One set of studies
examined the effect of the fungus \underline{Colletotrichum coccodes} on 
competition between its host, velvetleaf (\underline{Abutilon theophrasti}), and
soybeans in greenhouse (DiTommaso and Watson 1995) and field (DiTommaso et al.
1996) experiments.
The fungus generally is not lethal to adult velvetleaf.  It causes lesions
on leaves, which then are shed prematurely.  This results in stunted growth,
from which infected plants can eventually recover (DiTommaso et al. 1996).  
In velvetleaf
monocultures, inocculation with the fungus had a limited effect on yield.
However, when velvetleaf was grown with soybeans, inocculation led to
significant decreases in velvetleaf seed yield and slight increases in 
soybean yield relative to mixed cultures without the disease.  
DiTommaso et al. (1996) attributed this result to the
stunted growth of infected velvetleaf, which allowed the usually 
slower--growing soybeans to overtop and shade their competitors.

Similarly, low
levels of infection by the fungus \underline{Mycocentrospora acerina} on the
weed \underline{Viola arvensis} decreased leaf area, stunting growth and
reducing the host's competitive strength relative to wheat (Lawrie et al.
1999).  In a study of intraspecific competition between willows 
(\underline{Salix viminalis}), Verwijst (1993) found that infection by the 
fungus \underline{Melampsora epitea} changed the competitive relationship between
susceptible and non--susceptible clones.  Infection stunted growth during
the first year; in subsequent years
previously infected stools died because of their inability to compete
for light even though pathogen levels declined.  Lively et al. (1995) 
found that infection with the rust fungus
\underline{Puccinia recondita} had a significant effect on its host, \underline{
Impatiens campensis}, only at high host densities.  In the absence of strong
competition, the disease had little effect on growth. 

	The most comprehensive set of studies has examined the 
physiological effects of the rust fungus
\underline{Puccinia lagenophorae} on groundsel (\underline{Senecio vulgaris})
(Paul and Ayres 1984), and
the resulting impact on intraspecific (Paul and Ayres 1986, 1987b) and 
interspecific (Paul and Ayres 1987a, Paul 1989, Paul and Ayres 1990) 
competition.  This foliar pathogen affects its host in a number of ways,
including reduced photosynthesis by infected leaves, increased water
loss through transpiration (Paul and Ayres 1984), and diminished water 
uptake by the roots (Paul and Ayres 1987b).  The net result is that 
the growth of infected plants is
stunted and their relative competitive strength can be diminished.
The change in relative competitive strength can come either through
reduced yield by the host (Paul and Ayres 1987b, 1990) or through 
increased growth and yield by the
competitor (Paul and Ayres 1987a, 1989).  Moreover, the impact depends 
on the limiting resource.
Under high nutrient conditions, groundsel is competitively superior to
\underline{Capsella bursa-pastoris}, but the hierarchy is reversed in low
nutrient conditions.  The pathogen eliminated groundsel's advantage
when nutrients were plentiful, but had no effect on competition when growth 
was nutrient limited (Paul and Ayres 1990).  The impact of the fungus on 
intraspecific competition
was increased by drought conditions, since infected individuals were not
able to tolerate water stress well (Paul and Ayres 1987b).  These studies 
demonstrate that
sublethal diseases can affect competition in a variety of ways that depend
on the nature of the competitive interactions and the physiological
changes induced by infection.
 
	We still know little about how sublethal diseases impact complex
natural plant
communities.  Several factors contribute to a lack of information about 
what is probably a common phenomenon (Dobson and Crawley 1994).  
First, the very presence of pathogens
that have subtle effects on their hosts usually will go unnoticed in 
complex natural systems, unless there is a dedicated search for them.
Note that most of the known examples to date involve foliar pathogens such
as rust fungi; root--born or systemic pathogens can be much harder to detect
(Burdon 1987).
Second, even if the presence of the pathogens is known, their 
effects may be difficult to measure against the background noise of other
factors that influence competitive outcomes.
Competitive interactions between plants have been shown to
be mediated by a number of biotic and abiotic factors such as herbivory,
genetic heterogeneity, soil chemistry, and microclimate (Grace and
Tilman 1990).  Every population
is heterogeneous in the competitive strengths of individuals, and any study of
competition either averages such heterogeneities across the entire population
or samples only a subset of the conditions in the natural system.   
The sublethal effects of
pathogens may be yet another source of heterogeneity in the host population,
if infected and uninfected individuals differ in their competitive strength.
If the primary object of study is the plant community, it could be argued
that sublethal pathogens are simply one source of heterogeneity among many,
and that gaining detailed knowledge of the pathogen dynamics would not
contribute materially to our understanding of the community structure. 
This raises a second basic question: For which
systems and which questions do we need to track pathogen dynamics explicitly
rather than subsuming the disease into overall host population heterogeneity? 
	
	There are several a priori reasons why we may need to treat
sublethal infectious diseases
differently from other sources of host heterogeneity.
First, the existence of a population threshold below which the disease
fades out could introduce a structure not found with other factors (Kermack
and McKendrick 1927, Anderson and May 1979, May and Anderson 1979).  Second,
disease dynamics are often characterized by periodic outbreaks, giving the
system an internal source of cyclicity (May and Anderson 1979).  Third, 
localized disease
transmission can introduce a spatial structure, as outbreaks are clustered
within the host population (Real and McElhany 1996, Thrall and Burdon 1999). 
This spatial structure may interact with 
spatially localized competition in non-trivial ways, so that explicit
knowledge of the spatiotemporal dynamics of the disease is needed.

	In order to establish a baseline understanding of the dynamics
of plant competition mediated by sublethal diseases, I have studied a
simple, spatially explicit model of such a system.  The inclusion of 
space in the model make analysis more difficult.  However, since spatial
structure is important in determining the outcome of plant competition
(Durrett and Levin 1998, Bolker and Pacala 1999) and epidemics (Jeger 1989),
an understanding of the interaction between these processes must include
an account of the role of space.  

The model is in the
form of an interacting particle system; space is represented as an
infinite lattice,
and the dynamics evolve in continuous time according to stochastic 
interaction rules.  Interacting particle systems have proved useful in 
investigating a number of
basic ecological phenomena (Durrett and Levin 1994a,b).  Because populations
are discrete and the model is stochastic, processes such as disease outbreaks
that depend on spatially localized interactions are modeled robustly. However, 
use of a discrete lattice introduces an artificial geometry that makes
detailed realistic models of natural systems impossible.  In addition,
the use of continuous time makes realistic inclusion of age classes
or seasonality difficult.  Finally, the use of an infinite lattice masks
many aspects of stochasticity that may be biologically relevant.  One can
think of the infinite system as containing arbitrarily many independent 
copies of
the stochastic process.  Thus, a priori calculations such as the
expected density of occupied sites are reflected in any single realization
of the process (a feature known as ergodicity).  In a finite stochastic
process, different realizations can yield qualitatively different 
results.  
However, stochastic issues (such as whether an invasion or extinction
will occur) are essentially deterministic in the infinite particle system. 
Nevertheless, the model framework allows one to
take advantage of a growing body of theory on interacting particle
systems (Liggett 1985, Durrett 1995). 
While interacting particle systems do not yield detailed predictions
for particular populations, they are a useful framework for
investigating basic processes.  

I used a version of pair approximation
(Matsuda et al. 1992, Rand 1999) to study the model's dynamics.  This
technique allows one to derive a set of differential equations that
incorporate local spatial structure of the model.  By studying these
equations, one can efficiently map a model's parameter space and
gain insight into how local spatial structure drives dynamics.

	This simple, spatially explicit model
complements the only other theoretical work on such systems.  Gates et
al. (1986) studied a
deterministic, nonspatial model that tracked the amount of infected 
tissue for a 
foliar pathogen, using continuous growth within a growing season and
discrete yearly reproduction. In their model, diseased plant material 
stopped growing,
did not compete for resources, and had diminished reproductive yield.  
They found that such a disease can change the conditions for coexistence
of the plant species (or genotypes), 
and can lead to periodic or chaotic dynamics.  I have chosen 
to omit most of the physiological details they included, such as resource
flux, in order to focus
on the basic issues arising from the interaction of the spatially localized
processes of competition, reproduction, and transmission.

	I address four main questions in this study:

\begin{enumerate} 

	\item How does a nonlethal disease of a superior competitor affect
		the ability of an inferior competitor to invade?
	\item How does the presence of a competitor affect the persistence
		of the disease?
	\item Can  a nonlethal disease lead to host extinction?
	\item Does the impact of such a disease differ qualitatively from
		other sources of competitive heterogeneity?
\end{enumerate}
I answer these questions by examining three aspects of the model's 
behavior: the dependence of equilibrium population levels on parameter
values, critical parameter values for invasions and coexistence, and
spatio--temporal patterns.  The exploration of parameter space is
greatly facilitated by a simplifying technique known as pair approximation.
In order to answer the fourth question, I compare the results of the
disease model to an alternate one, in which ``host'' competitive strength
is determined at birth rather than through infection.

\section*{Model Structure}

	I developed a model for two plant species (or genotypes), a 
competitively
dominantly host (H) and an inferior competitor (C).  The host population is
subdivided into susceptible (S) and infected (I) individuals.  Each site
on an infinite  square lattice is occupied by one of these types or is 
empty (E).
All interactions occur between a site and the four adjacent sites.
Both species
reproduce by placing offspring onto neighboring empty sites, and all
individuals die at density--independent rates.  
Competition is manifested
by displacement: the superior competitor can place offspring on sites
occupied by the inferior competitor, displacing them. 
Infected hosts can transmit the disease to neighboring susceptible hosts.
I assume no vertical transmission, so that the offspring of infected
hosts are not infected.
The disease has
two effects on its host: infected hosts displace competitors at a reduced
rate, and are themselves subject to displacement by the competitor.
These disease effects roughly correspond to the experimental result that
infection can either inhibit growth by the host or promote growth by the
competitor (Ayres and Paul 1990).  In the absence of interspecific
competition, the disease has no effect on the host.  

The use of a single neighborhood structure for all interactions
is an important assumption, since differences
between the dispersal scales of the the plant species and the disease
may lead to qualitatively different results (Thrall and Burdon 1999).  
However, it greatly simplifies
the pair approximation technique, allowing a basic
understanding of the role that clustering plays in such systems. 

\subsection*{Disease Model}

	The structure of the disease model is summarized by the transitions
that can occur at a site (Figure 1a).
The probabilites that a randomly chosen site is in a particular state
(S, I, C, E) satisfy a set of ordinary differential equations.  Let
$P_{\sigma}$ be the probability that a random site is in state $\sigma$,
and $P_{\sigma | \sigma'}$ be the conditional probability that, given a
site is in state $\sigma'$, a randomly chosen neighbor is in state $\sigma$.
Then by considering the possible changes to a state, and the rates at
which they occur, we derive the ``master equations'':
\begin{eqnarray}
\dot{P}_{S} & = & (-\mu_{1} + \beta_{1} P_{E|S} + \alpha P_{C|S} - \lambda
	P_{I|S}) P_{S} + (\beta_{1} P_{E|I} + \hat{\alpha} P_{C|I}) P_{I} \\
\dot{P}_{I} & = & (-\mu_{1} + \lambda P_{S|I} - \gamma P_{C|I}) P_{I} \\
\dot{P}_{C} & = &  (-\mu_{2} + \beta_{2} P_{E|C} + \gamma P_{I|C}  
	 - \alpha P_{S|C} - \hat{\alpha} P_{I|C}) P_{C} \\
P_{E} & = & 1 - P_{S} - P_{I} - P_{C}. 
\end{eqnarray}
Here, $\mu_{1}$, and $\mu_{2}$ are the
death rates of the host and competitor, respectively.  Hosts place
an offspring on each neighboring empty site at rate $\beta_{1}/4$,
for a maximum birthrate of $\beta_{1}$ when surrounded by empty sites.
Healthy hosts displace
competitors at rate $\alpha/4$; infected hosts do so at rate $\hat{\alpha}/4$.
Of course, $\beta_{1} \ge \alpha \ge \hat{\alpha}$.  The transmission rate
of the disease between an infected plant and each susceptible neighbor 
is $\lambda/4$, and the competitor can displace infected
hosts at rate $\gamma/4 \le \beta_{2}/4$.

	In order to reduce the number of parameters, we can rescale time
by setting $\mu_{1} = 1$, so that one time unit is the expected lifetime
of a host.  In addition, I assume from now on that $\mu_{2}
= 1$, so that both plants have the same lifespan.  Because the disease
affects competition in two ways, we can clarify its role by introducing
parameters that measure its impact on competition.  
Let $\delta_{1}$ be the amount by which the disease
reduces the ability of the host to displace the competitor, and $\delta_{2}$
is the amount by which the disease increases the ability of the competitor
to displace the host.  That is, 
$\delta_{1} = (\alpha - \hat{\alpha})/\alpha$ and $\delta_{2}
= \gamma/\beta_{2}$, so that
 $0 \le \delta_{1}, \delta_{2} \le 1$.  When $\delta_{1} = 1$, the disease
completely prevents infected individuals from displacing competitors.
When $\delta_{2} = 1$, sites occupied by infected individuals are colonized
by the competitor as easily as empty sites.  The model
now involves three growth rates ($\beta_{1}, \beta_{2}, \lambda$) and
two disease effects ($\delta_{1}, \delta_{2}$).  

I note in passing that
this model reduces in a special case to the cyclic biased voter model, a
canonical model from the probability literature.  We obtain this
case by setting $\mu_{1} = \mu_{2} = 0$, so that there are no empty sites,
and $\delta_{1} = \delta_{2} = 1$.  Then the only possible transitions
are $S \rightarrow I \rightarrow C \rightarrow S$.  It is believed that
the cyclic biased voter model exhibits coexistence of all types for all
parameter values (Durrett 1999), although ``patient'' types (with low 
growth rates) have an advantage and tend to win out in finite simulations.

\subsection*{Alternate Model}

	In order to compare the effects of a pathogen to other sources of
competitive heterogeneity, I analyzed an alternate model from which the
disease is absent (Figure 1b).  In the alternate
model, the ``host'' species still has competitively strong and weak 
individuals (analogous to S and I in the disease model).
However, competitive strength is determined at birth rather than through 
infection.
Here, I let each offspring of a host be born into the class of inferior
competitors with some probability, $\rho$, 
independent of the status of the parent.  This
model uses the simplest way of
generating competitive heterogeneity; it does not incorporate a relationship
between the competitive strength of the parent and that of the offspring, as
may arise from genetic or environmental causes.
Thinking of S as strong individuals and I as inferior ones, we have:
\begin{eqnarray}
\dot{P}_{S} & = & (-\mu_{1} + \beta_{1}(1-\rho) P_{E|S} + \alpha(1-\rho) 
	P_{C|S}) P_{S} \nonumber \\
	& & + (\beta_{1}(1-\rho) P_{E|I} + \hat{\alpha}(1-\rho)
	P_{C|I}) P_{I} \\
\dot{P}_{I} & = & (-\mu_{1} + \beta_{1} \rho P_{E|I} + \hat{\alpha} \rho
	P_{C|I} - \gamma P_{C|I}) P_{I} \nonumber \\
	& & + (\beta_{1} \rho P_{E|S} + \alpha \rho P_{C|S}) P_{S}\\
\dot{P}_{C} & = &  (-\mu_{2} + \beta_{2} P_{E|C} + \gamma P_{I|C}
         - \alpha P_{S|C} - \hat{\alpha} P_{I|C}) P_{C} \\
P_{E} & = & 1 - P_{S} - P_{I} - P_{C}.
\end{eqnarray}
The interactions between the two species
are the same as in the main model; only the means of determining competitive
strength has changed.  Here, membership in the competitively weak class
of ``hosts'' comes at birth rather than through infection. 

\subsection*{Simulations}

	I performed all simulations on a $100 \times 100$ lattice using
asynchronous updating (Durrett 1995).  Figure 2 shows snapshots 
of the main and alternate models at equilibrium; the density of each type is
approximately equal in the two models.  Code for both models to run on
the simulation platform S3 (by T. Cox and R. Durrett) is available under 
the name DISCOM1 from http://gumby.syr.edu.

\section*{Pair Approximation}

	The master equations as given are exact; they describe the temporal 
evolution of the probability distribution for the state of a randomly
chosen site.  From the spatial ergodicity of the processes, these
quantities also give the spatial average of the density of each type
at any time during any single realization of the process.  However, each
system of equations is not closed; they depend on the quantities
$P_{\sigma|\sigma'}$ which must be determined.  There are two solutions
to this problem of closure.  The first approach is to use the approximation
$P_{\sigma|\sigma'} = P_{\sigma}$, known as the mean field approximation.
This removes spatial structure from the model; it can be thought of as
representing a situation in which the system is well--mixed by movement,
or in which dispersal and competition occur uniformly across arbitrarily
long distances (Durrett 1995).  The second approach is to try to compute 
$P_{\sigma|\sigma'}$ explicitly.  One can write down differential 
equations for
these conditional probabilities, but they incorporate terms involving
triplets of sites.  One then has to close these equations by approximating
the triplet probabilities in terms of pairs; this leads to the so--called 
pair approximation (Matsuda et al. 1992, Rand 1999).  The closed pair
equations give an approximation to the full system that includes local
spatial structure.  
In the appendix I present the differential equations for 
$P_{\sigma\sigma'}$, the probability that a random site is in state
$\sigma$ and a randomly chosen neighbor is in state $\sigma'$.
These equations describe the temporal evolution of the probability
distribution on the states of pairs of neighboring sites.  

	Three questions now arise: How do we close the pair equations, 
how do we use
them, and what aspects of the system do they capture?  Taking the last
question first, the equations only describe pairwise spatial structure.
As such, they do not capture larger scale spatial organization.  The
formation of large clusters of each type may be important
in determining the dynamics of the system, but the pair approximation
approach treats the system as well--mixed above the scale of pairs.
Nevertheless, the pair equations incorporate one aspect of spatial
structure and thus are an improvement over the mean field equations.
They estimate critical parameter values better than mean field
equations (Matsuda et al. 1992, Rand 1999) and can predict qualitative 
behaviors that depend on spatial structure (Sato et al.
1994, Rand 1999).  Most importantly, they allow us to understand how local
spatial structure determines the model's dynamics.

	The simplest use of the closed pair equations would be
to solve them numerically to determine critical parameter values and
equilibrium densities.  However, this approach sacrifices understanding
for calculation; it is difficult to relate outcomes to the form of
the equations or to infer the role that spatial structure plays. 
The alternate approach, which I use here, is to use the pair equations
to determine how spatial structure is incorporated into the master
equations.  This method, due to Matsuda et al. (1992), invokes a 
separation of
time scales argument.  Figure 3 shows time series for
the overall density $P_{C}$ and the conditional probability $P_{C|C}$
from a simulation.
The parameters are such that the competitor will eventually die out.
Notice that the conditional probability appears to achieve its equilibrium
much more quickly than the global density (the ``noisiness'' of the
conditional probability comes from the reduced sample size of this
quantity in a finite simulation).  There is a simple heuristic explanation
for this phenomenon: pairwise interactions should allow the system
to equilibrate more rapidly at the local than the global scale.
We can take advantage of this by computing the conditional probabilities,
then using them as fixed quantities in the master equations.
We can derive differential equations for the conditional 
probabilities from those for the pairs and singles:
\begin{equation}
\label{conditeqn}
\dot{P}_{\sigma|\sigma'} = \frac{d}{dt}\left(\frac{P_{\sigma\sigma'}}
	{P_{\sigma}}\right)
	= \frac{1}{P_{\sigma'}}\dot{P}_{\sigma\sigma'} - \frac{P_{\sigma|
	\sigma'}}{P_{\sigma'}}\dot{P}_{\sigma'}
\end{equation}
By computing equilibrium solutions of these conditional pair equations,
I obtain information about the spatial structure that I incorporate
as parameters in the master equations.  The spatial structure of the
populations is dynamic, but treating it as fixed in the master
equations allows efficient computation of threshold parameter values
for persistence.

	The final issue to resolve is how to close the pair (or, equivalently,
the conditional pair) equations.  The usual approach is to use the 
approximation $P_{\sigma|\sigma'\sigma''} = P_{\sigma|\sigma'}$.
Following Sato et al. (1994), I call this approach the ordinary pair 
approximation (OPA).
The assumption here is that ``the neighbor of my neighbor does not affect
my state."  In a variety of situations, this approximation is a good one
(Rand 1999).  However, the clustering of a rare type trying to invade 
introduces a problem.  As Figure 4 shows, the assumption 
$P_{\sigma|\sigma'\sigma} = P_{\sigma|\sigma'}$ breaks down when the
state $\sigma$ is rare.  Near the critical birthrate for successful
invasion, simulations show that $P_{C|IC} >> P_{C|I}$.  The few competitors 
are clustered together, so knowing that
another competitor is nearby is important.  Using the OPA assumption,
we would approximate $P_{C|IC} \approx P_{C|I}$, introducing a potentially
large error when the competitor is rare.
Since we are interested in determining the critical parameter values
for invasions, capturing this clustering of rare types may be important.

	One method of dealing with this was introduced by Sato et al. (1994), 
who called their approach improved pair approximation (IPA).  
The basic idea of IPA is to assume that $P_{\sigma|\sigma'\sigma''}
< P_{\sigma|\sigma'}$ in some cases, then to use conservation rules
for probabilities to compute $P_{\sigma|\sigma'\sigma}$.  
In their model of a sterilizing disease, Sato et al. (1994) assumed that 
$P_{S|EE}
= \epsilon P_{S|E}$ for some $\epsilon < 1$.  Other triplets were closed
using OPA.  It then follows that
$P_{S|ES} = 1 - P_{I|E} - \epsilon P_{E|E} = P_{S|E} + (1-\epsilon)P_{E|E}$.
Note that $P_{S|ES}/P_{S|E} \rightarrow \infty$ as S becomes rare,
since $P_{S|ES} \ne 0$ even when S is rare.  Also note
that OPA is a special case, obtained by setting $\epsilon = 1$.

	I follow Sato et al. (1994) in the general approach of IPA, although 
I use
slightly different clustering assumptions.  For each type of invasion
(host, disease, or competitor), I determined the nature of clustering
of the invader from simulations.  In general, I found two deviations
from the OPA closure assumptions: the segregation of the host and
competitor, and the clustering of infections within the host population.
Although the precise amount of clustering ($P_{\sigma|\sigma'\sigma}/
P_{\sigma|\sigma'}$) depends on the parameter values and states involved,
I used for simplicity the 
single parameter, $\epsilon = 0.8093$, that Sato et al. (1994) found
for the invasion of hosts into an empty system.  
Although some results may be sensitive to the value of $\epsilon$ 
(Sato et al. 1994), the most
important fact is that any $\epsilon < 1$ leads to increased clustering
of rare types.  As Figure 4 suggests, the value $\epsilon 
= 0.8093$ is conservative in many cases: here, it appears that
$P_{C|EI} \approx 0.4 P_{C|E}$ near the threshold. 

	Because the estimate of clustering used in IPA is valid for rare
types, this approach does not in general give good estimates of intermediate
equilibrium levels (far from the invasion threshold).  Thus, I used OPA
to determine the densities of resident types, and IPA to determine the
spatial structure of the invading competitor at low density.  Invoking the
separation of timescales argument, I included these quantities as fixed
values in the linearized master equation for the invader's density and 
computed the
critical parameter values for successful invasion.  The details of the 
closure schemes and the resultant IPA equations for the conditional 
probabilities are given the appendix.

\section*{Results}

	I present three types of results that describe the model's
dynamics and address the questions raised in the introduction.  First,
I show how equilibrium population levels depend on the growth rates of
the disease and two plant species.  Next, I use IPA to estimate the
threshold parameter values for disease and competitor invasions and
host extinction.  I compare these results to the alternate model and
explore how the results depend on local spatial structure.  Finally,
I examine the spatio--temporal dynamics of coexistence.

\subsection*{Growth Rates and Equilibria}

	Each plant species and the disease must have a sufficiently
high growth rate to invade the system.  A successful invasion by one
type comes at the expense of another type; for example, the invading
competitor displaces infected hosts.  Because of the cyclic structure
at the heart of the disease model, the impact of the invasion feeds
back on the invader.  
Figures 5, 6, and 7 illustrate
this phenomenon.
In these examples, the disease has just one effect: infected hosts can
be displaced, but their ability to displace competitors has not been
reduced.  The most interesting aspect of these results is that equilibirium
densities do not change monotonically with the growth rates.  

	As the host's birth rate increases, it first passes through
the critical value for the host to invade (Figure 5).  
The resulting invasion
suppresses the density of the competitor.  Further increase leads to
host densities sufficient to allow the disease to invade.  This invasion
in turn slows the growth of the host population and the decline of the
competitor.   Higher host birth rates lead to exclusion of the competitor; 
then the  equilibirium density of susceptible hosts remains constant, while
the density of infecteds grows with host birth rate.  

	For a fixed resident host density, the transmission rate
must be sufficiently large for the disease to invade (Figure 6).
Sufficiently high disease levels can allow the competitor to invade
where it previously could not.  This competitor invasion comes partly at
the expense of the host, so that increasing the transmission rate further
leads to decreases in the infected and overall host densities.  

	For a fixed host density with the disease endemic, a sufficiently
high birth rate allows the competitor to invade (Figure 7).
This invasion suppresses the disease level; 
an aggressive competitor can wipe out the disease, resulting
in a setback for the competitor.  This raises the interesting possibility
(discussed in the next section)
that after driving the disease out of the system, the competitor may itself
not be able to persist.  Once the disease disappears, the
competitor no longer has any effect on the host.

	These phenomena are due to the basic structure of the model,
which leads to a host density threshold for disease persistence as
well as negative feedback of the disease on itself via the competitor,
and vice versa.  As such, they do not depend on the spatial structure
of the model, or on the details of competition or disease transmission.
Notice that both the mean field and OPA methods capture the basic
effects of invasions on equilibrium population densities.
Still, spatial structure is important in determining the critical
parameter values for invasions.  In addition, the possibility of
host extinction induced by the disease depends critically on local
spatial structure (Sato et al. 1994).    
In order to gain a better understanding of how spatial structure
determines conditions for coexistence or extinction, we turn to maps of
parameter space generated using IPA.

\subsection*{Conditions for Coexistence}

	First I consider the question of how the disease and the
competitor affect each other's ability to persist.  
The interaction between infection and competition for a fixed
host birth rate is shown in Figure 8, with the thresholds
calculated using IPA.  When the disease has
the maximum possible effect on both directions of competitive 
displacement the thresholds for disease and competitor invasions
in  $\lambda,\beta_{2}$-space are shown in Figure 8a.  When
the disease is absent, the competitor can invade if $\beta_{2} \geq 3.1$.
Thus, the region between the lower curve and $\beta_{2} = 3.1$
describes conditions under which competitor invasion is only possible
because of the presence of the disease.  Conversely, the threshold for
disease invasion is increased from $\lambda = 4$ to the upper curve.
The presence of the competitor makes it much harder for the
disease to persist.  Because the disease threshold curve lies above that
for the competitor, we cannot have the situation alluded to above: 
an aggressive competitor which drives the disease from the system does
not thereby risk extinction.
(Of course, in a realization of the stochastic process in a finite population,
this outcome is possible.)  Notice that both thresholds depend more 
sensitively on the birth rate of the competitor than on the transmission 
rate of the disease.  The two curves should meet at the
threshold values indicated by the dashed lines; they miss because I
used different IPA closure schemes to compute the curves.  This illustrates
a drawback of the IPA technique: the various closure methods suggested
by simulations near different thresholds introduce inconsistencies ---
this is the price we pay to achieve more accurate overall threshold curves.

	Analogous threshold curves generated by mean field analysis (not 
shown) predict coexistence over a wider parameter set.  In particular, 
mean field theory predicts a much lower transmission rate needed for
disease invasion (Figure 6).  It is interesting to note
that mean field theory does a good job at estimating the threshold for
competitor invasion (Figure 7).  This occurs because spatial
structure has opposing effects.  One the one hand, it limits the ability of 
the competitor to reach empty sites and lowers the level of infection
in the host population.  On the other hand, it lowers the overall host
density and the host's ability to colonize sites occupied by the competitor.
The presence of opposing effects of spatial localization on a multitype
model, especially one with cyclic structure, suggests that spatial structure
should be treated explicitly to gain a realistic understanding of how
the dynamics depend on the parameters.
   
	When the host's competitive strength is set at birth (in the
alternate model) rather than through infection, there is no longer a
threshold population size or value of $\rho$ for the persistence of
competitively weak individuals in the population.  Thus, all phenomena
that depend on the presence of a disease threshold are absent from
the alternate model.  We can still ask how the presence of weak 
individuals in the host population affects the competitor's ability to
invade.  Figure 8b compares the competitor invasion thresholds
from the main and alternate models.  For the comparison, I adjusted
$\lambda$ and $\rho$ in the two models to achieve the same level of
``infection'' in the host population.  Note that it is more difficult
for the competitor to invade in the disease model than in the alternate
one.  Since the interactions between the host and competitor are the
same in the models, this difference is due to the spatial structure of
the populations.  As I discuss below, the disease leads to a segregation
of healthy and infected hosts; in the alternate model, weak individuals
are distributed randomly through the host population.  In turn, 
the spatial structure emerging during competitor invasion differs between
the models: the disease model develops higher $P_{S|C}$ and lower
$P_{I|C}$ than in the alternate model.  In other words, in the alternate
model invading competitors are more likely to encounter weak hosts
than they are in the disease model, even though the overall densities
are the same.

	So far I have assumed that the disease has the maximum possible
effect on both directions of competitive displacement.  It is interesting
to ask how the specific effects of the disease determine conditions
for coexistence.  Figure 8c shows the dependence of the 
competitor invasion threshold on the disease effects.  Note that the
growth rate parameters correspond to a point just above the competitor
invasion threshold in Figure 8a.  Thus, the competitor can
invade provided that a linear combination of the two disease effects
is sufficiently large.  Similarly, for growth rate parameters just above
the disease invasion threshold in Figure 8a, we see in 
Figure 8d that the disease can invade provided that a linear
combination of its effects is not too large.  The apparent linearity
of these thresholds in $\delta_{1},\delta_{2}$--space is an interesting
phenomenon.  In the mean field approximation, we see that the success
of the invasions depends linearly on these parameters (in fact, the 
ability of the disease to invade does not depend on $\delta_{1}$;
Equation 2).
However, in the pair approximation, the spatial structure emerging
during the invasion depends on $\delta_{1}$ and $\delta_{2}$; thus, 
these parameters are implicitly present in the master equations in 
nonlinear combinations.  It turns out that the spatial structure
of the competitor invasion ($P_{\sigma|C}$) is not strongly affected
by the specific disease effects.  The spatial structure of the disease
invasion ($P_{\sigma|I}$) is affected by the disease effects (hence
the nonzero slope of the threshold), but the resulting relationship
appears linear.  The causes of this very weak nonlinearity are not obvious,
but they might be clarified by more careful analysis of the pair equations.

\subsection*{Host Extinction}

	The next question to be addressed is whether a sublethal disease can
cause extinction of a host plant population.  Usually, we do not think of 
infectious diseases as a cause of host extinction, since the disease should
fade out of the population at low density.  However, Sato et al. (1994) used 
an IPA
analysis to argue that a sterilizing disease can cause host extinction
in a spatial model.  Because of the spatial clustering of the host, they
found that S--I contacts remained sufficiently high at low densities to
prevent any hosts from escaping infection.  This result depends critically
on the IPA closure (it can only occur for sufficiently small $\epsilon$),
and on the separation of time scales argument which bounds $P_{S|I}$ away
from zero even for arbitrarily low densities of S and I.

	In order to determine the threshold for host extinction, we need
to determine when the largest eigenvalue of the system of host equations
(1--2 or 5--6) is zero.  We proceed as before by first determining the
resident competitor density and spatial structure, then computing the
spatial structure of the host population at low density ($P_{\sigma|S}$ 
and $P_{\sigma|I}$) and including these as constants in the master
equations for S and I.  Figure 9a shows the result in the
case that disease only suppresses the ability of the host to displace
the competitor.  At high competitor densities, this is analogous to
the situation modeled by Sato et al. (1994) in which the disease suppresses 
host reproduction.  For sufficiently high transmission rates, IPA predicts
that host extinction will occur.  For intermediate transmission, full
coexistence occurs, and for sufficiently low transmission the disease,
the competitor, or both disappear.  Notice that the critical relative
transmission rate, $\lambda/\beta_{1}$, for host extinction first
decreases then increases with host birth rate.  This is in contrast to
the results in Sato et al. (1994), where the critical relative transmission 
rate
decreased monotonically with host birth rate.  That occured because
host density increased with birth rate, making it easier for the disease
to spread through the entire population.  However, in my model high
host birth rate threatens the competitor.  As the relative transmission
rate needed to save the competitor from exclusion increases, so must
the threshold for host extinction (since the disease only has an effect
in the presence of the competitor).

	Host extinction should be easier to
obtain with higher competitor birth rates; this is confirmed in 
Figure 9b.  In addition, the conditions for host extinction 
should depend on the specific disease effects.  
Figure 9c shows that host extinction occurs for 
intermediate combinations of the two disease effects.  For lower
disease effects, the damage done to the host population is not sufficient
to bring about extinction.  If the combination of disease effects
is too high, infected individuals are cleared from the population
more rapidly than they can transmit the disease at low densities, 
resulting in fadeout of the disease.  It is interesting to note that
the disease can cause extinction in this model by affecting either
direction of competitive displacement alone. In the alternate model,
host extinction is easier to obtain, because it does not rely on
the clustering of S and I individuals at low density.  If a sufficiently 
high proportion of ``host'' individuals are born into the weak class,
extinction will occur (Figure 9d).  As in the disease
model, the rate at which weak individuals are produced must be
sufficiently large to prevent the host from driving the competitor
extinct.

	Although IPA predicts that the disease can cause host extinction,
we must interpret this result with some care.  In simulations of the
model, host extinction is not a certainty in the region
of parameter space where IPA predicts it.  When the host population
reaches low density, the disease may fade out of the system, allowing
the surviving host population to rebound.  For a large enough system,
the probability of having an isolated host escape infection should
increase.  Thus, it is likely that disease-induced host extinction
is not possible in the interacting particle system on an infinite
lattice.  In a finite system, host extinction is possible, but it
depends on the particular distribution of S and I at low density
in any given realization of the process.  At some point, the separation
of time scales breaks down, because $P_{S|I}$ cannot remain constant at
arbitrarily low densities; then either host extinction or disease
fadeout may occur, depending on the realization.  Thus, a sensible
interpretation of the IPA result may be that it indicates parameter
values for which the disease can drive the host to very low densities,
at which point the host has a high risk for stochastic extinction events.

\subsection*{Spatiotemporal Dynamics of Coexistence}

	I have shown that the disease induces conditions for coexistence
different from those in the alternate model. It is also useful to
ask whether the spatiotemporal dynamics of coexistence, when it occurs,
differ between the two models.  Figures 10 and 11
show time series from simulations of the disease and alternate models
respectively, after transients have died out.  The parameters $\lambda$
and $\rho$ were chosen so that the time--averaged densities of the
two systems are approximately equal; all other parameters are the same.
Notice that the disease model displays large amplitude oscillations
that are absent from the alternate model.  This cycling behavior
is driven by the periodic outbreaks of the disease: the disease sweeps
through a patch of susceptible hosts, allowing the competitor to 
invade locally, after which healthy hosts move in.  The detection of
these oscillations depends on the spatial scale at which the system
is observed.  For an arbitrarily large system, oscillations far apart
will be independent and out of phase, yielding constant global densities
(Durrett 1999).
Presumably, observing the system at a carefully chosen spatial scale
that maximizes the deterministic signature (Keeling et al. 1997, Pascual and
Levin 1999) would give even clearer
evidence for cyclic dynamics.  In other models, limit cycle behavior
in the mean field (Durrett 1999) or pair (Rand 1999) equations has been 
used to explain
local cycling in the interacting particle system.  However, I have not
detected limit cycles in the mean field or pair equations for my model.
This suggests that the cycling is due to higher order spatial structure,
such as the buildup of large patches of susceptible hosts
subject to local epidemics.  These oscillations only occur in a subset of
the parameter values that allow coexistence.  I do not yet know how
the various parameters affect the cycling behavior, but it appears to
require high densities and rapid transmission. 

	As the snapshots of the two models in Figure 2 suggest,
the disease induces different spatial structure than would otherwise
occur.  Some aspects of this structure can be captured using spatial
covariances.  Let $P_{\sigma\sigma'}(r)$ be the probability that an
unordered pair of sites distance r apart are in the states $\sigma$
and $\sigma'$.  Then the spatial covariance of these states at distance
r is defined as $Cov_{\sigma\sigma'}(r) = P_{\sigma\sigma'}(r) - 
P_{\sigma}P_{\sigma'}$.  Positive covariance indicates aggregation
at that distance, while negative covariance indicates segregation.
Figure 12 shows covariances for the disease and alternate
models using the same parameter values as in Figures 2, 10, and 11.
Here, $r = 0$ denotes the four
nearest neighbors, and $r > 0$ denotes the sup norm (r = max$\{(x_{2}
- x_{1}),(y_{2} - y_{1})\}$).  In both models, each plant species
displays local aggregation, while the two species segregate from 
each other.  This segregation is noticeably stronger in the disease
model.  Because disease outbreaks shift the net competitive advantage
over entire regions, we tend to see the buildup
of larger patches of one species in the disease model than in the
alternate model where competitive advantage varies at a much finer
scale.  This suggests that segregation of competitors that is higher
than would otherwise be expected may indicate the presence of a
pathogen at work.  In addition, the disease model yields segregation
between susceptible and infected hosts, even though the host population
is aggregated.  Thus, if the competitive strength of individuals can
be estimated, segregation of weak and strong individuals within the
same population may also indicate the presence of a disease.  However,
this segregation may also arise from other causes such as dependence
of competitive strength on local environmental conditions. 

\section*{Discussion}

	Less is known about how infectious diseases structure natural
communities than other biotic processes such as predation and competition
(Dobson and Crawley 1994, Grenfell and Dobson 1995).
One reason for this is that many diseases have subtle effects on 
their hosts.  The sublethal damage caused by diseases can affect how
their hosts interact with their environment.  A classic example is
the increased risk of predation for diseased animals (Ives and Murray 1997).  
Here, I have
considered another important but often overlooked example: 
disease--induced competitive weakness in plants.  By studying a simple but
spatially explicit model of competing plant species, I have shown
that sublethal diseases can affect community structure and dynamics 
in novel ways.  The
competitive heterogeneity induced by sublethal diseases is distinguished
by the existence of a population threshold below which the disease
vanishes, and by the development of spatially localized, periodic
disease outbreaks.  Some features of the dynamics induced by sublethal
diseases may be shared by other biotic agents that affect competition
between plants.  Herbivores (Grace and Tilman 1990) and parasitic plants 
(Gibson and 
Watkinson 1986, Matthies 1996) can cause sublethal damage and 
change competitive interactions.  When these agents are host specific,
reproduce rapidly, and disperse locally, the conclusions of my
model may apply equally well to them. 

\subsection*{Disease--Competition Interactions in Space}

	The ability of the host, the competitor, and the disease to 
persist in my model depended on three growth rates and two disease
effects.  The presence of the disease allowed the competitor to
persist at lower birth rates than would otherwise have been possible.
Conversely, the competitor made it much more difficult for the
disease to invade the host population.  This effect would be increased
if the competitor displaced healthy hosts, since lower host density
would impede disease invasion.  The disease can put the host at
risk for extinction, provided it spreads quickly enough to persist
at low densities and does not have too severe consequences for infected
individuals.  The cyclic structure at the heart of the disease model
yielded a kind of law of decreasing returns for each growth rate;
aggressive growth by any of the types was punished in accordance with
the principle that ``the enemy of my enemy is my friend''.
These results should be robust, because they stem from the basic
assumption that the competitive relationship between plants can be
reversed by an infectious disease.  An alternate model in which 
competitive strength was determined at birth rather than by infection
did not display local cycling and lacked the feature of a critical
host density for disease persistence.  In the alternate model, 
host extinction was easier to obtain, and competitor invasions were
more likely to succeed for a given resident host population.  These
results illustrate that biotic agents such as infectious
diseases may shape plant community dynamics in ways that are 
qualitatively different from other sources of competitive heterogeneity.

	The spatial localization of both disease transmission and competition
was important in determining the quantitative and qualitative predictions
of the model.  I used pair approximations to compute the development of
local spatial structure on a fast time scale; this allowed me to compute
threshold parameter values efficiently and to explore how local spatial
structure influenced the dynamics. 
The structure of the
thresholds for disease and competitor invasions were the same in the
spatial and nonspatial (mean field) approximations.  However, the spatial
approximation
yielded higher growth rates needed for invasions, with coexistence
possible over a narrower range of parameter values. 
The possibility of host extinction in the disease model was intimately tied to
spatial structure, since it depended on the clustering of susceptible
and infected individuals at low density.  The disease induced 
different spatial structure from the alternate model of competitive
heterogeneity: segregation of the two plant species as well as of
weak and strong hosts was increased by the localized spread of the
disease.  This in turn led to the possibility of sustained oscillations
driven by periodic local epidemics.  This phenomenon was not captured
by studying pairwise interactions, so it appears to depend on higher
order spatial structure.

	From this model it is clear that sublethal diseases can cause
population
dynamics and spatial structure that differ qualitatively from other
causes of competitive heterogeneity.  Whether this means that explicit
understanding of the disease dynamics is necessary depends on the system 
being studied and the questions being asked.  In systems where disease
incidence is fairly uniform in space and time, the equilibrium structure
of the plant community may be understood by incorporating disease
effects into overall competitive heterogeneity.  
This may be the case in most unperturbed natural communities, where disease
incidence tends to be low and relatively constant (Dinoor and Eshed 1984).
On the other hand, 
when competitive dominance seems to shift between species in space or
time, infectious diseases are a possible explanation.  Large epidemics
are made more likely by human activities that introduce novel pathogens or 
decrease community diversity (Dobson and Crawley 1994). In addition, the 
existence
of a host population threshold for disease invasion or persistence has
important implications for how plant communities affected by sublethal
diseases will respond to 
perturbations.  When a host population passes through such a threshold,
the competitive hierarchy in the community can change.  In addition,
predictions about the success of plant invasions that come from controlled
experiments may give inaccurate results if the conditions or sample sizes
do not replicate disease incidence in the natural system.  Thus, while
sublethal diseases may be viewed as a minor source of heterogeneity in many
equilibrium natural communities, there are important situations in which
explicit understanding of disease dynamics may be critical.

\subsection*{Open Questions}

	Better understanding of the effects of sublethal diseases on 
community structure will require a combination of generic models like 
this one, system--specific detailed models, and experiments.  Simple
models can be used to answer a number of basic questions that I did
not address.  For example, different dispersal scales of the pathogen and
plants may be important (Thrall and Burdon 1999).  Long--range dispersal 
of the pathogen
may lead to spatially homogeneous disease incidence, changing the impact
of the disease on community structure.  Spatiotemporal variability in
environmental conditions may reinforce or confound the effects of the
disease.  Many pathogens, such as fungi, are sensitive to environmental
conditions like light or moisture (Burdon 1987). The destabilizing effects 
of diseases
may be increased if outbreaks are driven by environmental heterogeneity.
On the other hand, environmental variation can cause spatially and
temporally correlated changes in plant competitive strengths, which could
be confused with the effects of a disease.  These issues could be 
clarified by incorporating environmental heterogeneity into generic
models.  Another question raised by my model is how higher order spatial
structure differs between disease and non-disease models, and how this
structure influences dynamics.  There is a clear need for analytic
approximations for such models above the scale of pairs.  Finally, the
role of stochasticity in determining the success of invasions needs
to be examined more carefully.  Although my model is stochastic, the
assumption of an infinite population size removed many important effects
of stochasticity.  In an infinite population, ergodicity implies that
if, for example, $P_{I}(t) = .01$, then every realization of the
process will have this density of infected individuals at that time.
In a finite population, the density of infected individuals at time t will 
be a random variable with mean $.01$; in some realizations the
disease will be present at higher densities, and in some it will have
faded out entirely.  Thus, it is important to study finite stochastic
models to determine the probability distributions of
events that depend critically on low densities, such as host extinction
and disease fadeout (Keeling and Grenfell 1997, Keeling 2000).

	More realistic system--specific models will be useful for 
determining how the community effects of sublethal diseases depend
on the mechanisms of pathogen transmission and plant competition.
Further understanding of these systems will require the incorporation
of specific disease--induced physiological changes at particular
life stages; diseases that induce leaf drop in adults may have very
different consequences from those that slow seedling growth.
Ultimately, true understanding of the role of sublethal diseases
will have to come from experimental studies that show not only that
such diseases can impact competition in simple, controlled systems,
but that they actually contribute to the structure of real plant
communities.

\section*{Acknowledgments}

This research was conducted with support from NSF DBI-9602226, the 
Research Training Grant -- Nonlinear Dynamics in Biology, awarded to the
University of California, Davis.  I thank Alan Hastings for his
guidance and Carole Hom for her comments on the manuscript.

\section*{Appendix}

\subsection*{Pair Equations}

        Let $P_{\sigma\sigma'}$ be the probability that a randomly chosen
site is in state $\sigma$ and a randomly chosen neighbor of it is in state
$\sigma'$, so that $P_{\sigma\sigma'}=P_{\sigma|\sigma'}P_{\sigma'}$.
Note that we are treating pairs as ordered, but that
$P_{\sigma\sigma'}=P_{\sigma'\sigma}$. 
Also, let $P_{\sigma|\sigma'\sigma''}$ be the probability that, given a
random site is $\sigma'$ and a randomly chosen neighbor is $\sigma''$, another
randomly chosen neighbor is $\sigma$.  Then by considering the changes
that can occur to a pair, and the rates at which they occur (and remembering
that in continuous time only one change to a pair can occur at a time),
we obtain for the disease model:
\begin{eqnarray}
\dot{P}_{SS} & = & 2\{\frac{\beta_{1}}{4}(1+3P_{S|ES})P_{ES} + \frac
        {3\beta_{1}}{4}P_{I|ES}P_{ES} +  \frac{\alpha}{4}(1+3P_{S|CS})P_{CS}
        \nonumber \\
        & & + \frac{3\hat{\alpha}}{4}P_{I|CS}P_{CS}\} - 2\{\mu_{1}
        + \frac{3\lambda}{4}P_{I|SS}\}P_{SS} \\
\dot{P}_{SI} & = & \frac{3\beta_{1}}{4}P_{S|EI}P_{EI} + \frac{\beta_{1}}{4}
        (1+3P_{I|EI})P_{EI} + \frac{3\alpha}{4}P_{S|CI}P_{CI}
        \nonumber \\
        & & + \frac{\hat{\alpha}}{4}(1+3P_{I|CI})P_{CI} + \frac{3\lambda}{4}
        P_{I|SS}P_{SS}
        \nonumber \\
        & & - \{2\mu_{1} + \frac{\lambda}{4}(1+3P_{I|SI}) + \frac{3\gamma}{4}
        P_{C|IS}\}P_{SI} \\
\dot{P}_{SC} & = & \frac{3\beta_{1}}{4}P_{S|EC}P_{EC} + \frac{3\beta_{1}}{4}
        P_{I|EC}P_{EC} + \frac{3\alpha}{4}P_{S|CC}P_{CC}
        \nonumber \\
        & & \frac{3\hat{\alpha}}{4}P_{I|CC}P_{CC} + \frac{3\beta_{2}}{4}
        P_{C|ES}P{ES} + \frac{3\gamma}{4}P_{C|IS}P_{IS}
        \nonumber \\
        & & - \{\mu_{1} + \mu_{2} + \frac{3\lambda}{4}P_{I|SC} +
        \frac{\alpha}{4}(1+3P_{S|CS}) + \frac{3\hat{\alpha}}{4}P_{I|CS}\}
        P_{SC} \\
\dot{P}_{II} & = & 2\{\frac{\lambda}{4}(1+3P_{I|SI})P_{SI}\} -
        2\{\mu_{1} + \frac{3\gamma}{4}P_{C|II}\}P_{II} \\
\dot{P}_{IC} & = & \frac{3\lambda}{4}P_{I|SC}P_{SC} + \frac{3\beta_{2}}{4}
        P_{C|EI}P_{EI} + \frac{3\gamma}{4}P_{C|II}P_{II}
        \nonumber \\
        & & -\{\mu_{1} + \mu_{2} + \frac{\gamma}{4}(1+3P_{C|IC})
        + \frac{\hat{\alpha}}{4}(1+3P_{I|CI})
        \nonumber \\
        & & + \frac{3\alpha}{4}P_{S|CI}\}P_{IC} \\
\dot{P}_{CC} & = & 2\{\frac{\beta_{2}}{4}(1+3P_{C|EC})P_{EC}
        + \frac{\gamma}{4}(1+3P_{C|IC})P_{IC}\}
        \nonumber \\
        & & - \{\mu_{2} + \frac{3\alpha}{4}P_{S|CC}
        + \frac{3\hat{\alpha}}{4}P_{I|CC}\}P_{CC} \\
P_{E\sigma} & = & \frac{1}{2}\{P_{\sigma} - P_{\sigma\sigma} -
        2\sum_{\sigma' \neq \sigma,E} P_{\sigma\sigma'}\}.
\end{eqnarray}

\subsection*{Structure of IPA}

	The problem of closure is to
estimate $P_{\sigma|\sigma'\sigma''}$ in terms of $P_{\sigma|\sigma'}$,
when $\sigma$ is at low density.  In general, these terms are not equal.
Conservation rules for probabilities allow one to assume certain forms
for some of the triplets, then deduce the forms for the rest.  The key 
relationship is that:
\begin{eqnarray}
P_{\sigma\sigma'} & = & P_{\sigma\sigma'\sigma} + \sum_{\sigma''\ne\sigma}
                        P_{\sigma\sigma'\sigma''}
                \nonumber \\
                & = & P_{\sigma|\sigma'\sigma}P_{\sigma\sigma'} +
                \sum_{\sigma''\ne\sigma} P_{\sigma|\sigma'\sigma''}
                P_{\sigma'\sigma''}.
\end{eqnarray}
Dividing by $P_{\sigma'}$, we get:
\begin{equation}
P_{\sigma|\sigma'} = P_{\sigma|\sigma'\sigma}P_{\sigma|\sigma'} +
                \sum_{\sigma''\ne\sigma} P_{\sigma|\sigma'\sigma''}
                P_{\sigma''|\sigma'},
                \nonumber
\end{equation}
 so that:
\begin{equation}
\label{conseqn}
P_{\sigma|\sigma'\sigma} = 1 - \frac{1}{P_{\sigma|\sigma'}}
                \sum_{\sigma''\ne\sigma} P_{\sigma|\sigma'\sigma''}
                P_{\sigma''|\sigma'}.
\end{equation}

	As explained below, I assume $P_{\sigma|\sigma'\sigma''} = \epsilon
P_{\sigma|\sigma'}$ in certain cases.  
In order to estimate the value of $\epsilon$, Sato et al. (1994) used an
independent estimate of the critical birth rate for the disease-free process.
Solving the pair equations for the critical birth rate, they found:
\begin{equation}
\beta_{c} = \frac{4}{3\epsilon}.
\end{equation}
An independent estimate of $\beta_{c}/4=0.4119$ due to Katori and Konno (1990)
leads to
the estimate $\epsilon=0.8093$.  Thus, this value of $\epsilon$ captures
the clustering of occupied sites invading the empty system near the
critical birth rate for survival.

	As an example of my procedure for computing the
threshold parameter values for invasion, consider the case of competitor
invasion.  First, I compute the host densities using:
\begin{eqnarray}
\dot{P}_{S} & = & 0 \nonumber \\
\dot{P}_{I} & = & 0 \nonumber \\
\dot{P}_{SS} & = & 0 \nonumber \\
\dot{P}_{SI} & = & 0 \nonumber \\
\dot{P}_{II} & = & 0 \nonumber \\
E & = & 1 - S - I, \nonumber
\end{eqnarray}
dropping all terms involving C.
Then I determine the threshold for competitor invasion by solving:
\begin{eqnarray}
\frac{\dot{P}_{C}}{P_{C}} & = & 0 \nonumber \\
\dot{P}_{C|C} & = & 0 \nonumber \\
\dot{P}_{S|C} & = & 0 \nonumber \\
\dot{P}_{I|C} & = & 0, \nonumber
\end{eqnarray}
with $P_{C} = P_{C|S} = P_{C|I} = P_{C|E} = 0$.  
When all parameters are specified, we can solve the last three equations
for the spatial structure early in the invasion.  By the separation of
timescales argument, we can include this spatial structure as fixed
parameters in the master equations.  For threshold calculations, we can
specify all parameter values but one, then solve the four equations
simultaneously for the critical value of the parameter of interest.

        Next I present the specific closure assumptions I used and derive
the resulting differential equations for the conditional probabilities
involved with each type of invasion. 

\subsection*{Competitor Invasion}

	I found that $P_{C|\sigma I} = \epsilon P_{C|\sigma}$ for $\sigma
\ne C$. That is, the presence of an I individual nearby lowers the 
probability of finding a C.  This makes sense, since the mutual displacement
between these types will lead to segregation.  After closing other triplet
terms using OPA, we find from equation~\ref{conseqn} that:
\begin{equation}
P_{C|\sigma C} = P_{C|\sigma} + (1-\epsilon) P_{I|\sigma}. 
\end{equation}
Furthermore, from the relation $\sum_{\sigma'} P_{\sigma'|\sigma C} = 1$,
it follows that:
\begin{equation}
P_{I|\sigma C} = \epsilon P_{I|\sigma}
\end{equation}
for $\sigma \ne C$.

	Now, we use equation~\ref{conditeqn} with these closure assumptions
to derive equations for the conditional probabilities.  Since we are
interested in the behavior at low competitor density, we use $P_
{C|\sigma} = 0$ for $\sigma \ne C$ to get:
\begin{eqnarray}
\dot{P}_{S|C} & = & \frac{3}{4}\beta_{1} P_{S|E}P_{E|C} + \frac{3}{4}
		\beta_{1}\epsilon P_{I|E}P_{E|C} + \frac{3}{4}\alpha
		P_{S|C}P_{C|C} \nonumber \\
              & & + \frac{3}{4}\hat{\alpha} P_{I|C}P_{C|C}
	        + \frac{3}{4}\beta_{2}P_{E|C}P_{S|E} + \frac{3}{4}
		\gamma P_{I|C}P_{S|I} - \mu_{1}P_{S|C} \nonumber \\ 
	      & & - \frac{3}{4}\lambda\epsilon P_{I|S}P_{S|C} 
		- \frac{\alpha}{4}P_{S|C} - \beta_{2}P_{E|C}P_{S|C}
		\nonumber \\
	      & & - \gamma P_{I|C}P_{S|C} + \frac{\alpha}{4}P_{S|C}^{2}
		+ \frac{\hat{\alpha}}{4}P_{I|C}P_{S|C}  \\
\dot{P}_{I|C} & = & \frac{3}{4}\lambda\epsilon P_{I|S}P_{S|C} +
	      \frac{3}{4}\beta_{2}\epsilon P_{E|C}P_{I|E} +
	      \frac{3}{4}\gamma\epsilon P_{I|C}P_{I|I} -
 	      \mu_{1} P_{I|C} \nonumber \\
	      & & - \frac{\gamma}{4} P_{I|C} - \frac{3}{4} \gamma (1 - 
		\epsilon) P_{I|I} P_{I|C} - \frac{\hat{\alpha}}{4} P_{I|C}
		- \beta_{2} P_{E|C} P_{I|C} \nonumber \\
	      & & - \gamma P_{I|C}^{2} + \frac{\alpha}{4} P_{S|C} P_{I|C}
		+ \frac{\hat{\alpha}}{4} P_{I|C}^{2}  \\
\dot{P}_{C|C} & = & \frac{\beta_{2}}{2}P_{E|C} + \frac{3}{2}\beta_{2}
		(1-\epsilon)P_{I|E}P_{E|C} + \frac{\gamma}{2}P_{I|C}
		\nonumber \\
	      & & + \frac{3}{2}\gamma(1-\epsilon)P_{I|I}P_{I|C}
		- \mu_{2}P_{C|C} - \frac{\alpha}{2}P_{S|C}P_{C|C}
		\nonumber \\
              & & - \frac{\hat{\alpha}}{2}P_{I|C}P_{C|C}
	          - \beta_{2}P_{E|C}P_{C|C} - \gamma P_{I|C}P_{C|C}.
\end{eqnarray}
I used the same closure approximations to derive analogous equations for
the alternate model.

\subsection*{Disease Invasion}

For disease invasion, the predominant feature was the segregation of
susceptible and infected hosts.  Thus, I used the closure $P_{I|\sigma S}
= \epsilon P_{I|\sigma}$ for $\sigma \ne I$.  Using the conservation
laws as before, we find that:
\begin{equation}
P_{I|\sigma I} = P_{I|\sigma} + (1-\epsilon) P_{S|\sigma},
\end{equation}
and
\begin{equation}
P_{S|\sigma I} = \epsilon P_{S|\sigma}
\end{equation}
for $\sigma \ne I$.

	This closure scheme, along with the low density assumption that
$P_{I|\sigma} = 0$ for $\sigma \ne I$ leads to the equations:
\begin{eqnarray}
\dot{P}_{S|I} & = & \frac{\beta_{1}}{4}P_{E|I} + \frac{3}{4}\beta_{1}
		P_{S|E}P_{E|I} + \frac{3}{4}\alpha\epsilon P_{S|C}P_{C|I}
	      	\nonumber \\
	      & & + \frac{\hat{\alpha}}{4}P_{C|I} + \frac{3}{4}\hat{\alpha}
		(1-\epsilon)P_{S|C}P_{C|I} - \mu_{1}P_{S|I} 
		- \frac{\lambda}{4}P_{S|I} \nonumber \\
	      & & - \frac{3}{4}\lambda(1 - 2\epsilon)P_{S|S}P_{S|I}
		- \lambda P_{S|I}^{2} + \gamma P_{C|I}P_{S|I} \\
\dot{P}_{I|I} & = & \frac{\lambda}{2}P_{S|I} + \frac{3}{2}\lambda
		(1-\epsilon)P_{S|S}P_{S|I} - \mu_{1}P_{I|I} \nonumber \\
	      & & - \frac{\gamma}{2}P_{C|I}P_{I|I} - \lambda P_{S|I}P_{I|I} \\
\dot{P}_{C|I} & = & \frac{3}{4}\lambda P_{S|I}P_{C|S} + \frac{3}{4}\beta_{2}
		P_{C|E}P_{E|I} - \mu_{2}P_{C|I} \nonumber \\
	      & & - \frac{\gamma}{4}P_{C|I} - \frac{\hat{\alpha}}{4}P_{C|I}
		- \frac{3}{4}\hat{\alpha}(1-\epsilon)P_{S|C}P_{C|I} 
		\nonumber \\
	      & & - \frac{3}{4}\alpha\epsilon P_{S|C}P_{C|I}
		- \lambda P_{S|I}P_{C|I} + \gamma P_{C|I}^{2}.  
\end{eqnarray}

\subsection*{Host Extinction}

	Here, both S and I are at low density.  As with competitor invasion,
I found that the host and competitor segregate.  In addition, infected
individuals are clumped within the host population.  Thus, I used the
closure approximations $P_{S|\sigma C} = \epsilon P_{S|\sigma}$ for
$\sigma \ne S$, as well as $P_{I|CC} = \epsilon P_{I|C}$, $P_{I|EC} = 
\epsilon P_{I|E}$, and $P_{I|SS} = \epsilon P_{I|S}$.  Then it follows that:
\begin{equation}
P_{S|\sigma S} = P_{S|\sigma} + (1-\epsilon)P_{C|\sigma}
\end{equation}
and
\begin{equation}
P_{C|\sigma S} = \epsilon P_{C|\sigma}
\end{equation}
for $\sigma \ne S$.  Also, we have:
\begin{eqnarray}
P_{I|EI} & = & P_{I|E} + (1-\epsilon)P_{C|E} \\
P_{I|CI} & = & P_{I|C} + (1-\epsilon)P_{C|C} \\
P_{I|SI} & = & P_{I|S} + (1-\epsilon)P_{S|S} \\
P_{C|EI} & = & \epsilon P_{C|E} \\
P_{C|CI} & = & \epsilon P_{C|C} \\
P_{S|SI} & = & \epsilon P_{S|S}.
\end{eqnarray}

	This closure scheme, along with the low density assumption that
$P_{S|\sigma} = P_{I|\sigma} = 0$ for $\sigma \ne S,I$ leads to the
equations:
\begin{eqnarray}
\dot{P}_{S|S} & = & \frac{\beta_{1}}{2}P_{E|S} + \frac{3}{2}\beta_{1}
		(1-\epsilon)P_{C|E}P_{E|S} + \frac{\alpha}{2}P_{C|S}
		+ \frac{3}{2}\alpha (1-\epsilon)P_{C|C}P_{C|S}
		\nonumber \\
	      & & - \mu_{1}P_{S|S} - \frac{3}{2}\lambda\epsilon
		P_{I|S}P_{S|S} - \beta_{1}P_{E|S}P_{S|S} - \alpha
		P_{C|S}P_{S|S} \nonumber \\
	      & & + \lambda P_{I|S}P_{S|S} - \beta_{1}P_{E|I}P_{S|S} 
		\frac{P_{I|S}}{P_{S|I}} - \hat{\alpha}P_{C|I}P_{S|S}
		\frac{P_{I|S}}{P_{S|I}} \\
\dot{P}_{I|S} & = & \frac{P_{I|S}}{P_{S|I}} \{ \frac{\beta_{1}}{4}P_{E|I}
		+ \frac{3}{4}\beta_{1}(1-\epsilon)P_{C|E}P_{E|I}
		+ \frac{\hat{\alpha}}{4}P_{C|I}  \nonumber \\
	      & & + \frac{3}{4}\hat{\alpha}(1-\epsilon)P_{C|C}P_{C|I}
		 - \mu_{1}P_{S|I} - \frac{\lambda}{4}P_{S|I} 
		 - \frac{3}{4}\lambda (1-2\epsilon)P_{S|S}P_{S|I} 
		\nonumber \\
	      & & - \frac{3}{4}\gamma\epsilon P_{C|I}P_{S|I}
		 - \beta_{1}P_{E|S}P_{S|I} - \alpha P_{C|S}P_{S|I}
		 + \frac{\lambda}{4}P_{I|S}P_{S|I} \nonumber \\
	      & & - \beta_{1}P_{E|I}P_{I|S} - \hat{\alpha}P_{C|I}P_{I|S} 
		\} \\
\dot{P}_{C|S} & = & \frac{3}{4}\beta_{1}\epsilon P_{E|S}P_{C|E}
		+ \frac{3}{4}\beta_{1}\epsilon P_{E|I}P_{C|E}
		\frac{P_{I|S}}{P_{S|I}} + \frac{3}{4}\hat{\alpha}
		\epsilon P_{C|I}P_{C|C}\frac{P_{I|S}}{P_{S|I}} \nonumber \\
	      & & \frac{3}{4}\beta_{2}\epsilon P_{C|E}P_{E|S}
		+ \frac{3}{4}\gamma\epsilon P_{C|I}P_{I|S} - \mu_{2}P_{C|S}
		- \frac{\alpha}{4}P_{C|S} \nonumber \\
	      & & - \frac{3}{4}\alpha(1-2\epsilon)P_{C|C}P_{C|S}
		- \beta_{1}P_{E|S}P_{C|S} - \alpha P_{C|S}^{2}
		+ \frac{\lambda}{4}P_{I|S}P_{C|S} \nonumber \\
	      & & - \beta_{1}P_{E|I}P_{C|S}\frac{P_{I|S}}{P_{S|I}}
		- \hat{\alpha}P_{C|I}P_{C|S}\frac{P_{I|S}}{P_{S|I}} \\ 
\dot{P}_{S|I} & = & \frac{\beta_{1}}{4}P_{E|I} + \frac{3}{4}\beta_{1}
		(1-\epsilon) P_{C|E}P_{E|I} + \frac{\hat{\alpha}}{4}P_{C|I}
		\nonumber \\ 
	      & & + \frac{3}{4}\hat{\alpha}(1-\epsilon)P_{C|C}P_{C|I}
		- \mu_{1}P_{S|I} - \frac{\lambda}{4}P_{S|I}
		- \lambda P_{S|I}^{2} \nonumber \\
	      & & + \frac{3}{4}\gamma P_{C|I}P_{S|I} - \frac{3}{4}\lambda
		P_{I|S}P_{S|I} - \frac{3}{4}\lambda (1-2\epsilon)
		P_{S|S}P_{S|I} \\
\dot{P}_{C|I} & = & \frac{3}{4}\lambda P_{S|I}P_{C|S} + \frac{3}{4}\beta_{2}
		\epsilon P_{C|E}P_{E|I} + \frac{3}{4}\gamma P_{C|I}P_{I|I}
		- \mu_{2}P_{C|I} \nonumber \\
	      & & - \frac{\gamma}{4}P_{C|I} - \frac{\hat{\alpha}}{4}P_{C|I}
		- \frac{3}{4}\hat{\alpha}(1-\epsilon)P_{C|C}P_{C|I}
		\nonumber \\
	      & & - \lambda P_{S|I}P_{C|I} + \frac{\gamma}{4}P_{C|E}^{2} \\
\dot{P}_{I|I} & = & \frac{\lambda}{2}P_{S|I} + \frac{3}{2}\lambda P_{I|S}
		P_{S|I} + \frac{3}{2}\lambda(1-\epsilon)P_{S|S}P_{S|I}
		- \mu_{1}P_{I|I} \nonumber \\
	      & & - \frac{\gamma}{2}P_{C|I}P_{I|I} - \lambda P_{S|I}P_{I|I}.
\end{eqnarray}

	To derive the analogous equations for host extinction in the alternate
model, different closure approximations are warranted.  Here, infected
individuals are not clustered within the host population.  The only significant
deviation from OPA in simulations was the segregation of the host and 
competitor species.  Letting $H = S + I$ denote all hosts, we have
$P_{H|\sigma C} = \epsilon P_{H|\sigma}$ for $\sigma = C,E$.  Then we find 
that:
\begin{equation}
P_{S|\sigma H} = P_{S|\sigma} + (1-\epsilon) P_{C|\sigma}
	\frac{P_{\sigma|S}P_{S|I}}{P_{\sigma|S}P_{S|I} + P_{\sigma|I}P_{I|S}}
\end{equation}
\begin{equation}
P_{I|\sigma H} = P_{I|\sigma} + (1-\epsilon) P_{C|\sigma}
	\frac{P_{\sigma|I}P_{I|S}}{P_{\sigma|I}P_{I|S}+P_{\sigma|S}P_{S|I}}
\end{equation}
\begin{equation}
P_{C|\sigma H} = \epsilon P_{C|\sigma}
\end{equation}
for $\sigma = C,E$.

\pagebreak

\section*{Figure Captions}

Figure 1: State transitions in the disease (a) and alternate (b)
                models. Arrows indicate possible changes in state
                at a site.
\newline
\newline
Figure 2: Snapshots of the disease (a) and alternate (b) models at
$t = 1000$.  States are: S = dark gray, I = light gray, C = white, E = black.
Parameters are $\beta_{1}=\alpha=4, \beta_{2}=6, \lambda=20$ or $\rho=.54,
\delta_{1}=\delta_{2}=1$.
\newline
\newline
Figure 3: Separation of timescales between global and conditional densities.
                Initial conditions are $P_{S}=0.2$,
                $P_{I}=0.4$, $P_{C}=0.05$.  Parameter values are
                $\beta_{1}=\alpha=2$, $\beta_{2}=1.75$, $\lambda=10$,
                $\delta_{1}=\delta_{2}=1$.
\newline
\newline
Figure 4: Ratios of conditional probabilities near competitor
                invasion threshold. (a) Near the threshold birthrate
                ($\beta_{2} \approx 1.95$),
                competitors are rare and clustered.  A neighbor of an I
                is much more likely to be C if another neighbor is C than
                otherwise.  (b--d) A neighbor of a site is less likely
                to be C if another neighbor is I than otherwise.
                IPA captures this type of clustering, but OPA does not.
                Parameters are $\beta_{1} = \alpha = 2,
                \lambda = 10, \delta_{1} = \delta_{2} = 1$; results from
                simulation are at $t=10$.
\newline
\newline
Figure 5: Equilibrium densities (at $t = 500$) as a function of host
        birthrate. Parameters are $\alpha=\beta_{1}, \beta_{2}=5, \lambda=10,
        \delta_{1}=0, \delta_{2}=1.$
\newline
\newline
Figure 6: Equilibrium densities (at $t = 500$) as a function of disease
        transmission rate. Parameters are $\beta_{1}=\alpha=2, \beta_{2}=5,
        \delta_{1}=0, \delta_{2}=1.$
\newline
\newline
Figure 7: Equilibrium densities (at $t = 500$) as a function of competitor
        birthrate. Parameters are $\beta_{1}=\alpha=2, \lambda=10,
        \delta_{1}=0, \delta_{2}=1.$
\newline
\newline
Figure 8: a) Thresholds for disease and competitor invasions.
                Disease invades below top curve; competitor invades
                above bottom curve.  Dashed lines indicate threshold
                values for disease invasion without competitor and
                for competitor invasion without disease.  Parameters
                are $\beta_{1}=\alpha=2, \delta_{1}=\delta_{2}=1.$
                b) Comparison of competitor invasion thresholds in
                main and alternate models.  Competitor invades above
                curve.
                Parameters are $\beta_{1}=\alpha=4,
                \delta_{1}=\delta_{2}=1.$ We match host infection level
                by varying $\lambda$ and $\rho.$
                c) Dependence of competitor invasion on disease effects.
                Competitor invades above line.  Parameters are
                $\beta_{1}=\alpha=2, \beta_{2}=2.5, \lambda=10.$
                d) Dependence of disease invasion on disease effects.
                Disease invades below line.  Parameters are
                $\beta_{1}=\alpha=2, \beta_{2}=5.5, \lambda=10.$
\newline
\newline
Figure 9: a) Host extinction.  Host goes extinct above solid
                line.  Disease invades above dotted line.  Competitor
                invades above dashed line.  Parameters are
                $\beta_{2}=10, \delta_{1}=1, \delta_{2}=0.$
                b) Effect of competitor birth rate on host extinction.
                Host goes extinct above curves.  Other parameters as in (a).
                c) Dependence of host extinction on disease effects.  Host
                goes extinct in region between curves.  Parameters are
                $\beta_{1}=\alpha=5, \beta_{2}=10, \lambda/\beta_{1}=30$.
                d) Host extinction in alternate model.  Host goes extinct
                above solid line.  Competitor invades above dashed line.
                Parameters are as in (a).
\newline
\newline
Figure 10: Times series from a simulation of the disease model
                on a $100\times 100$
                lattice.  Parameters are $\beta_{1}=\alpha=4,
                \beta_{2}=6, \lambda=20, \delta_{1}=\delta_{2}=1$.
\newline
\newline
Figure 11: Times series from a simulation of the alternate model
                on a $100\times 100$
                lattice.  Parameters are $\beta_{1}=\alpha=4,
                \beta_{2}=6, \rho=.54, \delta_{1}=\delta_{2}=1$.
\newline
\newline
Figure 12: Covariances as a function of distance, from simulation of models
                on a $100\times 100$ lattice.
                Parameters are as in Figures 2, 10 and 11.

\pagebreak

\section*{Figures}

Figure 1:

\begin{figure}[h]

        \includegraphics{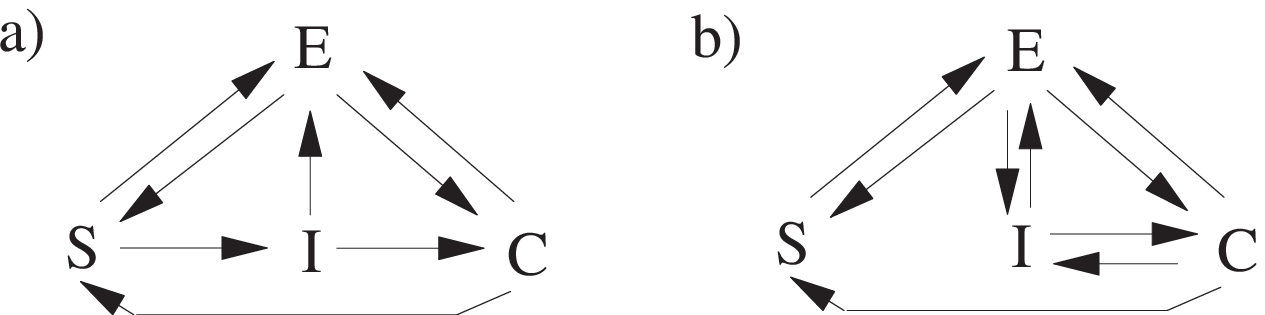}

\end{figure}

\pagebreak

Figure 2a:

\begin{figure}[h]

        \scalebox{.5}{\includegraphics{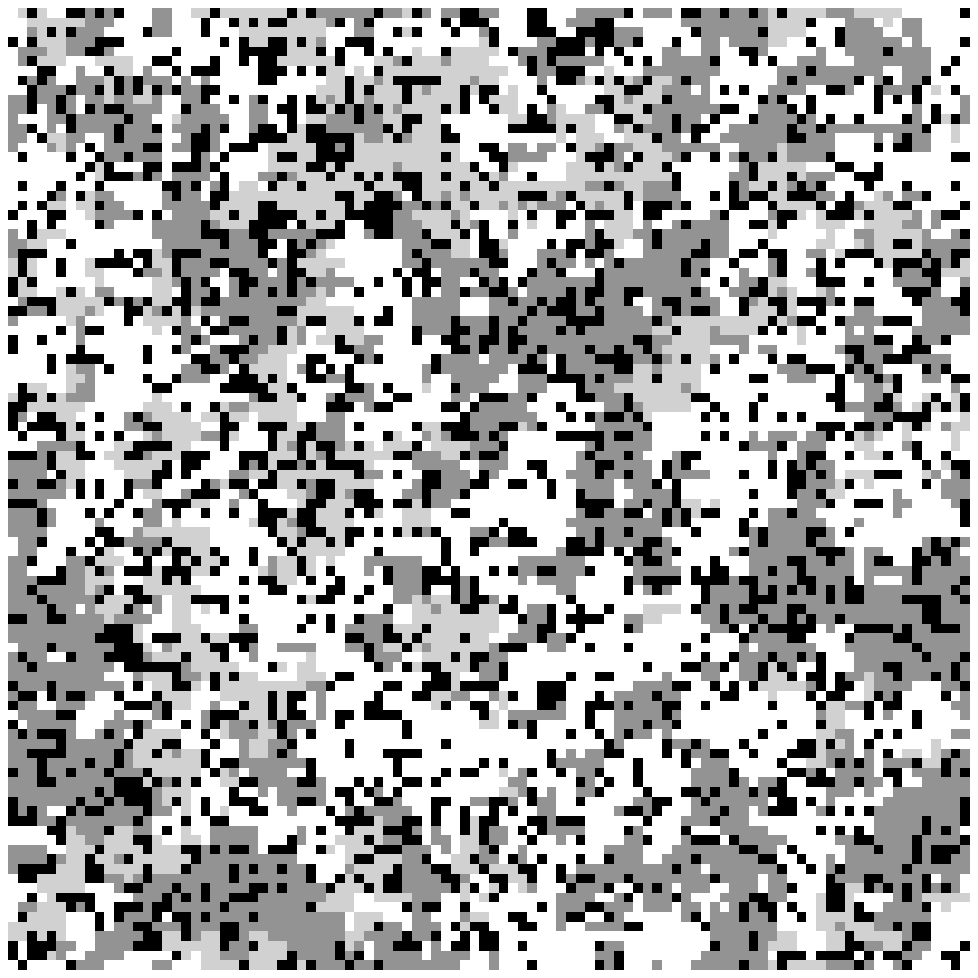}}

\end{figure}

Figure 2b:

\begin{figure}[h]

        \scalebox{.5}{\includegraphics{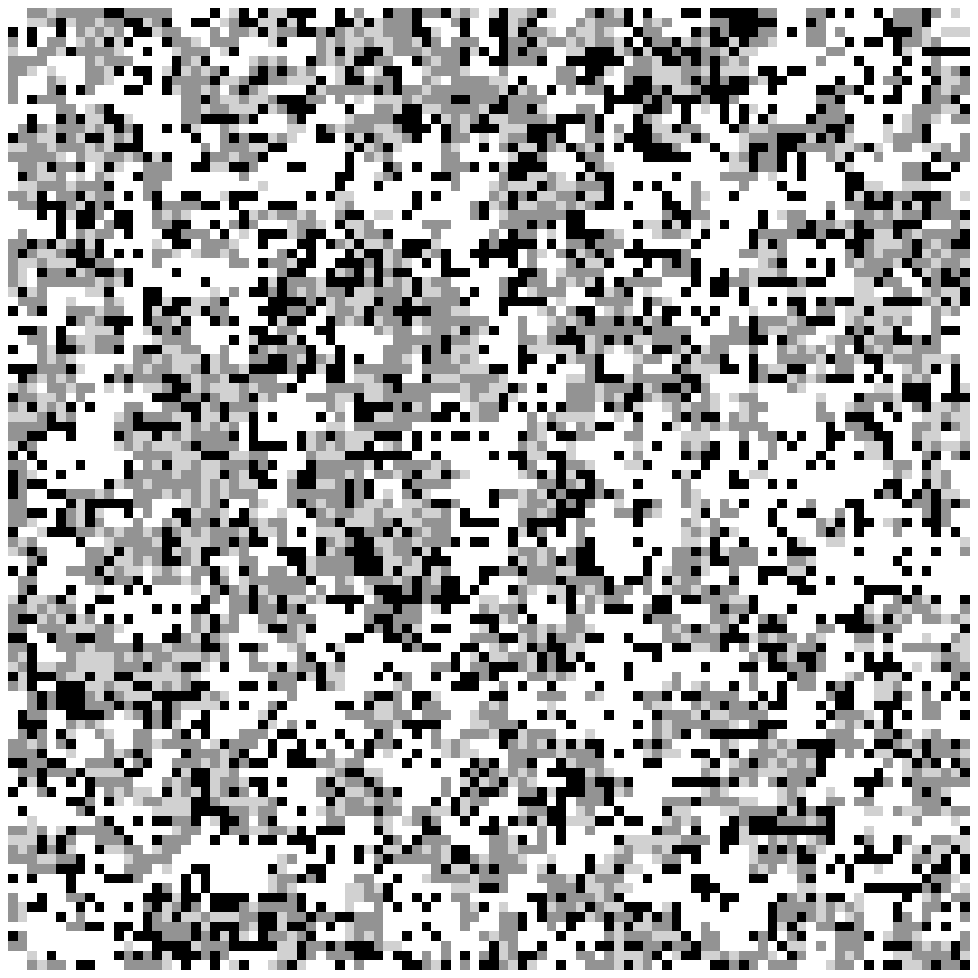}}

\end{figure}

\pagebreak

Figure 3:

\begin{figure}[h]

        \includegraphics{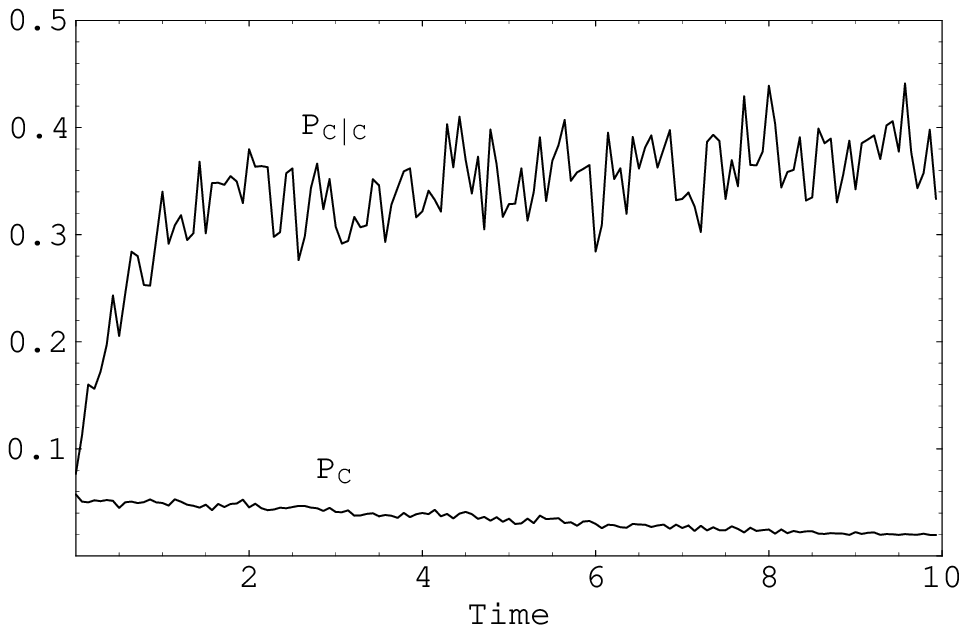}

\end{figure}

\pagebreak

Figure 4a:

\begin{figure}[h]

        \includegraphics{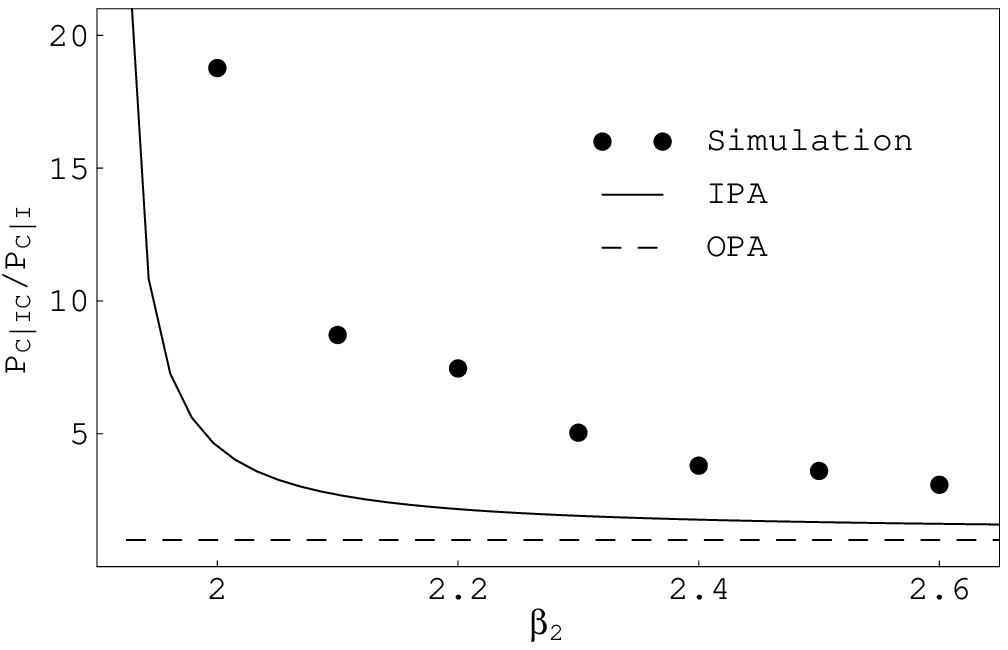}

\end{figure}

Figure 4b:

\begin{figure}[h]

        \includegraphics{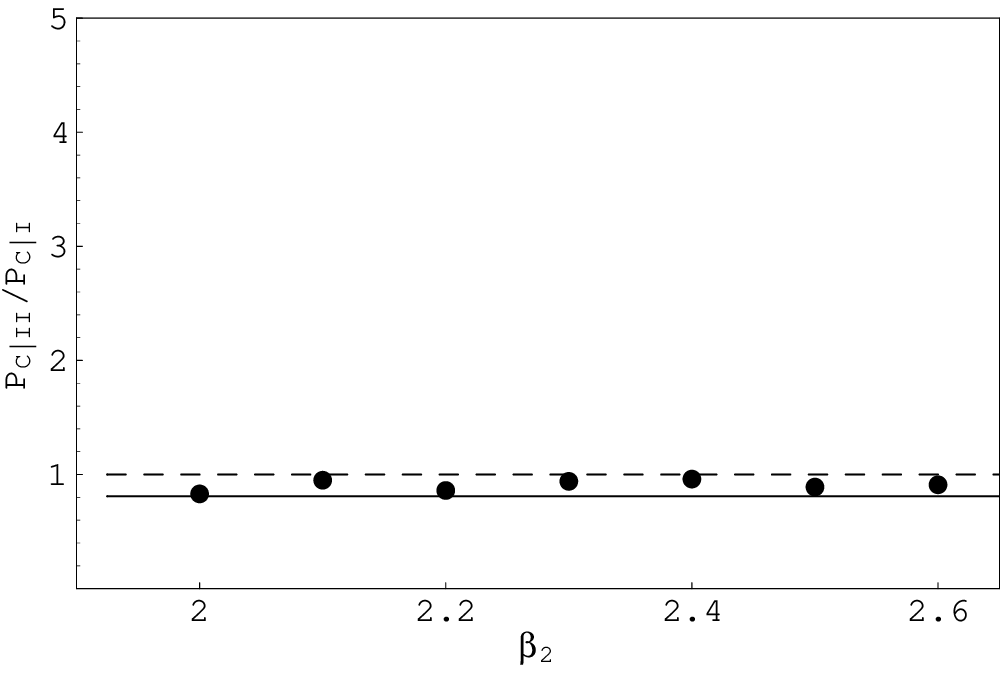}

\end{figure}

\pagebreak

Figure 4c:

\begin{figure}[h]

        \includegraphics{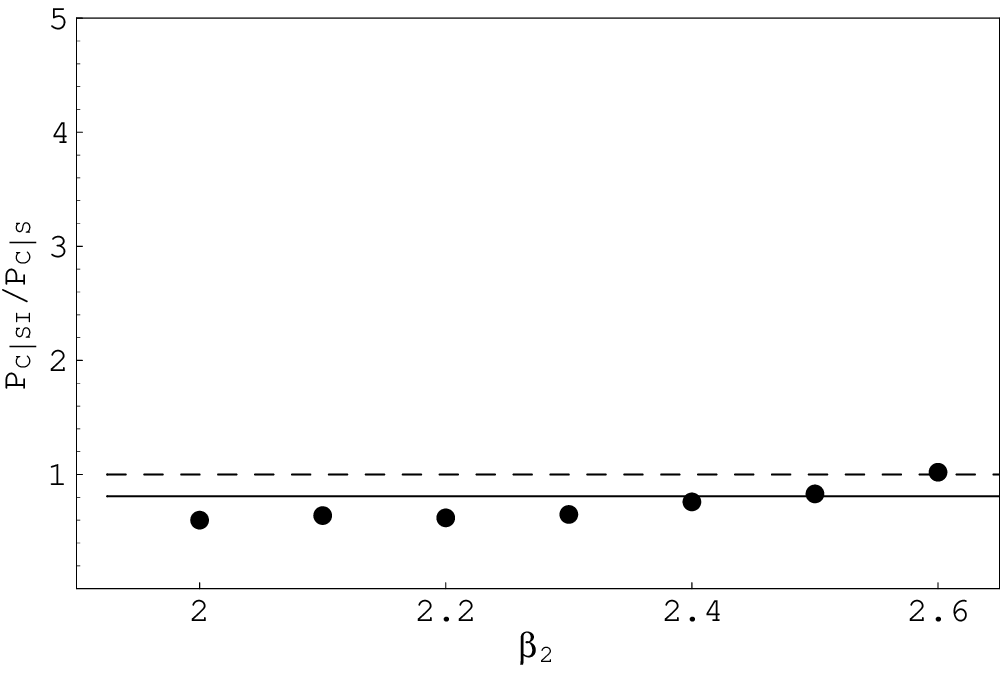}

\end{figure}

Figure 4d:

\begin{figure}[h]

        \includegraphics{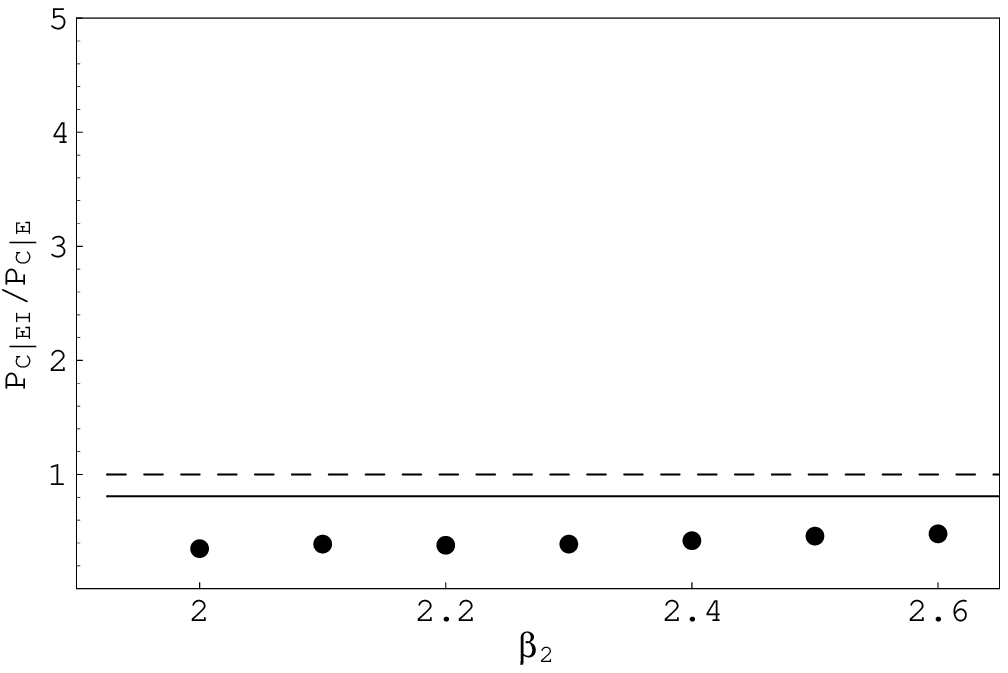}

\end{figure}

\pagebreak

Figure 5a:

\begin{figure}[h]

        \includegraphics{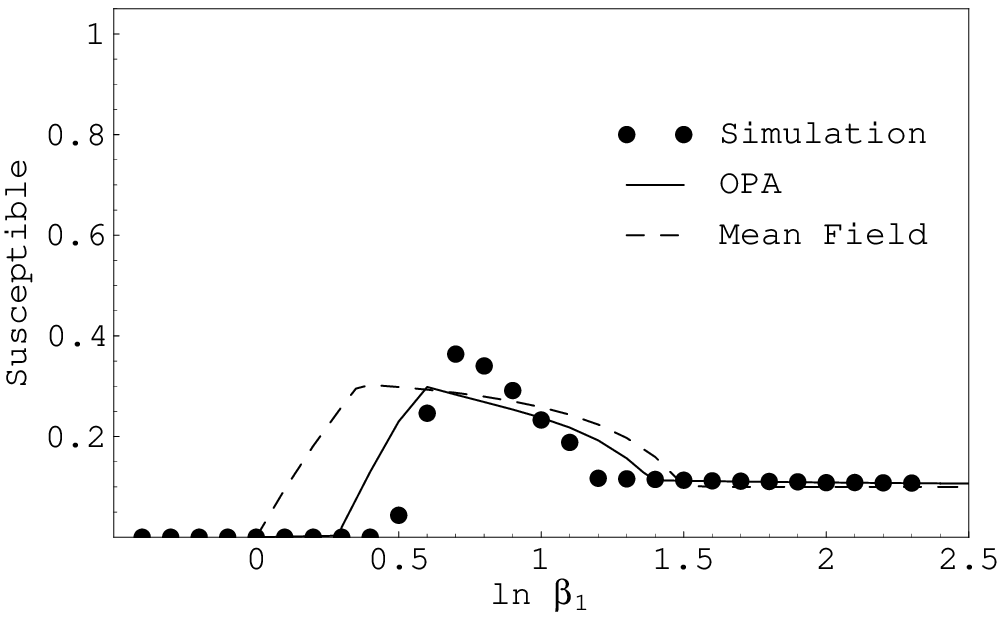}

\end{figure}

Figure 5b:

\begin{figure}[h]

        \includegraphics{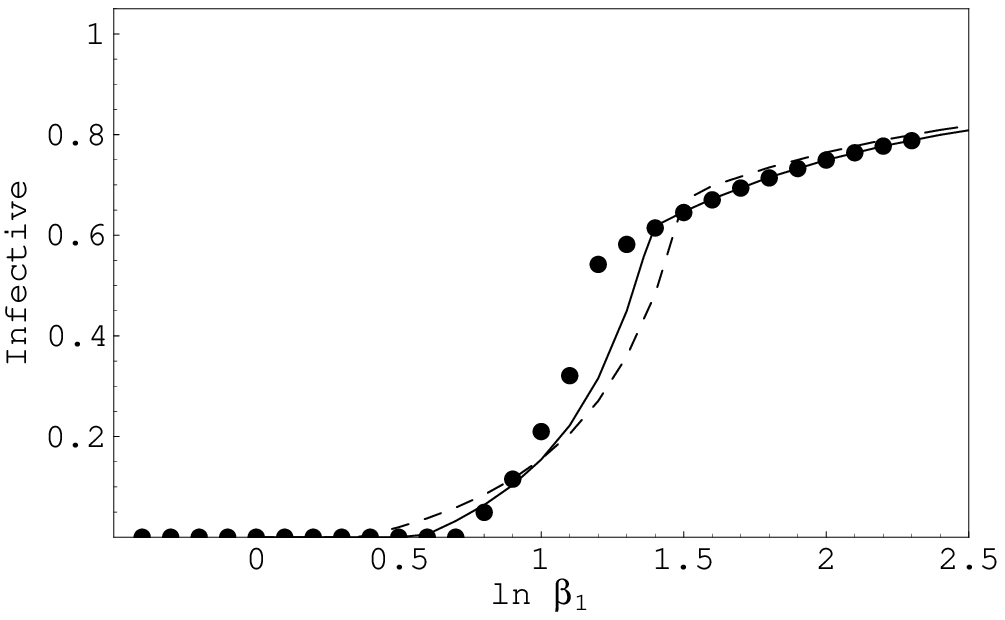}

\end{figure}

\pagebreak

Figure 5c:

\begin{figure}[h]

        \includegraphics{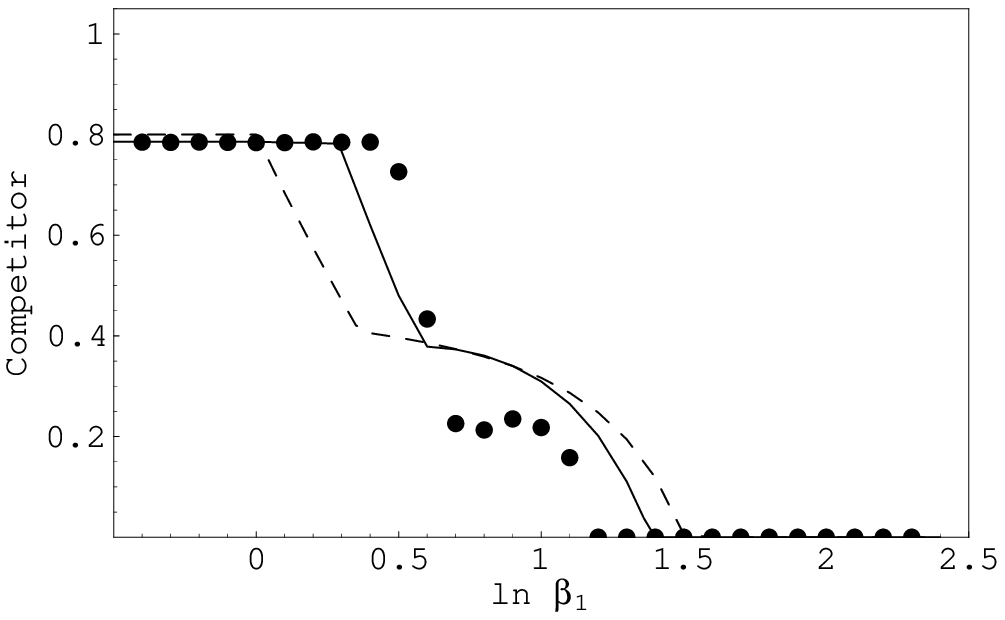}

\end{figure}

Figure 5d:

\begin{figure}[h]

        \includegraphics{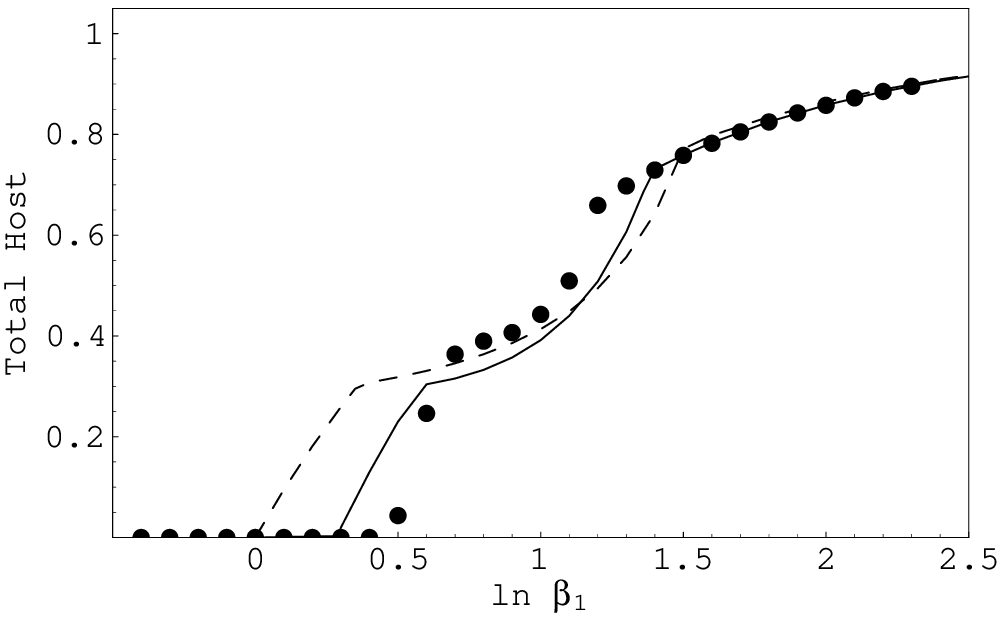}

\end{figure}

\pagebreak

Figure 6a:

\begin{figure}[h]

        \includegraphics{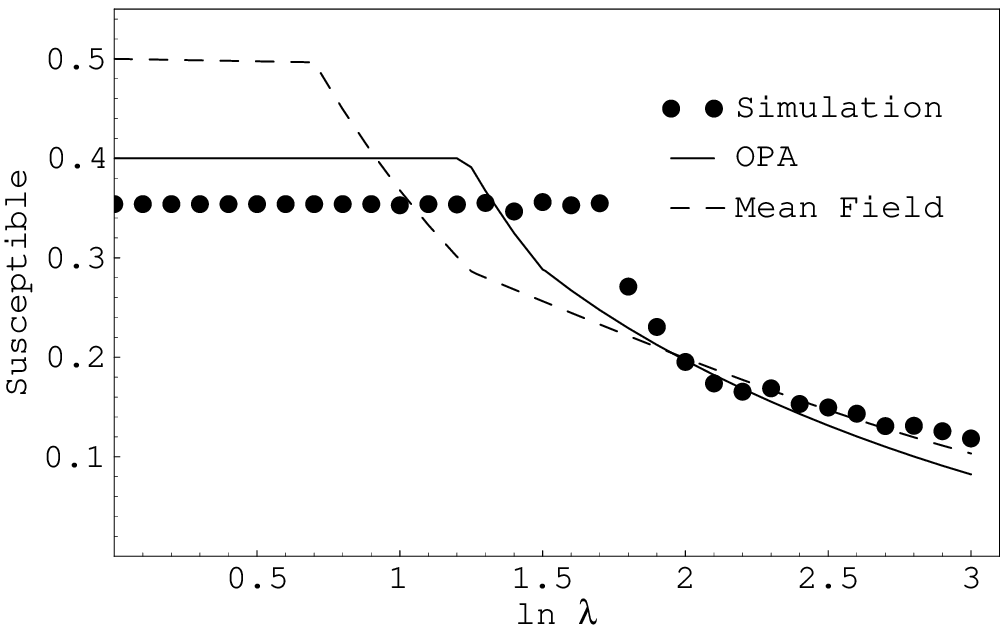}

\end{figure}

Figure 6b:

\begin{figure}[h]

        \includegraphics{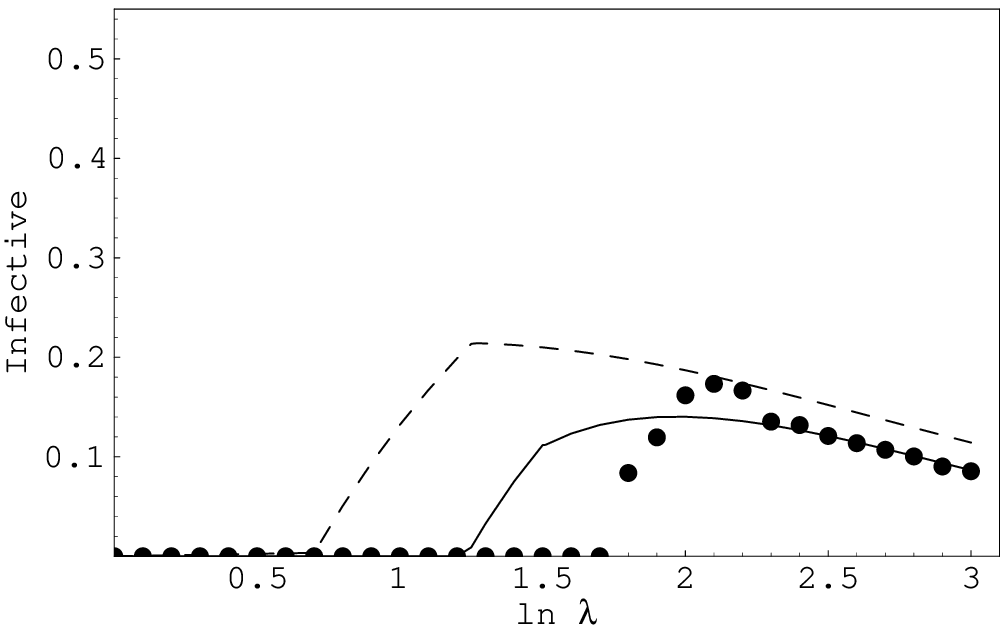}

\end{figure}

\pagebreak

Figure 6c:

\begin{figure}[h]

        \includegraphics{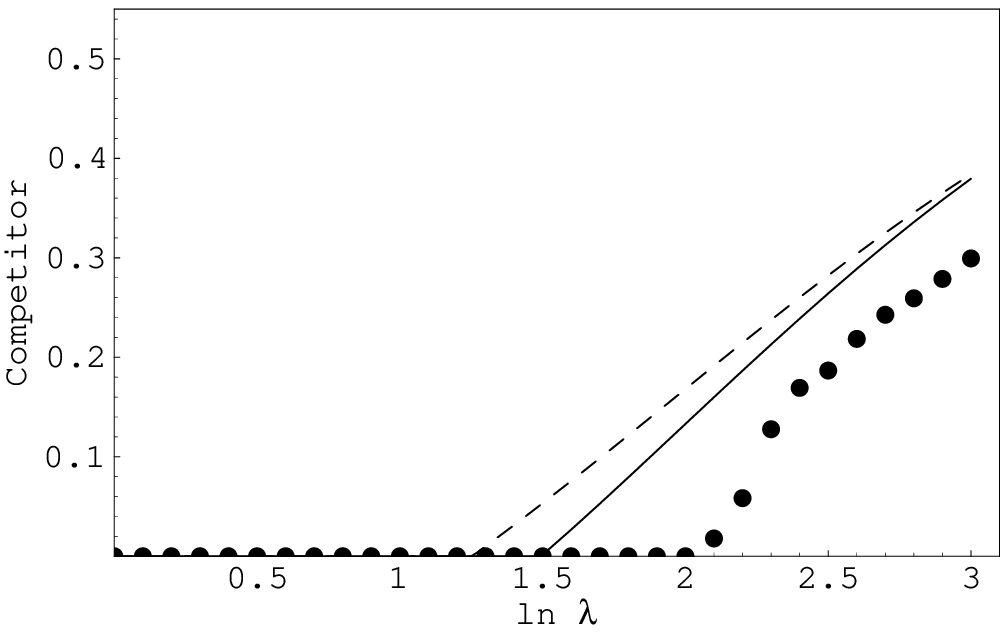}

\end{figure}

Figure 6d:

\begin{figure}[h]

        \includegraphics{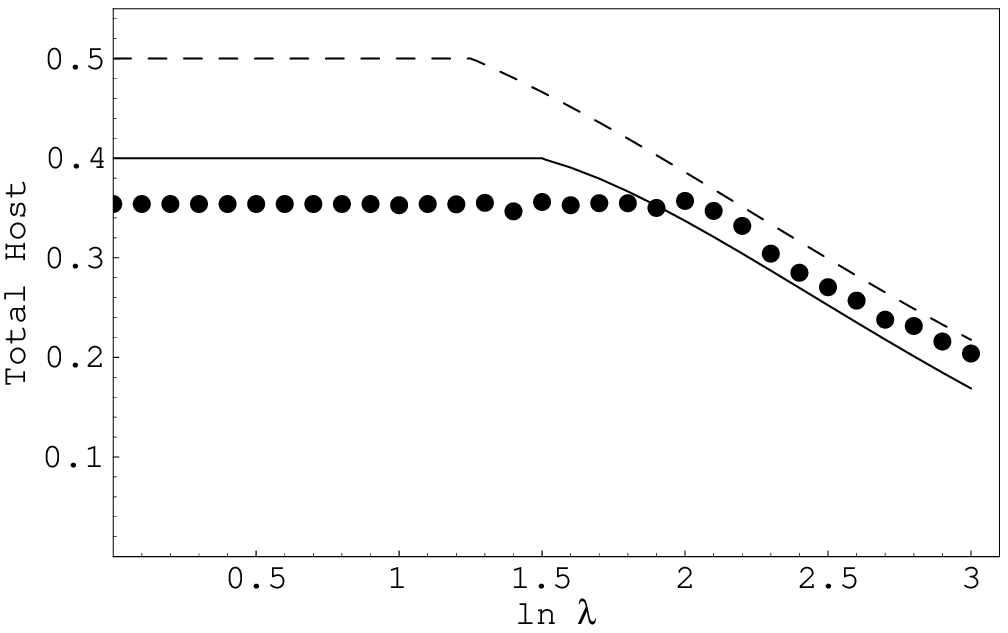}

\end{figure}

\pagebreak

Figure 7a:

\begin{figure}[h]

        \includegraphics{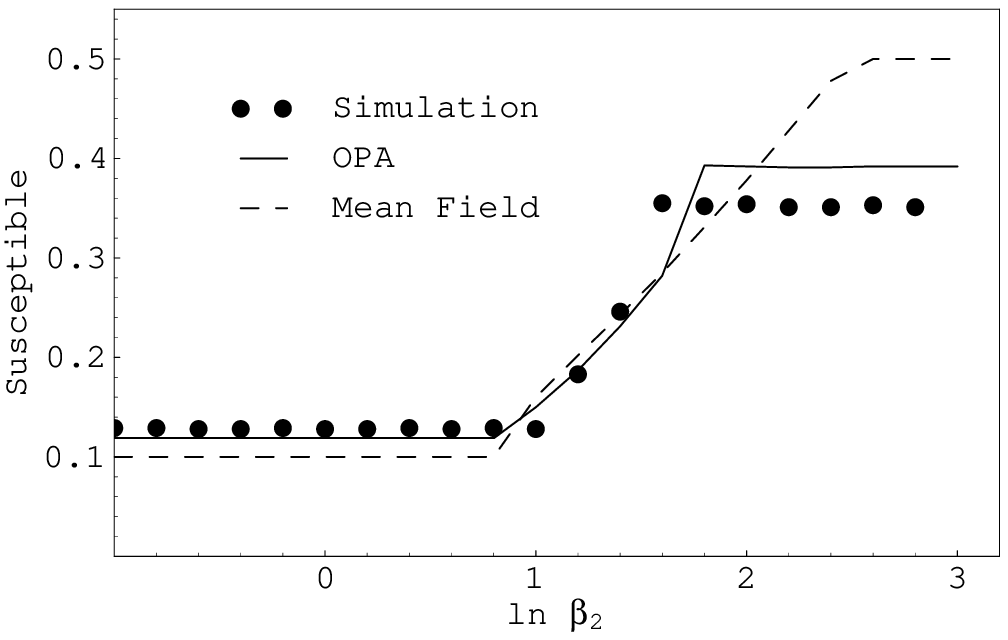}

\end{figure}

Figure 7b:

\begin{figure}[h]

        \includegraphics{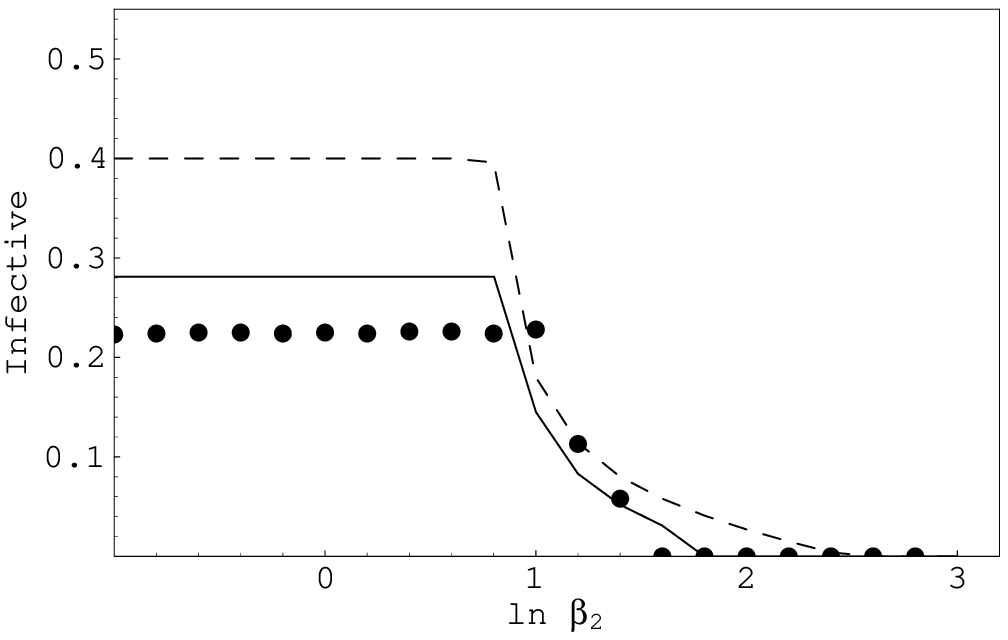}

\end{figure}

\pagebreak

Figure 7c:

\begin{figure}[h]

        \includegraphics{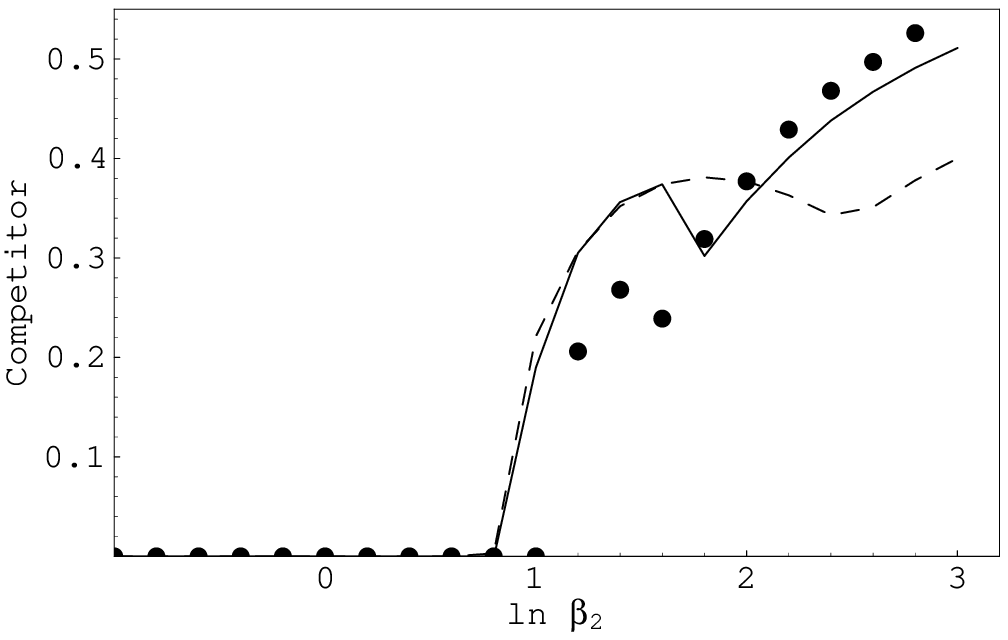}

\end{figure}

Figure 7d:

\begin{figure}[h]

        \includegraphics{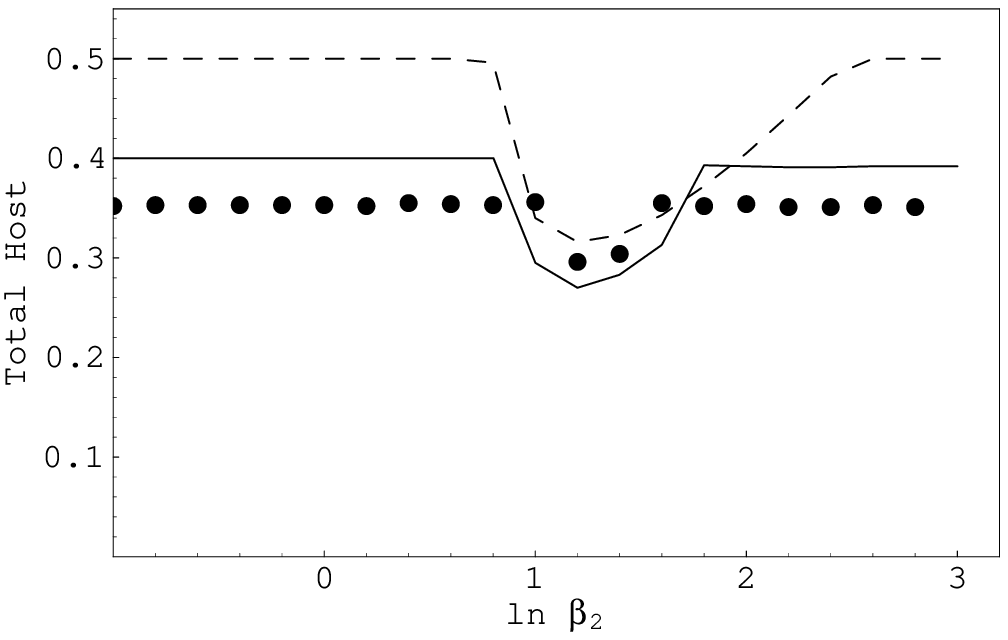}

\end{figure}

\pagebreak

Figure 8a:

\begin{figure}[h]

        \includegraphics{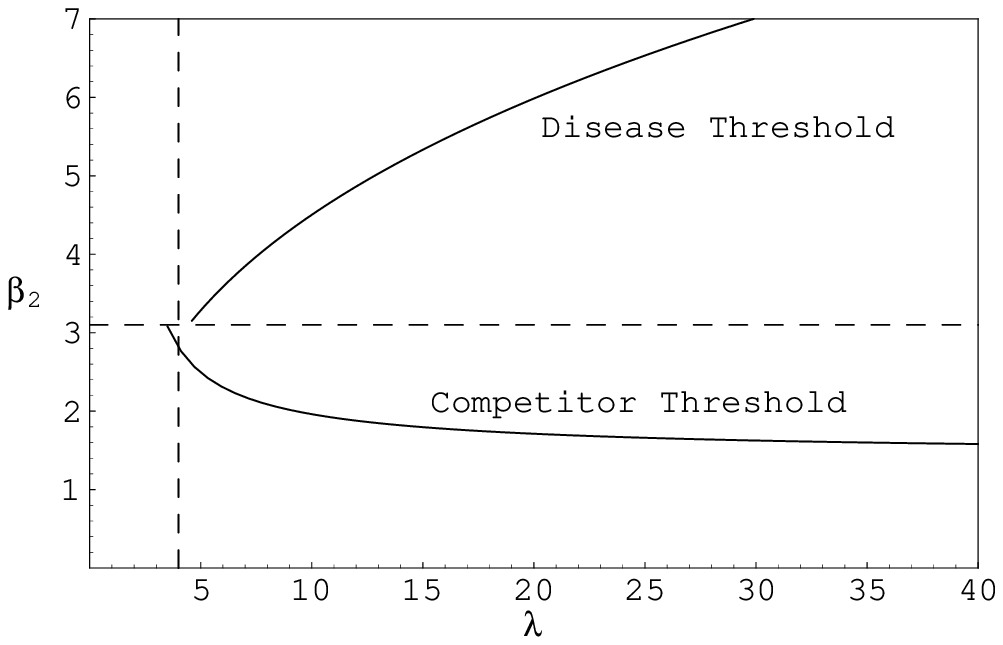}

\end{figure}

Figure 8b:

\begin{figure}[h]

        \includegraphics{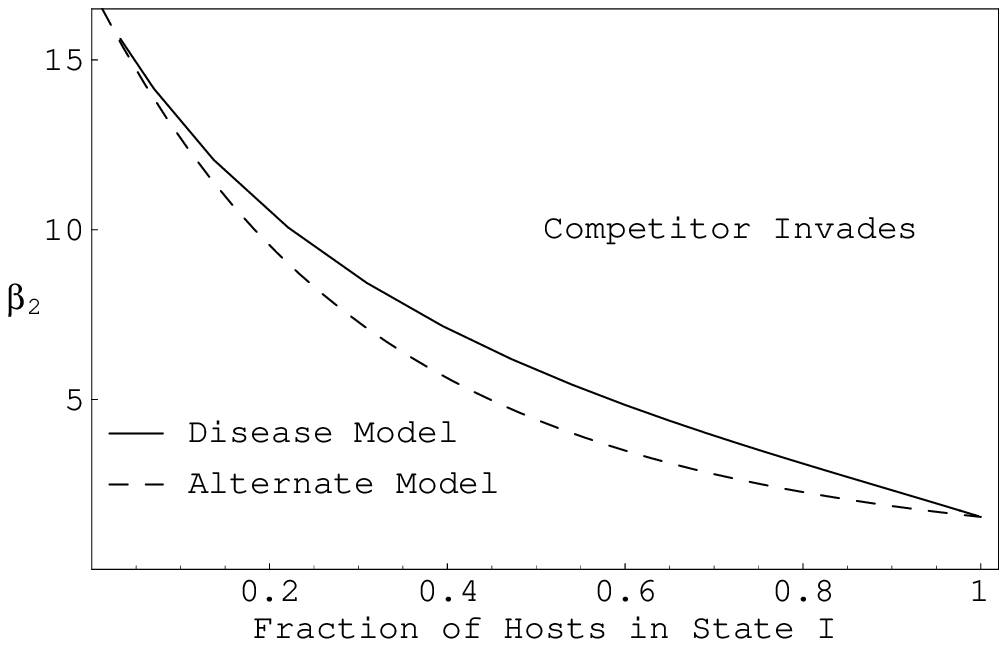}

\end{figure}

\pagebreak

Figure 8c:

\begin{figure}[h]

        \includegraphics{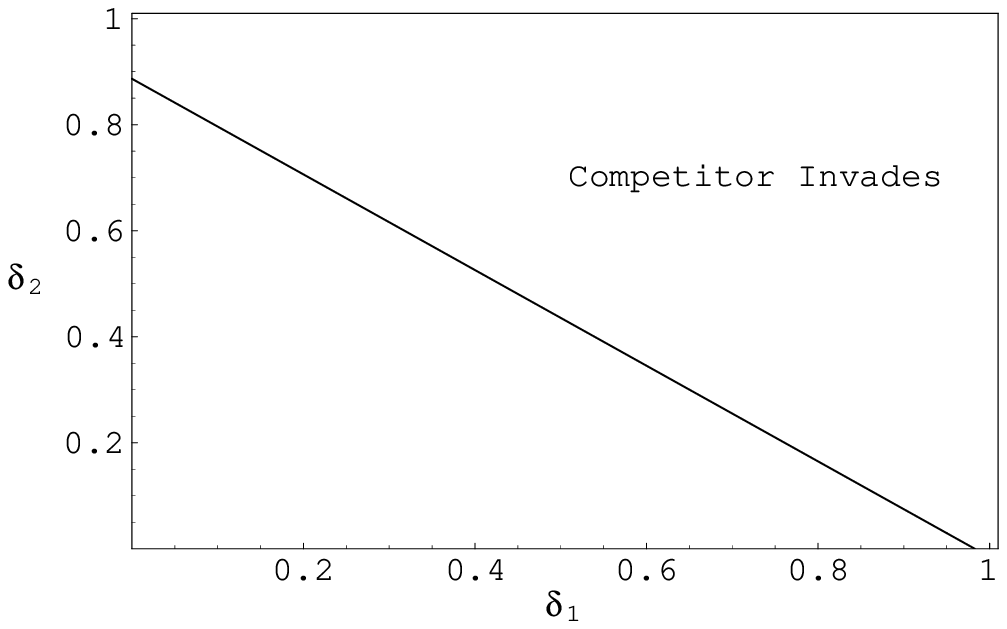}

\end{figure}

Figure 8d:

\begin{figure}[h]

        \includegraphics{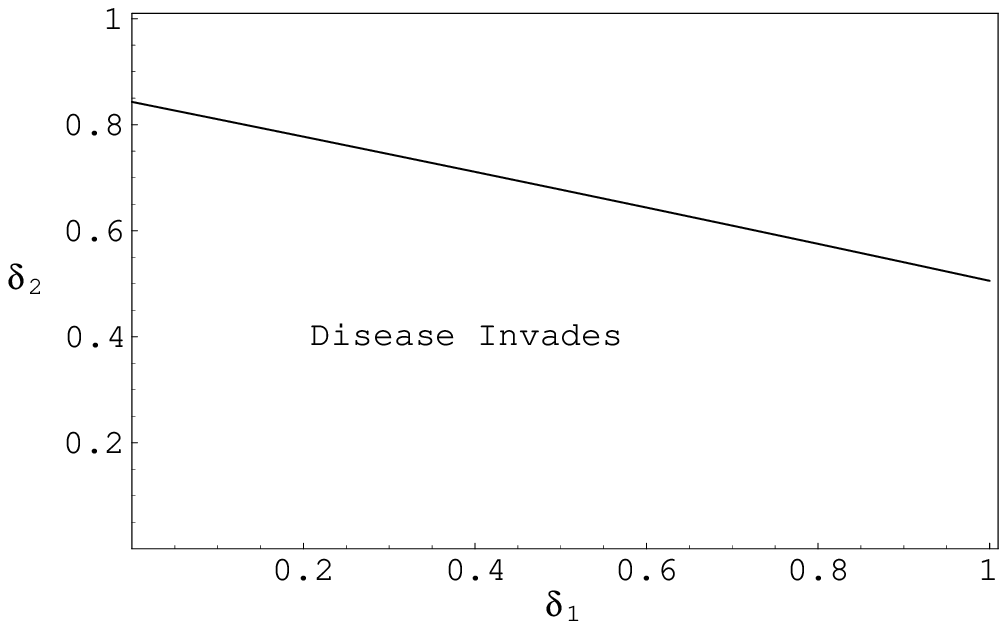}

\end{figure}

\pagebreak

Figure 9a:

\begin{figure}[h]

        \includegraphics{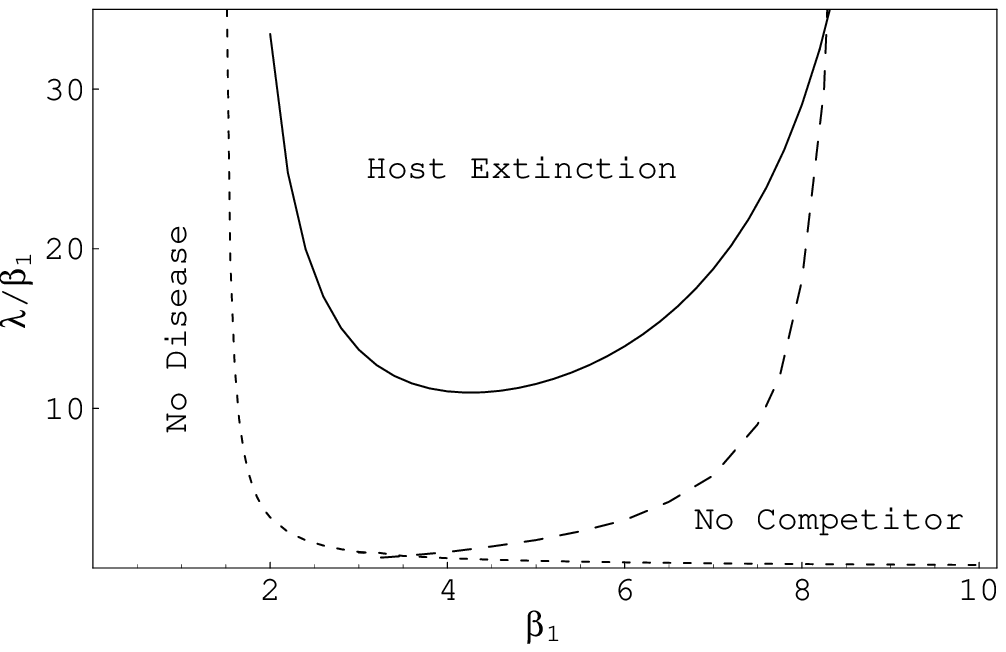}

\end{figure}

Figure 9b:

\begin{figure}[h]

        \includegraphics{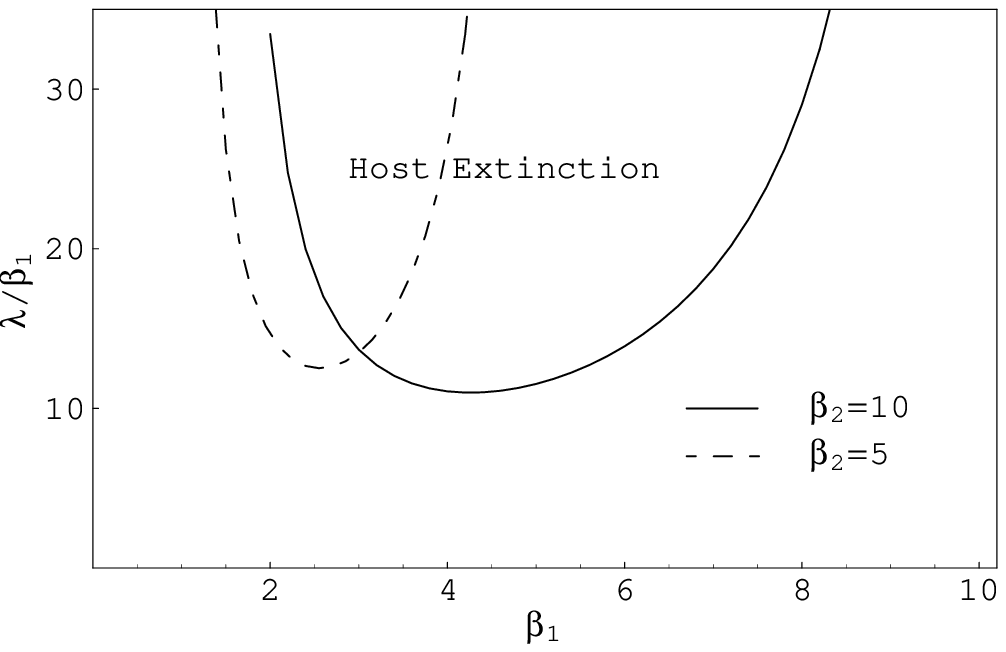}

\end{figure}

\pagebreak

Figure 9c:

\begin{figure}[h]

        \includegraphics{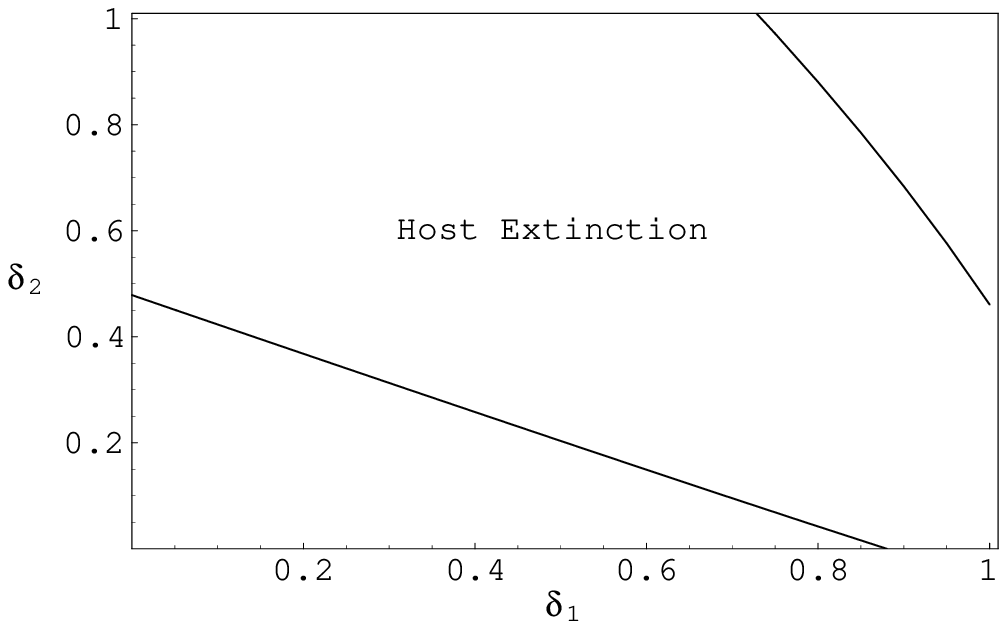}

\end{figure}

Figure 9d:

\begin{figure}[h]

        \includegraphics{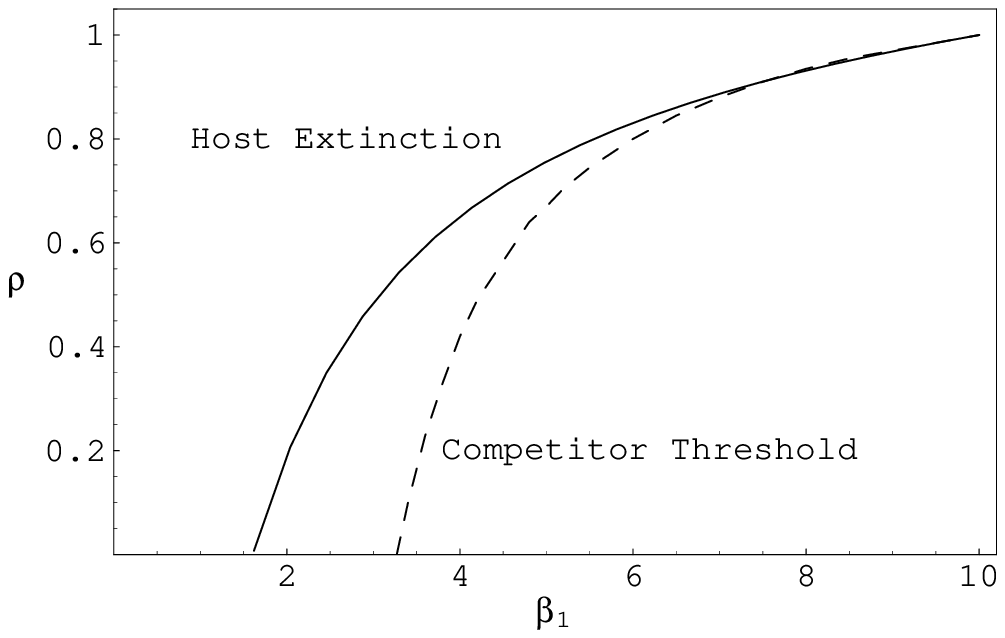}

\end{figure}

\pagebreak

Figure 10a:

\begin{figure}[h]

        \includegraphics{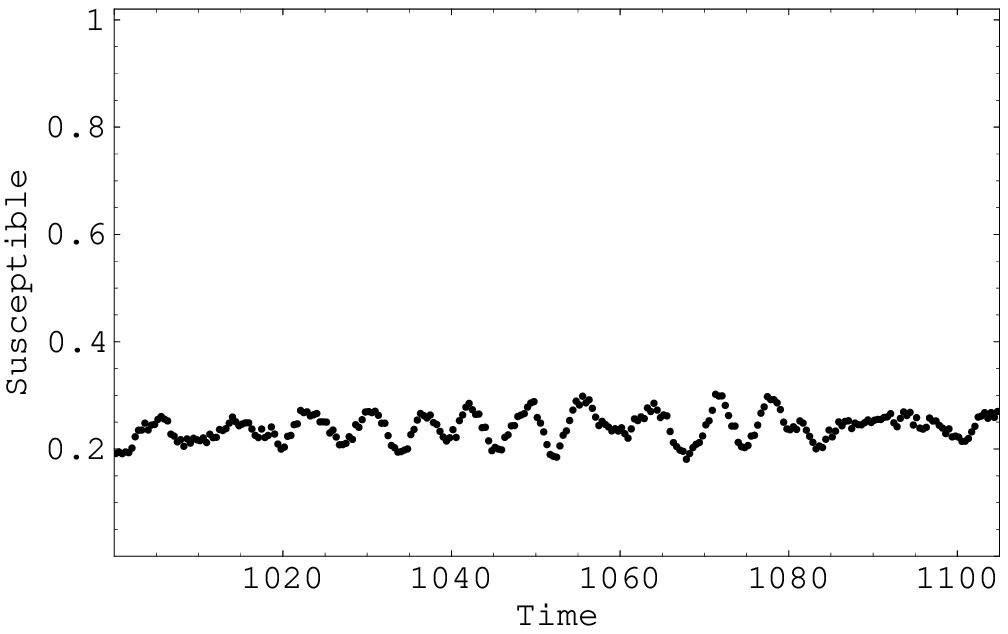}

\end{figure}

Figure 10b:

\begin{figure}[h]

        \includegraphics{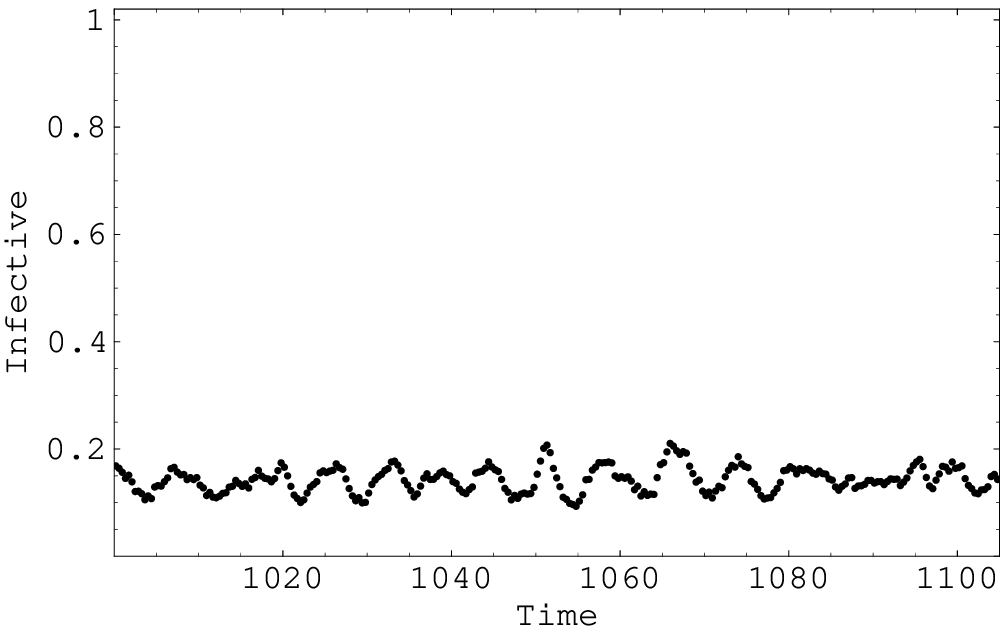}

\end{figure}

\pagebreak

Figure 10c:

\begin{figure}[h]

        \includegraphics{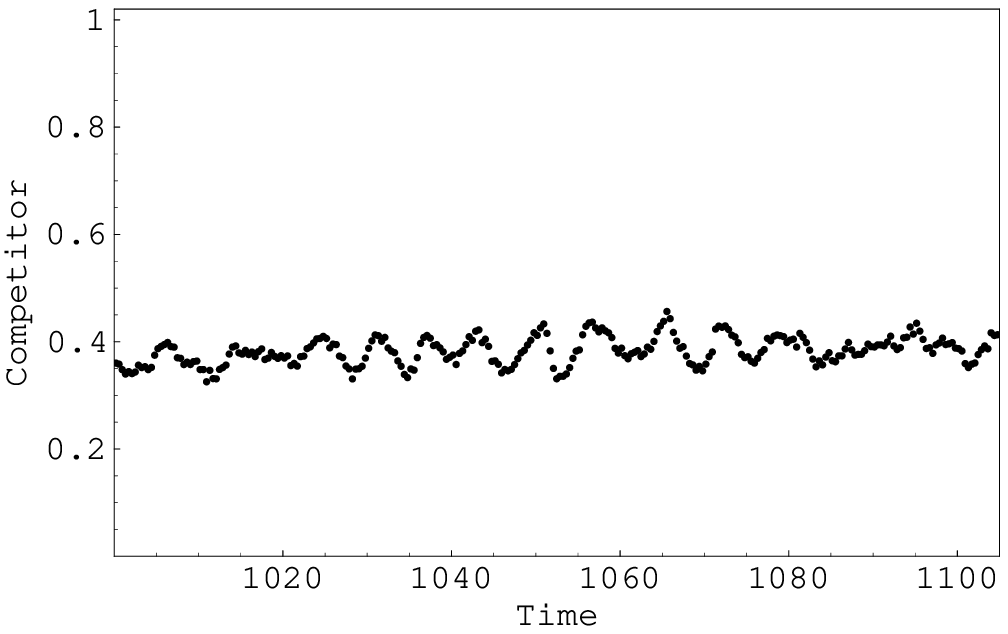}

\end{figure}

Figure 10d:

\begin{figure}[h]

        \includegraphics{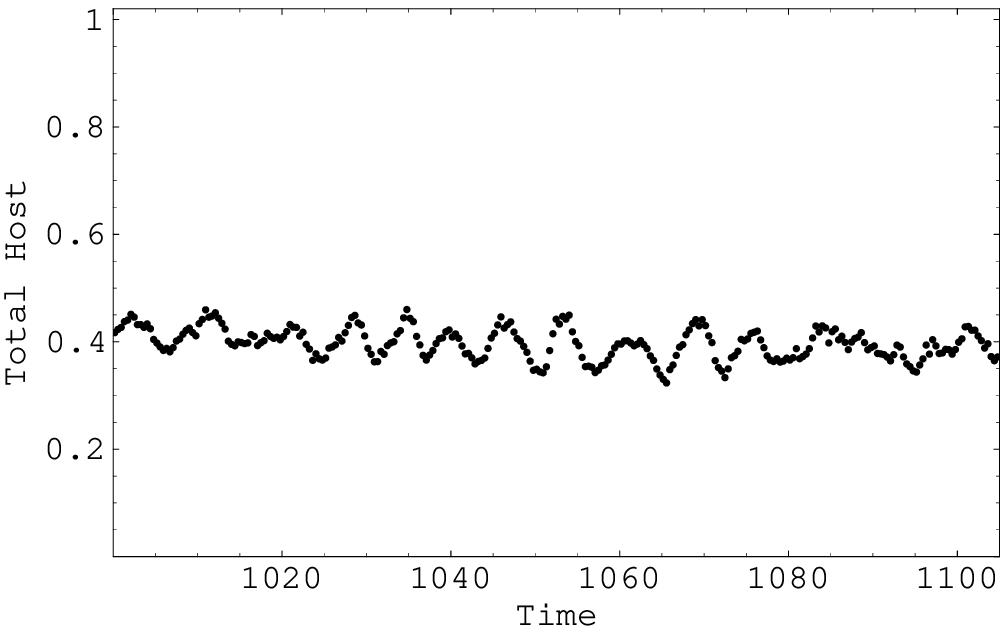}

\end{figure}

\pagebreak

Figure 11a:

\begin{figure}[h]

        \includegraphics{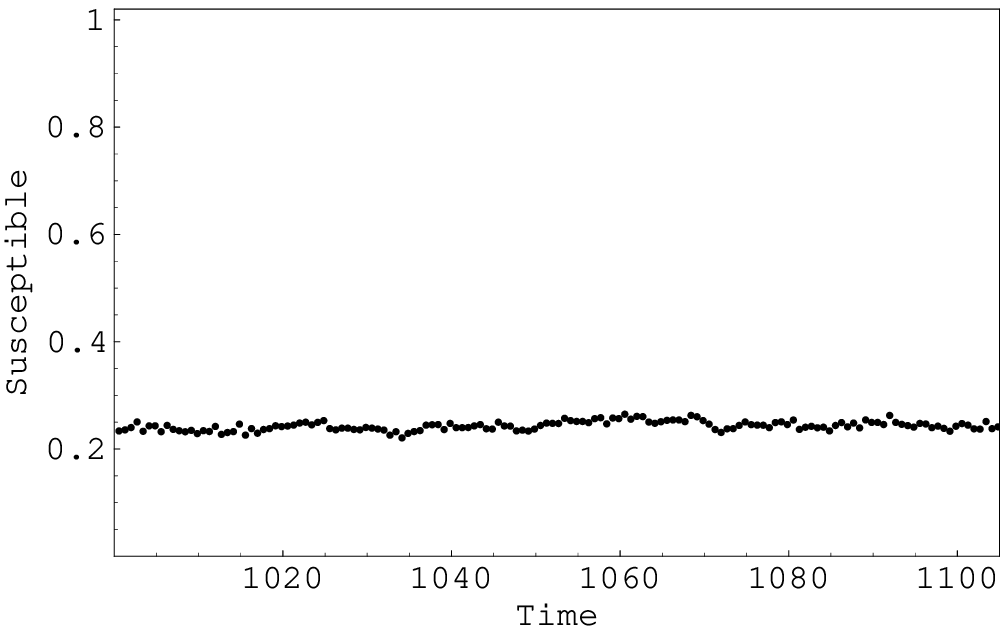}

\end{figure}

Figure 11b:

\begin{figure}[h]

        \includegraphics{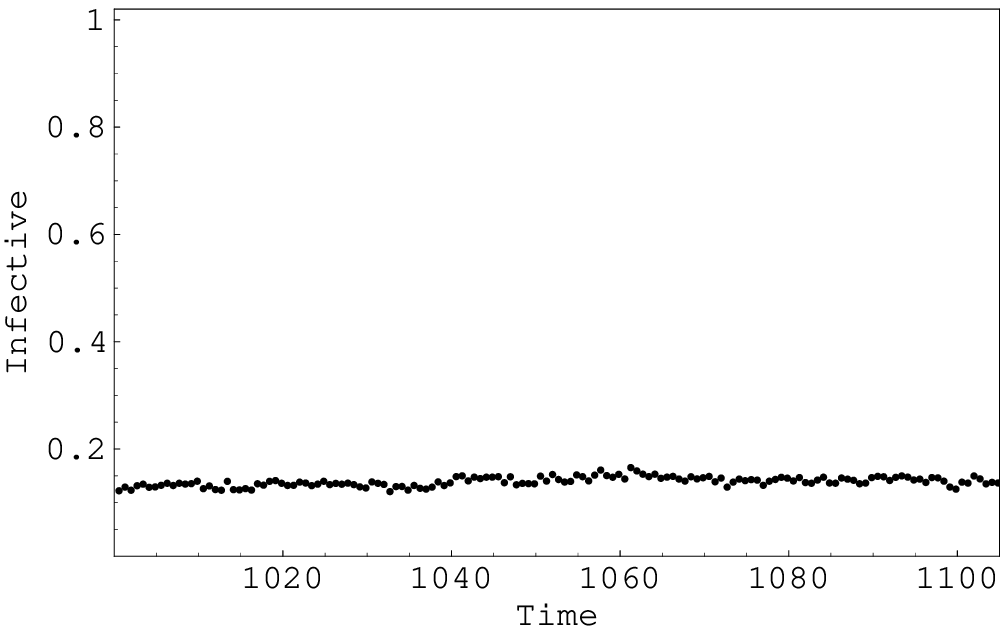}

\end{figure}

\pagebreak

Figure 11c:

\begin{figure}[h]

        \includegraphics{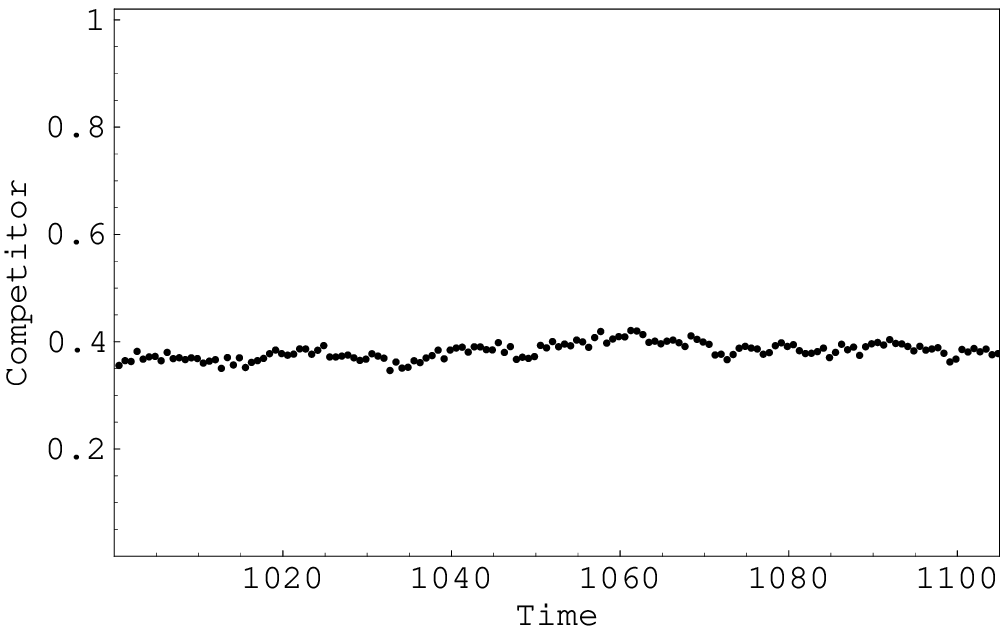}

\end{figure}

Figure 11d:

\begin{figure}[h]

        \includegraphics{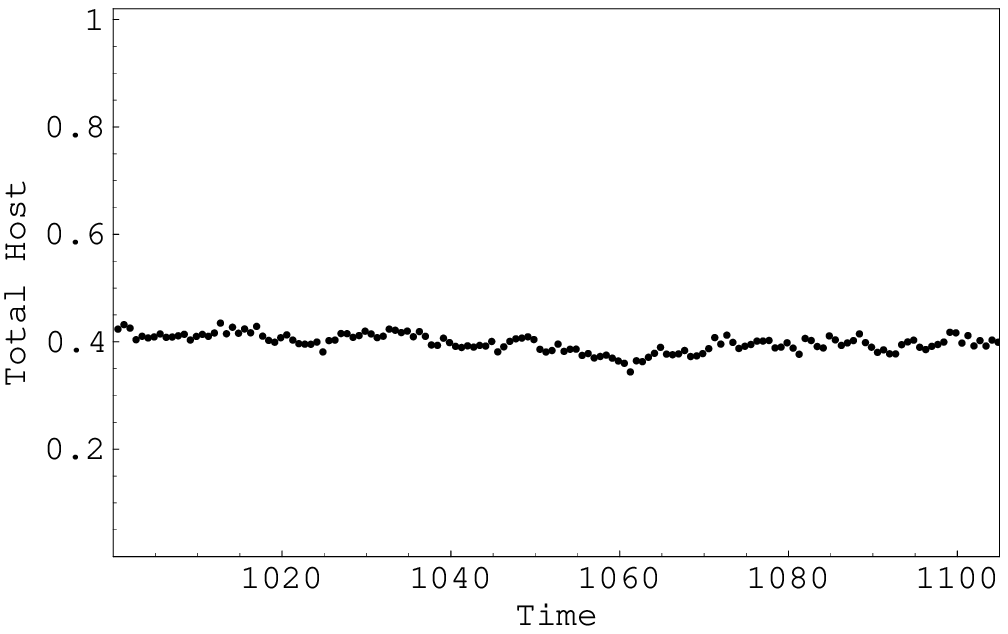}

\end{figure}

\pagebreak

Figure 12a:

\begin{figure}[h]

        \includegraphics{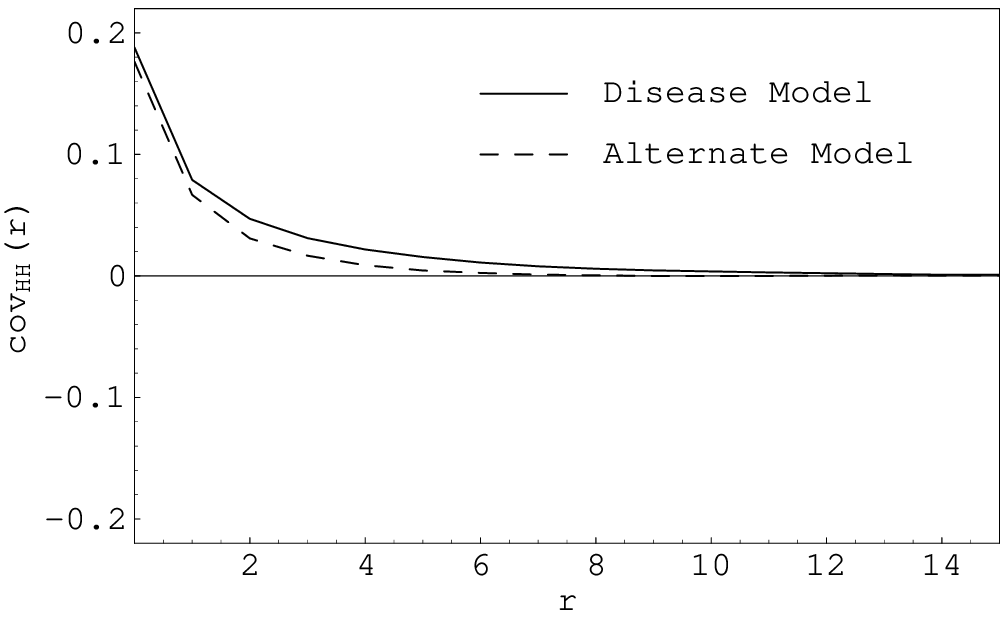}

\end{figure}

Figure 12b:

\begin{figure}[h]

        \includegraphics{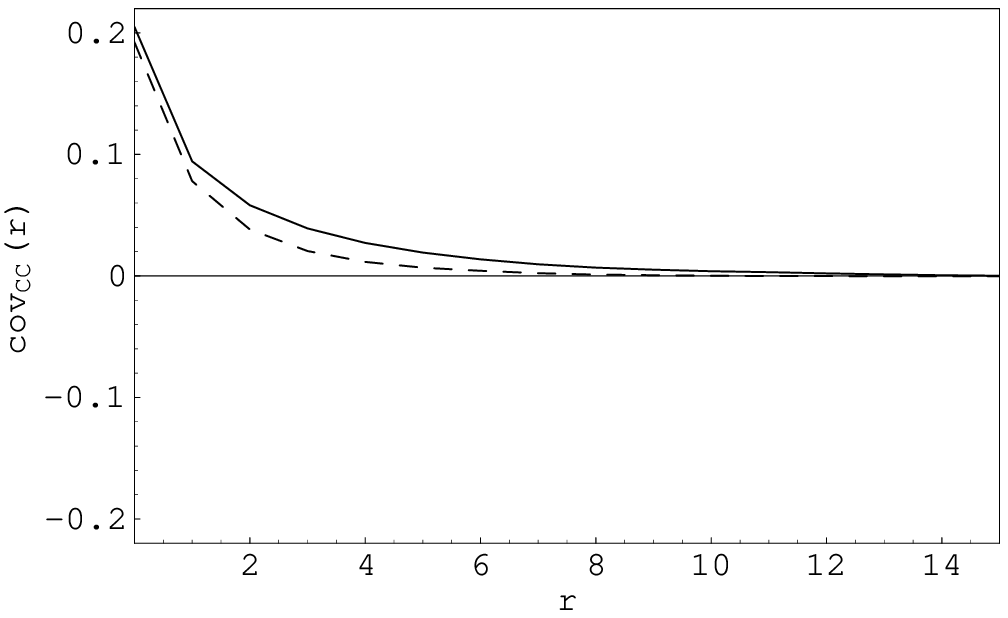}

\end{figure}

\pagebreak

Figure 12c:

\begin{figure}[h]

        \includegraphics{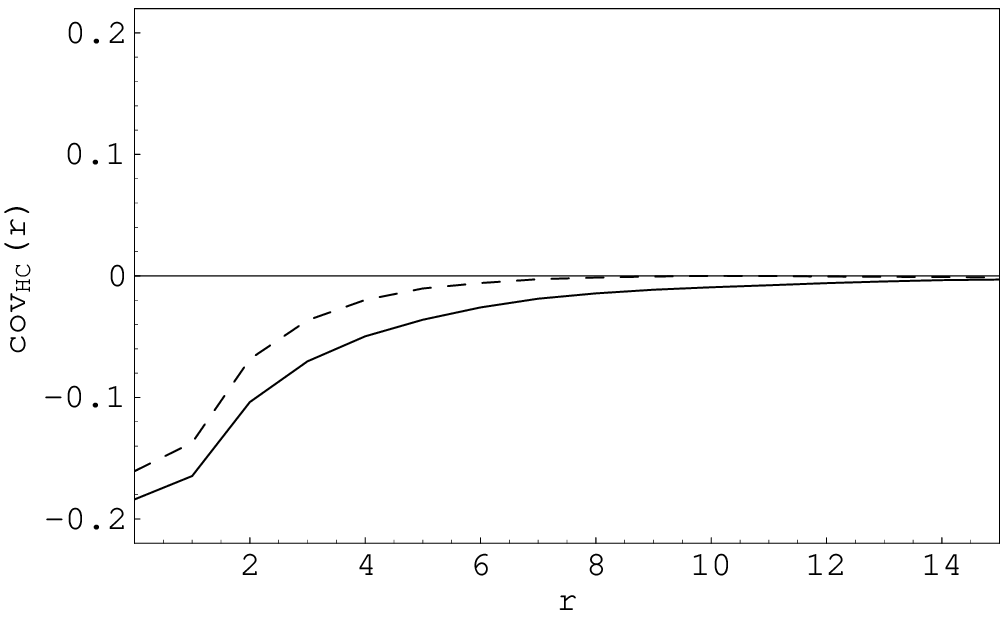}

\end{figure}

Figure 12d:

\begin{figure}[h]

        \includegraphics{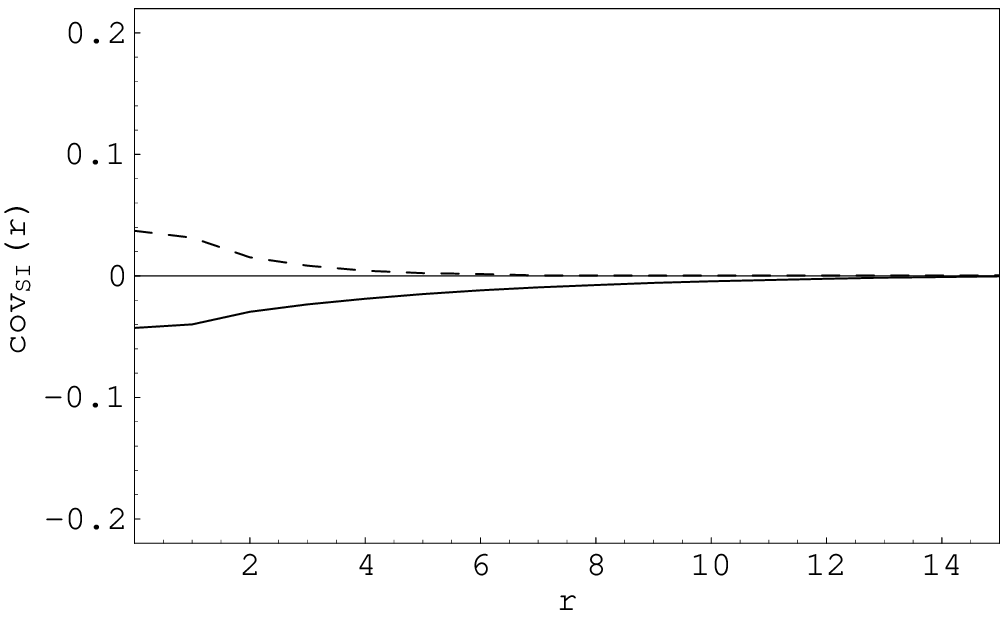}

\end{figure}

\chapter[Epidemic Threshold in Space]{Epidemic Threshold in Space}
\thispagestyle{myheadings}
\markright{}

\section*{Abstract}

	For an epidemic to occur in a closed population, the transmission rate
must be above a threshold level.  In plant populations, the threshold
depends not only on host density, but on the distribution of hosts in
space.  This paper presents an alternative analysis to an epidemic
model in continuous space proposed by Bolker (1999).  A type of moment
closure is used to determine the dependence of the epidemic threshold
on host spatial distribution and pathogen dispersal.  Local correlations
that arise during the early phase of the outbreak determine whether a
true global epidemic will occur. 

\section*{Introduction}

	One of the most important concepts to arise from epidemiological
theory is the existence of an epidemic threshold for infectious diseases
(Kermack and McKendrick, 1927).
In its most basic form, this theory states that a pathogen can only
cause an epidemic (i.e. increase from low levels) if the host population
is sufficiently large (or dense).  Equivalently, for a given host population,
a pathogen can only invade if the transmission rate is sufficiently high.
For a directly transmitted pathogen
which makes the host infectious for a finite time (after which the
host dies or recovers), the simple SIR model yields a threshold condition
that depends only on the transmission and recovery rates and on the
host population size.  The threshold criterion has been extended to
include a number of complicating factors, such as free--living parasite
stages, host behavioral heterogeneity, vector transmission, genetic 
heterogeneity,
and stochastic effects (Anderson, 1991; Nasell, 1995; 
Keeling and Grenfell, 2000; Madden \emph{et al.}, 2000).
In general, the threshold criterion can be stated as: an epidemic will
occur if and only if $R_{0} > 1$, where $R_{0}$ is the expected number
of new infections caused by a single infective individual placed in
a totally susceptible population until it recovers.  Thus, an epidemic
can occur if and only if the initial infectives more than replace themselves
before they recover.  The dependence of $R_{0}$ on various details of
disease transmission and host behavior or ecology is therefore of intense
interest.

 	The models that underlie these insights were developed primarily
for diseases of humans and other animals.  The importance of formulating
epidemic threshold criteria for diseases of plants has also been recognized
(Jeger, 1986; May, 1990; Onstad, 1992; Jeger and van den Bosch, 1994).
An essential underlying assumption of the models developed for animals
is that of mass action: it is assumed that the population is sufficiently
well mixed that, at least within subclasses, any individual is equally 
likely to come into contact with any other.  There are clearly limitations 
of this assumption for plants and other sessile organisms.  As long
as the pathogen has spatially localized dispersal (i.e. it cannot travel
from an infected host to any other in the population with equal liklihood),
some plants are more likely than others to become infected at any time.
Both the spatial structure of the host population and the dispersal pattern
of the pathogen could potentially determine whether a disease can increase
from low density in a plant population (Real and McElhany, 1996).  There is 
as yet no general theory of
how the fine scale distribution of hosts and pathogen dispersal affect the 
epidemic threshold
in plants.  In this paper, I use a simple stochastic SIR model in
continuous space to address two related questions:
\begin{enumerate}
	\item How does the epidemic threshold depend on the spatial
		distribution of host plants?
	\item How does the epidemic threshold depend on the dispersal
		distance and kernel shape of the pathogen?
\end{enumerate}

	The role of spatial structure in diseases of plants has received
a great deal of attention from experimentalists and theoreticians
(Jeger, 1989).
Despite this, there does not appear to be any clear empirical demonstration 
of spatial structure
affecting an epidemic threshold in plants.  In fact, experimentalists
appear to be less concerned with the existence of a threshold than
theoreticians.  The question of whether an epidemic has occurred is
qualitative, while experiments more readily provide quantitative
measures such as the amount of diseased host material.  As a result, 
experimental
studies on the effects of host spatial distributions have focused on
the size of epidemics, and have not demonstrated the ability of spatial
factors to switch a system between being able or unable to support an
epidemic.  Nevertheless, experiments that show an effect of spatial
structure on epidemic sizes do support the hypothesis that spatial structure
can affect the epidemic threshold.  Burdon and Chilvers (1976) manipulated
the spatial structure of a host plant population while keeping the overall host
density constant.  They found that for clumped hosts, epidemics progressed
more quickly at first, then later more slowly, than for uniformly distributed
hosts. They attributed this to the higher availability of susceptible
neighbors early in the clumped population, followed by the difficulty of
spreading from one clump to another.  The importance of spatial structure
for epidemics in plants has also been demonstrated by
Mundt and coworkers, who studied the effects of changing the
size of monoculture stands in intercropped plants, using experiments and
detailed computer models (Mundt and Browning, 1985; Mundt, 1989; Brophy and
Mundt, 1991).

	The effect of spatial structure on the epidemic threshold has
been investigated using several modeling frameworks.  In one approach,
the host population density is thought of as a continuous variable, a sort
of fluid medium through which the disease travels. 
This has given rise to a number of reaction--diffusion,
integro--differential, and focus--expansion models which incorporate 
different assumptions about pathogen dispersal  
(reviewed in Minogue (1989) and Metz and van den Bosch (1995)). 
When the host density is uniform, the threshold criterion is
unchanged from nonspatial models: an epidemic will occur if and only
if it would occur with global host dispersal (Holmes, 1997).  
The object of interest
then is the speed with which the disease travels through the population
from an initial focus.  More generally, when
host density varies in space, there is a ``pandemic'' threshold: if the
host density is sufficiently high everywhere, the disease will cause an
epidemic that reaches every region (Kendall, 1957; Thieme, 1977; 
Diekmann, 1978).  This framework is useful for
studying many aspects of disease spread at the geographic scale, or in
agricultural systems for which uniformly high density is the norm.  
However, it does not address spatial structure at the scale of
individuals, which can be especially important in natural systems
(Alexander, 1989).  Moreover, the role of spatial structure in models is 
often manifested only when individuals are treated as discrete units
(Durrett and Levin, 1994; Levin and Durrett, 1996; Holmes, 1997).   
 
	The epidemic threshold can depend on spatial structure at the
scale of individuals,
as demonstrated in a number of lattice based models
(Sato \emph{et al.}, 1994; Durrett, 1995; Levin and Durrett, 1996; Holmes, 
1997; Filipe and Gibson, 1998; Keeling, 1999; Kleczkowski and Grenfell,
1999).  In a lattice
model, each location (in discretized space or in a social network) is 
occupied by a single individual of some type (or perhaps is empty).  Pathogen
transmission can then only occur between individuals that lie within
some neighborhood, or are otherwise connected.  The key insight from 
these models is that local pathogen transmission causes local buildup
of high densities of infectives.  This local saturation of infection can
prevent a global epidemic from occuring if infectives are essentially
surrounded by too many other infectives, without enough susceptible 
neighbors to infect (Keeling, 1999).  As a result, the rate of transmission 
needed to
cause an epidemic may be much higher than that in an analogous mass action
model (Durrett, 1995; Levin and Durrett, 1996; Holmes, 1997; Keeling, 1999).  
These results are instructive for plant diseases, since they 
demonstrate that local spatial processes can have a strong impact on
the epidemic threshold criteria.  However, lattice models are limited in
the kind of information they can provide for plant populations.  The
fact that points on the lattice typically either are or are not neighbors 
does not
allow us to study implications of the rich variety of spatial structures
found in plant populations (Alexander, 1989), or of the shapes of pathogen 
dispersal kernels (McCartney and Fitt, 1987; Minogue, 1989).  

	Metapopulation models treat spatial processes at a larger scale
than that of lattice models (Real and McElhany, 1996; Thrall and Burdon, 1997;
Thrall and Burdon, 1999).  In a metapopulation approach, the host
population is thought of as broken into distinct patches. Within each
patch, the population is treated as well mixed; only the distribution
of patches in space affects the disease's progress.  This yields useful
information about how spatial structure at the landscape scale influences
epidemics, but it does not address issues at the scale of individual
plants.  For pathogens whose dispersal scale is comparable to the
spacing of individual hosts, we must consider spatial structure at a much 
smaller scale than that of a metapopulation.

	Another approach to studying the epidemic threshold in plants
was introduced in a nonspatial model by Gubbins \emph{et al} (2000).
They distinguished between primary inoculum from a free--living pathogen
stage, and secondary infection from contact between infected and
susceptible tissue.  In addition, they incorporated general functional
forms for the dependence of the transmission rates on the densities of
host and pathogen.  In principle, the effects of space can be incorporated
in the functional forms; for example, the effect of local saturation
of infectives could be described by a transmission rate that decreases
as the density of infectives increases.  However, the explicit dependence
of the transmission rates on the spatial structure of the hosts and
the dispersal of the pathogen is difficult to predict.  Like the lattice
models, this spatially implicit approach indicates that spatial structure
may be important, but does not provide details on how spatial processes
affect the epidemic threshold.
   
	A useful spatially explicit model framework for studying epidemics 
in plant
populations was introduced by Bolker (1999).  In his model, known as a point
process, individual plants are treated as discrete units, but
their locations are specified in continuous space rather than on a lattice.
The probability of disease transmission between two individuals is
governed by the pathogen dispersal kernel, a function of the distance
between them.  This framework allows one to study arbitrary spatial
distributions of hosts and arbitrary pathogen dispersal kernels at a fine
scale.  Because
individuals are discrete, similar issues of local disease saturation
as seen in the lattice models occur in the point process model 
(Bolker, 1999).  Bolker 
studied the SI (no recovery or death) and SIR versions of his model
using an approximation technique known as moment closure.  In this approach,
one writes down differential equations for the mean densities and
spatial covariances of susceptible and infected individuals. 
The covariances themselves depend on higher order spatial statistics,
but one achieves a closed system of equations by assuming that the
higher order statistics can be approximated in terms of means and
covariances.  Bolker found that epidemics in randomly scattered host
populations proceeded much slower than in mass action models, as local
pathogen dispersal limited the availability of susceptible hosts near
disease breakouts.
He also found that when he increased the clustering of the initial host
population, the epidemic could initially grow faster than in a mass action
model, but eventually slowed as it was limited by transmission between clusters.
Thus, mass action models will generally
overestimate the rate at which a disease invades a plant population, except
in cases where host clustering is sufficient to accelerate the epidemic.

Despite the success of his moment closure equations at predicting epidemic
dynamics over a range of conditions, Bolker was not able to use them
to compute epidemic threshold criteria.  In order to compute the threshold,
one needs to compute the spatial structure of the initial phase of the
(potential) epidemic.  In point process and lattice models, it is often
the case that the spatial structure of an invading population reaches a
pseudoequilibrium quickly, long before the overall densities equilibrate
(Matsuda \emph{et al.}, 1992; Bolker and Pacala, 1997; Keeling, 1999; 
Dieckmann and Law, 2000).
Heuristically, this occurs because the system reaches equilibrium at the
local scale more quickly than at the global scale when interactions are
localized.  Thus, if one can compute this pseudoequilibrium spatial
structure, one can use it to determine whether a global invasion can
proceed.  This was the approach used by Keeling (1999) to compute the epidemic
threshold for a lattice SIR model; in that context, the analog of moment
closure is called pair approximation.  However, there is no \emph{a priori}
guarantee that the moment equations will converge to a pseudoequilibrium
early in the invasion, and this failure to converge prevented Bolker (1999)
from computing the threshold criteria.

	In this paper, I present an alternative analysis of Bolker's
SIR model which allows computation of the epidemic threshold.  My analysis
is based on another version of moment closure, which uses different
assumptions about the relationship between the higher moments and the
means and covariances.  The resulting moment equations do have the
pseudoequilibrium behavior needed for threshold calculations.  I use them
to show how the question of whether or not an epidemic will occur
depends on the host population structure by using populations that are either
Poisson (randomly) distributed, clustered, or overdispersed.  I also
show how the threshold depends on the dispersal distance of the pathogen
and on the particular form of the dispersal kernel by comparing results
using exponential, Gaussian, and ``fat--tailed'' kernels.  Because the
moment equations are nonlinear, analytic solutions are not available.
However, this approach allows efficient numerical calculation of threshold
transmission rates as spatial parameters are varied. 

\section*{Model Formulation}

\subsection*{Stochastic Model}

	The model I use is identical to the SIR model introduced by Bolker
(1999); I review its formulation briefly here.  The model (called a point
process) treats both space and time as
continuous variables.  Interactions are local and stochastic, but the system
behaves deterministically at large spatial scales (a phenomenon known as
spatial ergodicity).  Since space is treated as homogeneous, \emph{a priori}
calculations do not depend on location; spatial structure depends on the
distance between points, but not on the locations themselves.

	Individual host plants are located at randomly chosen points in
two dimensional space, with initial density $S_{0}$.  Since we are studying
the rapid development of epidemics, no births or deaths (except due to
disease) are included.  A small fraction of the plants are initially
infected; they are chosen at random from the host population.  We ignore
any latent period, so that an infected plant is immediately infective.
An infective (I) plant can infect any susceptible (S) plant; the rate
at which this happens depends on the rate of production of pathogen
particles and the distance between the two plants.  Infected plants
die or recover at a constant rate (so that the infective period is
exponentially distributed); dead or recovered plants (R) have no bearing
on the rest of the system and are thus ignored.  

	To calculate the rate at which a given susceptible plant becomes
infected, we integrate over all space the contributions of the 
infected plants in the population.   A host at location $x$ becomes
infected at rate $\lambda \int D(|x-y|) I(y) dy$.  Here, $I(y)$ is the
density of infected plants at location $y$, and $D(|x-y|)$ is the
dispersal kernel of the pathogen.  It is normalized to be a probability
density function ($\int D(|x-y|) dy = 1$), so it weights the rate of infection
by the distance between the infective and susceptible hosts.  The rate
parameter $\lambda$ is analogous to the contact rate in mass action
models.  It is phenomenological, incorporating rate of pathogen production,
survival of the pathogen in the environment, and probability of successful
infection when a host is encountered.

	Spatial structure is incorporated into the model in two ways: the
dispersal kernel of the pathogen and the initial distribution of the hosts.
In this paper, I will assume that the dispersal kernel is a radially 
symmetric, decreasing function of distance (so that polar coordinates will
be used from now on).  This matches the dispersal patterns found for a
number of plant pathogens with various dispersal mechanisms (McCartney
and Fitt, 1987; Minogue, 1989).  However, note
that factors like vector behavior, advection, and spore aerodynamics
can give rise to different types of dispersal kernels (Aylor, 1989; 
McElhany \emph{et al.}, 1995).  Even restricted
to radially symmetric decreasing functions, a number of different dispersal
kernels for plant pathogens can be used.  I choose three simple
kernels that I believe illustrate the ways that kernel shape influences
the epidemic threshold.  As a baseline, I use a negative exponential kernel.
I compare it with a normal (Gaussian) kernel which decays more rapidly
with distance, and a ``fat--tailed'' kernel which decays more slowly.
I use normalized kernels (the integral over all space is one), so that
the kernel is a probability density function for the distance traveled
by a viable pathogen particle.
In order to compare kernels of different types, 
I follow Bolker (1999) in using the ``effective area'' of the kernel,
\begin{equation} 
A = \left(\int_{0}^{2 \pi} \int_{0}^{\infty} [D(r)]^2 dr d\theta \right)^{-1}.
\end{equation}
The term ``effective area'' comes from the fact that if the kernel is
constant on a finite disk (and zero outside it), this formula gives the
area of the disk.  Thus, I say that two kernels have the same spatial
scale if they have the same effective area. The three kernels and their
summary statistics are given in Table 1.

	I also use three qualitatively different patterns for the initial
distribution of hosts in space.  The simplest configuration is given by
a spatial Poisson process, in which the locations of plants are chosen
independently of one another.  A Poisson population has a constant probability
per unit area of having a plant, regardless of the positions of other plants.
Clumped host patterns are generated by
a Poisson cluster process (Diggle, 1983; Bolker, 1999).  In this process,
``parent'' sites are chosen by a spatial Poisson process with intensity
$\gamma$.  Around each parent site we independently place $n_{C}$ ``daughter''
plants using a host distribution kernel $H(r)$.  We discard the parent
sites, yielding a population with density $S_{0} = \gamma n_{c}$.  The
locations of the plants are no longer independent, since within clusters
the local density is higher than the overall density of the population.
Deviating from the basic Poisson distribution in the other direction, I
used a simple inhibition process (Diggle, 1983) to generate an overdispersed
host population.  Again we begin with a Poisson process of intensity
$\gamma$; this time we eliminate all plants that are within a distance
$a$ of another individual.  The resulting population has density
$S_{0} = \gamma \exp(-\pi \gamma a^{2})$.  The imposition of a
minimum possible distance between plants crudely captures patterns which
can arise from competition (Bolker and Pacala, 1997) or ``pathogen shadows'' 
(Augspurger, 1984).

	With the pathogen dispersal kernel, host distribution, and transmission
and recovery rates specified, the model can be simulated on a computer.  
Figure 1 shows snapshots from the early stages of epidemics in Poisson,
clustered, and overdispersed host populations.  The three examples use the
same host densities and pathogen dispersal.  It is clear that the host
spatial structure is playing an important role in determining the success
of the pathogen invasion.  We are limited in what we can learn about how
spatial structure shapes the epidemic threshold from simulations alone.
To gain further insight, we turn to equations which describe the temporal
evolution of the densities and spatial structure of an emerging epidemic.

\subsection*{Main Equations}

	Let $p_{SI}(r)$ be the joint density of S and I at distance $r$; that
is, it is the limiting probability of finding an S and an I individual in
small regions distance $r$ apart, as the area of the regions goes to zero.
Then since each new infection is the result of an interaction between an
S--I pair, the global densities satisfy the differential equations:
\begin{eqnarray}
\dot{I} & = & \lambda\int\int D(r) p_{SI}(r) r dr d\theta - \mu I\\
\dot{S} & = & -\lambda\int\int D(r) p_{SI}(r) r dr d\theta,
\end{eqnarray}
where $\mu$ is the recovery (death) rate.  The first thing we should do
is nondimensionalize the equations.  Since individuals are discrete, we
cannot rescale how we count them; we only need to rescale time and space.
We can rescale time by defining $\tau = \mu t$, so that
one time unit corresponds to the expected lifetime of an infected individual.
We can rescale space by defining $\rho^{2} = r^{2} S_{0}$, so that the
unit of space is that which yields an initial host density of one.
Formally, the equations can be rewritten in terms of the following
dimensionless quantities: $\hat{S} = S/S_{0}$, $\hat{I} = I/S_{0}$,
$\hat{p}_{SI}(\rho) = p_{SI}(r)/S_{0}^{2}$, $\hat{D}(\rho) = D(r)/S_{0}$,
$\hat{\mu} = 1$, and $\hat{\lambda} = \lambda S_{0}/\mu$.  For notational
simplicity, I will use the same notation as in the original equations,
using $\mu = 1$ and $S_{0} = 1$, with the understanding that all quantities
have been nondimensionalized by the procedure above.  

	Next, we define the spatial covariance $c_{SI}(r) = p_{SI}(r) - SI$, 
and the related spatial correlation, $\mathcal{C}_{SI}(r) = c_{SI}(r)/SI$.
(Note that this is a slight abuse of the usual meaning of correlation, since
we use the mean densities rather than their variances.)  Also, we define
the weighted covariance and correlation by: $\bar{c}_{SI} =
\int\int D(r) c_{SI}(r) r dr d\theta$ and $\bar{\mathcal{C}}_{SI} =
\int\int D(r) \mathcal{C}_{SI}(r) r dr d\theta$.  With this notation,
the nondimensionalized equations can be written as:
\begin{eqnarray}
\dot{I} & = & \lambda (1 + \bar{\mathcal{C}}_{SI}) S I - I\\
\dot{S} & = & -\lambda (1 + \bar{\mathcal{C}}_{SI}) S I.
\end{eqnarray}
When the correlations are zero, spatial structure disappears from the
model and we have the mass action SIR model.  Thus, the weighted correlation
summarizes the deviation from the mass action approach; it captures the
population structure ``seen'' from the point of view of an individual
using a given dispersal kernel.

	We can also summarize the spatial structure of the initial host
population in terms of spatial correlations.  For a Poisson process, 
$\mathcal{C}_{SS} = 0$ for all $r$.  For a Poisson cluster process
with density $S_{0}$, we have:
\begin{equation}
\mathcal{C}_{SS}(r) = \frac{n_{c} (n_{c} - 1)}{S^{2}} (H * H)(r),
\end{equation}
where $H * H$ denotes the convolution of the host dispersal kernel with
itself (Diggle, 1983).  
Note that the correlation is positive at all distances, and if $H(r)$ 
decreases monotonically to zero, then so does the correlation.  Finally,
the inhibition process yields:
\begin{equation}
\mathcal{C}_{SS}(r) = \left\{ \begin{array}{ll}
                        -1 & r < a \\
                        (\gamma/S)^2 \exp(-\gamma U(r)) - 1 &
                                a < r < 2 a \\
                        0 & r > 2 a,
                      \end{array}
                      \right.
\end{equation}
where $U(r) = 2 \pi a^{2} - 2 a^{2} \cos^{-1}(r/(2 a))
+ r \sqrt{a^{2} - r^{2}/4}$ is the area of the union of two circles of
radius $a$ and centers distance $r$ apart (Diggle, 1983).  Note that
the correlation is negative up to distance $a$, after which it is
positive and decreases to zero at distance $2 a$.

	As the disease invades, the SI correlations evolve.  However, when
we are computing the threshold criterion for a successful invasion, only
the initial behavior of the model is relevant.  When the epidemic is
started by randomly infecting host plants, it appears that we should use
the initial host correlation for the SI correlations in the main equations.
In this case, we would predict that the success of the invasion depends
only on the host distribution, and not on the further clustering of infected
individuals within the population.  Moreover, for a Poisson host distribution,
we would predict that the spatial threshold is the same as the mass action
one, since the host correlation is zero.  As Bolker (1999) pointed out,
this approach would be analogous to incorporating other forms of host
heterogeneity, via the coefficient of variation in the host population
(May and Anderson, 1989).  However, this approach does not capture the
full effect of spatial structure on epidemics;
the evolution of SI correlations early in the invasion is crucial in
determining whether or not a true epidemic will occur.  As a result,
even though we are focusing on the threshold criteria, it is necessary
to understand the dynamics of the correlations.

\subsection*{Correlation Equations}

	The main equations as given above exactly describe the evolution
of the mean densities; however, they include the unknown correlations.
In order to arrive at a closed model, we need to specify the dynamics
of the correlations.  One approach is to assume that $\mathcal{C}_{SI}(r) = 0$ 
for all $r$.  This is the so--called mean field assumption, and it yields
the nonspatial mass action model.  The mean field model can be seen as
the limiting behavior of the spatial model as dispersal becomes global
or the population is well mixed.  Alternatively, it can be seen as a
first approximation to the behavior of the system with local dispersal.
As Figures 2 and 3 show, the mean field assumption is a poor one when
dispersal distances are not very far; it generally overestimates the
size of an epidemic.  Moreover, it fails to capture the fact that changing
the pathogen dispersal distance can make the difference between the
success (Figure 2) and failure (Figure 3) of an epidemic.  

	In order to include spatial structure in the dynamics, we can
write down differential equations for the joint densities:
\begin{eqnarray}
\dot{p}_{SI}(|x-y|) &  =  & \lambda \int_{z\neq x} D(|y-z|) p_{SSI}(x,y,z) dz
                        \nonumber \\
                    & & - \lambda \int_{z\neq x} D(|y-z|) p_{ISI}(x,y,z) dz
                        \nonumber \\
                    & & - \lambda D(|x-y|) p_{SI}(|x-y|) - p_{SI}(|x-y|).
\end{eqnarray}
Here, $p_{SSI}(x,y,z)$ is the joint density of S at $x$, S at $y$, and I at 
$z$.  The derivation of this equation follows the standard procedure
described in Bolker (1999).  Essentially, we compute the dynamics of pairs 
of sites by
following changes to one member of the pair at a time; these changes may
be density independent, due to interaction with the other member, or
due to interactions with a third individual (hence the ``triplet''
densities).  In this case, the first term describes the creation of an
SI pair from the infection of one member of an SS pair; the second term
describes the destruction of an SI pair by infection of the S by a third
plant; the third term describes infection within the pair; the last
term describes the death of the infected plant.  Now, of course we
face the problem that the triplet densities are not known.  We arrive
at a closed model by assuming that the triplet densities can be written
in terms of mean densities and pairs.  This process, known as moment
closure, yields an approximation to the true dynamics that we hope
captures the important aspects of spatial structure.  
   
	There are several \emph{a priori} plausible ways to approximate
the triplet densities; the closure must be chosen based on the accuracy and 
utility of the resulting system (Dieckmann, 2000).  Bolker (1999) closed 
the SIR model by
assuming, for example, that $p_{SSI}(x,y,z) = S p_{SI}(|y-z|) + 
S p_{SI}(|x-z|) + I p_{SS}(|x-y|) - 2 S^{2} I$.  This approach, called
a power--1 or central moment closure, yields a system of linear
integro-differential equations for the pair densities.  Bolker found that
the system gave close approximation to the dynamics of the SIR model
provided that dispersal distances were not too short.  However, the
equations did not possess a pseudoequilibrium spatial structure during the
initial phase of the epidemic.  Thus, they could not be used to calculate
the epidemic threshold criterion.  

	As an alternative approach, I use the following asymmetric power--2
closure assumption: $p_{SSI}(x,y,z) = p_{SS}(|x-y|) p_{SI}(|y-z|)/S$,
and $p_{ISI}(x,y,z) = p_{IS}(|x-y|) p_{SI}(|y-z|)/S$.
This closure ignores any relationship between the two neighbors of the
central individual.  Since we are interested in the rate at which the
central individual changes state, we hope that this omission does not
introduce too great an error.  Law and Dieckmann (2000) used this version
of moment closure to analyze the dynamics of point process models of
plant competition.  Note that this closure assumption is
the continuous space analogue of the usual pair approximation used to
study lattice models.  Keeling's (1999) analysis of the threshold structure
of a lattice based SIR model used an approximation that includes the
ordinary pair approximation as a special case (a system with no
``triangles'', or $\phi = 0$).

	Having decided on a particular moment closure, one must decide
in what form to study the resulting equations for spatial structure.
One can write down equations for pair densities, covariances, or correlations;
when allowed to evolve over time, all give the same information.  However,
they differ in what they can tell us about invasion criteria.  Notice
that in the main equations, the SI correlation appears as a correction to
the transmission rate.  Thus, if the correlation reaches an equilibrium
early in the epidemic, we can incorporate it as a parameter rather than
a variable in the main equations; the resulting model has the same form
as the nonspatial SIR model, and threshold calculations are straightforward.
This motivates writing down differential equations for the correlations,
in order to compute their equilibria.  Closing the pair equations with
the asymmetric power--2 closure, 
and using $\dot{\mathcal{C}}_{SI}(r) = \frac{1}{SI} (\dot{p}_{SI}(r) 
- (\mathcal{C}_{SI}(r) + 1)(S \dot{I} + I \dot{S}))$ yields:
\begin{eqnarray}
\dot{\mathcal{C}}_{SI}(r) & = & \lambda\left[S (\bar{\mathcal{C}}_{SI} + 1)
                                (\mathcal{C}_{SS}(r) - \mathcal{C}_{SI}(r))
                                - D(r) (\mathcal{C}_{SI}(r) + 1)\right] \\
\dot{\mathcal{C}}_{SS}(r) & = & 0.
\end{eqnarray}

	The early spatial structure of the invasion can be found by
solving the SI correlation equation for a pseudoequilibrium using
$I=0$, $S=1$, and the initial host correlation for $\mathcal{C}_{SS}$.
This pseudoequilibrium satisfies:
\begin{equation}
\mathcal{C}_{SI}(r) = \frac{(\bar{\mathcal{C}}_{SI} + 1) \mathcal{C}_{SS}(r)
                        - D(r)}{\bar{\mathcal{C}}_{SI} + 1 + D(r)}
\end{equation}
Notice that this equation always has the spurious solution 
$\mathcal{C}_{SI}(r) = -1$.  I cannot solve the equation analytically for
the pseudoequilibrium correlation.  However, the equation can
be solved numerically as a fixed point problem using the 
$\mathcal{L}_{1}$ norm.  That is, we plug a trial solution into the right
hand side, and iterate until the integral of the difference between
the left and right hand sides is as small as we want.  From this
solution we calculate the weighted SI correlation that we include in
the main equations as a parameter; the correlation depends on both
the initial host structure and the early spread of the disease.  I
call the resulting system the mean field correlation model (MFC).
As a check on the accuracy of the moment closure assumption itself,
I also integrate the correlation equations over time; I call this approach
the moment model.   Notice that the pseudoequilibrium correlation
depends only on the spatial parameters of the model; the rate parameters
($\mu$ and $\lambda$) affect how quickly the spatial structure develops,
but not its form.

\section*{Results}

\subsection*{Model Performance}
 
	Before I present the MFC model's threshold predictions, we must
consider general features of the spatial SIR system and the MFC model's 
performance in predicting the outcome of an invasion.  First consider the
case of a Poisson distributed host.  Here, nonzero SI correlations arise
only because of local pathogen dispersal, which creates pockets of high
disease density near initial foci.  This local disease saturation in turn
slows the growth of the epidemic by limiting the supply of susceptibles
``seen'' by infective plants.  As a result, epidemics can be much smaller
than predicted by a mean field model (Figure 2).  For moderate pathogen
dispersal distances (an effective area of 20 in this example), there is
little difference between the moment and MFC models.  Both correctly 
predict that the epidemic will be smaller than predicted by mean field,
but they still overestimate the size of the epidemic.  This occurs because
the negative SI correlations that arise in the moment and MFC models are
much less severe than are found in the simulations.   

	For sufficiently short pathogen dispersal, the local supply of
susceptibles may be completely exhausted almost immediately; in this case,
the epidemic fails to move beyond its initial foci (Figure 3).  Changing
the pathogen dispersal distance has no effect on the mean field prediction,
which is qualitatively incorrect.  Now a difference between the moment
and MFC models has emerged, with the former predicting a small global epidemic
and the latter correctly predicting the invasion's failure.  The
pseudoequilibrium correlation used in the MFC is of the correct magnitude,
and is sufficient to halt the disease's spread; in the moment model, the
correlations do not grow to the pseudoequilibrium level and thus are not
as effective at slowing the epidemic.

	Next, consider the effect of changing the initial host distribution.
A clustered host population can support a larger epidemic than a Poisson
distributed host (Figure 4, compare with Figure 2).  In this case, the
clustering of hosts is sufficient to cause positive SI correlations early
in the epidemic.  This accelerates the disease's spread, an effect correctly
predicted by the moment and MFC models.  However, these models overestimate
the positive SI correlation and thus the size of the epidemic.  In fact,
since the pseudoequilibrium correlation is positive, the MFC model in
this case performs the poorest, since it fails to incorporate the eventual
depletion of susceptibles and resulting negative SI correlation.  For
this set of parameters, the epidemic promoting effects of host clustering
and epidemic inhibiting effects of local pathogen dispersal nearly cancel
out, so that the mean field prediction is quite accurate.  However, this
balance is not a general phenomenon, so that we cannot rely on ignoring
spatial structure to predict the outcome of an epidemic.  

	In these three examples, the accuracy of the moment and MFC models
was limited by the fact that they underestimate the buildup of negative
SI correlations.  This is a general feature of the models, and it stems from
the form of moment closure I am using.  As noted above, this closure
assumption ignores any relationship between two neighbors of a focal plant.
However, since the epidemic creates high local densities of infectives,
a neighboring plant is more likely to be infective if other infectives
are nearby than if not.  That is, we should expect that $p_{ISI}(x,y,z)$
is greater and $p_{SSI}(x,y,z)$ is smaller than our approximations when
$|x - z|$ is small.  Correcting this would result in a decrease in the
density of SI pairs (equation 8), yielding larger negative correlations.
Relationships between the neighbors can be incorporated by using
alternative closure assumptions or tuning the current one, as discussed
below.  However, this comes at the cost of greater computational cost
and the possible failure of the pseudoequilibrium approach through
non--convergence.  Since my moment closure assumption appears to capture
the qualitative effects of spatial structure on the success of an epidemic,
I will use it to compute threshold criterion.  It should be remembered
throughout that the MFC model's underestimate of negative SI correlations
means that the actual transmission rate needed to cause an epidemic will
in general be greater than predicted (Figure 5).  Thus, to the extent that 
it predicts spatial deviations from the mean field model, the MFC is 
essentially conservative.  The exception to this is when MFC predicts
that host clustering will yield larger epidemics than mean field predictions;
here it overestimates the epidemic--promoting spatial effects.

	One other feature of the moment model can be seen from 
equation 10, which predicts that spatial correlations between susceptible
individuals are unchanged during the epidemic. 
In fact, however, 
correlations between susceptibles always increase during an epidemic,
because those individuals not infected tend to be clustered in areas that
the epidemic has not yet hit (Bolker, 1999).  Thus, this version of
moment closure fails to capture the effect of an epidemic on the spatial
structure of the surviving host population.  As a result, the MFC and
moment models are limited in their utility for predicting the progress
of a given epidemic.  Rather, the MFC model should be seen strictly as
a tool for calculating the threshold structure.

\subsection*{Threshold Criteria}

	When the initial host density and recovery rate are scaled to 1, 
the mean field model predicts that an epidemic will occur if and only
if $\lambda > 1$.  In the spatial SIR model, the epidemic criterion also
involves two spatial factors: the distribution of hosts and dispersal of
pathogens.  From the MFC model, we see that an epidemic will occur if
and only if $\lambda (1 + \bar{\mathcal{C}}_{SI}) > 1$.  Thus, we can
compute the epidemic threshold by varying the parameters governing
host distribution and pathogen dispersal, computing the pseudoequilibrium
SI correlation, and solving the invasion criterion for the critical
transmission rate $\lambda$.  

	First consider the case of clustered hosts.  Figure 6a shows the
epidemic threshold when the pathogen dispersal and host clustering
kernels are exponential functions.  Where the threshold surface lies
above the plane $\lambda = 1$, MFC predicts that epidemics are harder to start
than in the mean field system; where the surface is below $\lambda = 1$, 
spatial structure makes epidemics more likely.  When the pathogen dispersal 
distance is large,
the threshold converges to the mean field case regardless of host distribution,
as we would expect.  When the host clustering distance is large, we approach
a Poisson distribution (clustering is very weak).  In this limit, the
spatial threshold is strictly greater than the mean field threshold;
epidemics are easiest with global pathogen dispersal and become more
difficult to achieve as dispersal decreases.  
	
	When the host is clustered, epidemics may be either more or less likely 
than mean field
theory predicts, depending on the pathogen dispersal distance.  Figure 6b
shows a typical cross--section of the threshold surface with constant
host clustering.  (Each cross--section of surface has this form, if
continued to sufficiently small values of $A_{d}$.)  As pathogen dispersal
distance decreases from the global case, initially epidemics are more likely
to occur.  This occurs because localized dispersal allows the pathogen
to take advantage of the local abundance of susceptible hosts in a 
clustered distribution. However, if the pathogen dispersal distance is
too short, it quickly depletes the supply of susceptibles even in a
clustered population, preventing a true epidemic.  As a result, for any
given clustered host population, there is an intermediate pathogen
dispersal distance at which epidemics are the easiest to obtain.  The
more tightly clustered the host population, the shorter this ``optimal''
dispersal distance will be.

	The qualitative prediction that for clustered hosts, epidemics
are most likely when pathogen dispersal is intermediate can be confirmed by
simulations (Figure 7).  Here, the transmission rate is slightly below
1.  For global pathogen transmission, no epidemic occurs.  As we decrease
the dispersal distance, we pass through the threshold, and an epidemic
occurs, infecting around 15\% of the hosts before it runs its course.
For extremely short distance pathogen dispersal, there is an initial
burst of infections, but the disease quickly burns out without spreading
to an appreciable portion of the population; we have passed back through
the threshold. 

	For a fixed pathogen dispersal distance, the effect of changing
the host clustering is simple.  For any given (finite) $A_{d}$, the critical 
value of $\lambda$ decreases as $A_{h}$ decreases.  The reason for this
is clear: as long as pathogen dispersal drops off with distance (as in 
the kernels I employ here), increasing host proximity makes transmission
more likely.  Since the dispersal kernel has a maximum at $r = 0$, an epidemic
would be most likely if all plants occupied the same point in space.  This
qualitative prediction is also confirmed by simulations (Figure 8).  With
local pathogen dispersal and $\lambda = 1$, we find that no epidemic
occurs in the Poisson distributed host.  As we decrease $A_{h}$ so that
the host is clustered, we pass through the threshold and obtain epidemics
that increase is size as clustering increases.   There is an important
caveat regarding the effect of host clustering on epidemics.  Although
the threshold transmission rate decreases as host clustering increases,
one cannot assume that the size of the epidemic increases monotonically
with host clustering.  Indeed, if hosts are packed into tight groups that
are far from one another, the disease may find it difficult to spread
between clusters.  In that case, the final size of the epidemic will be
limited by the number of clusters that are initially infected.  Thus,
tight clustering may promote the occurence of an epidemic while simultaneously
limiting its final size.  Watve and Jog (1997) found that an intermediate
cluster size minimized the size of an epidemic because of this tradeoff
between within--cluster and between--cluster spread. 

	Next, consider the case when hosts are overdispersed; i.e. there is
a minimum distance $a$ between them.  When this inhibition distance is 0,
we have a Poisson distributed host.  As the inhibition distance increases,
the MFC model predicts that epidemics become more difficult (Figure 9). 
Note that this effect is weak unless $A_{d}$ is very small, since over
most of the parameter range, the mean pathogen dispersal distance is much
greater than the inhibition distance.  For a given initial host density,
there is an upper limit to the inhibition distance we can impose and
still achieve the required density; for $S_{0}$ scaled to 1, the maximum
inhibition distance is $\frac{1}{\sqrt{e \pi}} \approx 0.34$.  For a
fixed host distribution, the threshold's dependence on $A_{d}$ is 
qualitatively like the Poisson case; the threshold transmission rate is
strictly greater than 1 and decreases to 1 as dispersal becomes global.
Simulations support the prediction that epidemics decrease as the inhibition
distance increases (Figure 10), although it is not clear in this case
that we have passed through a threshold.    
	
	Thus far, we have used exponential kernels both for host clustering
and pathogen dispersal.  Next, we consider the effect of changing the
kernels' shapes.  The shape of the dispersal kernel has been shown to 
be important in determining such aspects of an invasion as the speed and
form of 
a traveling wave (Kot \emph{et al.}, 1996; Lewis and Pacala, 2000), 
with kernels that decay faster than exponentially
(thin--tailed) and kernels that decay slower than exponentially (fat--
tailed) producing qualititatively different results.  It is not clear
\emph{a priori} whether kernel shape will be important in determining
the epidemic threshold, since the pseudoequilibrium correlations may
only depend on some measure of the kernel such as effective area or
mean dispersal distance.  

	To test whether kernel shape does in fact
matter, I computed threshold surfaces using a fat--tailed kernel and
the normal kernel (which is thin--tailed) (Table 1).  The results for a Poisson
distributed host are given in Figure 11, and for clustered hosts in
Figure 12.  In the case of clustered hosts, I used the same type of
kernel for host distribution as for pathogen dispersal.  The MFC model
predicts that kernel shape is indeed important in determining the
epidemic threshold.  When kernels are scaled to have the same effective
area, thin--tailed kernels yield larger correlations, and thus a larger
deviation from mean field predictions, than fat--tailed kernels.   The
same result occurs if one uses mean dispersal distance rather than
effective area to equate kernels.  The prediction that kernel shape can
make the difference between the success or failure of an epidemic is
confirmed by simulations (Figure 13).  With highly localized dispersal
and a transmission rate well above the mean field threshold, we see
that the epidemic never gets started with a normal kernel, is checked
by strong SI correlation with an exponential kernel, and only spreads
to a significant portion of the population with a fat--tailed kernel.

\section*{Discussion}

	This analysis of the threshold structure of a spatial SIR model
by a type of moment closure gives insight both into the biology of
epidemics in plant populations and into the technique of moment
closure.

\subsection*{Moment closure}

	I was able to investigate different aspects of Bolker's (1999)
spatial SIR model than he did, by using a different moment closure
assumption.  The asymmetric power--2 closure I used yielded correlation
equations that possess a pseudoequilibrium; including this
correlation in the mean field equations yielded a model that incorporates
the spatial structure of an emerging epidemic into the simple mean
field SIR framework. Numerical computation of the pseudoequilibrium
via a fixed--point method was straightforward; this allowed me to explore
the space of spatial parameters efficiently.  The qualitative predictions
of the mean field correlation approach were confirmed by simulations.
One could attempt to determine the threshold structure directly by
relying only on simulations rather than the MFC approximation.  However,
this approach would be computationally expensive and introduce other
difficulties.  Since any simulation uses a finite population, stochasticity
can be important in determing the outcome of an invasion.  Moreover, 
criteria must be established to determine whether or not a given simulation
run qualifies as an epidemic.  This could involve investigating how results
scale with the size of the simulation, a tedious prospect.  Finally, 
simulations do not offer explanations for observed phenomena; by contrast,
the MFC approach allows us to interpret results in terms of a simple
measure of spatial structure during the early phase of an epidemic.  By using
the MFC model to compute threshold structure, we are sacrificing some accuracy
for efficiency, clarity, and explanatory power.

	Although the MFC model correctly predicts the qualitative dependence
of the epidemic threshold on spatial factors, its predictions are not
as accurate quantitatively as one might hope.  Simulations indicate that
the true threshold transmission rate needed for an epidemic is generally
much higher than predicted by MFC.  In a addition, the size of an epidemic
predicted by MFC is often a gross overestimate.  Because it only incorporates
the early spatial structure, the MFC model is perhaps inappropriate for
predicting the full time series of an epidemic.  However, even the full 
moment model (the correlation equations evolving over time) is less accurate
than the time series predictions of Bolker's (1999) approach.  It appears
that the moment closure assumption I have used systematically overestimates
the size of an epidemic.  This probably occurs because I ignore relationships
between the two neighbors of a focal individual; since infections occur
in clusters, the states of two neighbors of an individual are in fact
likely to be correlated.  By ignoring such ``inter--neighbor'' correlations,
this moment closure effectively assumes that infective individuals are
spread more evenly throughout the population than they really are. 

	There are at least two ways one could try to correct this error
of the moment closure.  First, one could ``tune'' the closure assumptions
according to the states of the individuals involved.  For example, 
consider the effect of two infective neighbors on a focal susceptible
individual.  We could assume that $p_{ISI}(x,y,z) = (1+\epsilon) p_{IS}(|x-y|) 
p_{SI}(|y-z|)/S$ for some positive $\epsilon$.  That is, a neighbor
of the susceptible plant is more likely to be infective if another
neighbor is infective.  This is analogous to the ``improved pair
approximation'' introduced by Sato et al. (1994), who found it to be useful
for predicting the quantitative and qualitative outcomes of epidemics in
a lattice model.   However, it could be difficult to determine an
appropriate value of $\epsilon$, especially since it should in principle
depend on the distance $|x-z|$ between the two neighbors.  Alternatively,
one could explicitly include inter--neighbor terms (in this example,
$p_{II}(|x-z|)$) in the closure approximation.  Such terms are included
in Bolker's (1999) power 1 method; they are also included in the symmetric
power 2 and power 3 methods described by Dieckmann and Law (2000).    
These alternative closure assumptions are worth investigating,
although it is not clear that they will improve accuracy or possess
the pseudoequilibrium property needed for threshold calculations.
Note that including inter--neighbor terms in the closure approximation
yields second moment equations that involve convolution integrals.
Evaluating convolution integrals at each step in an iterative
procedure to calculate the pseudoequilibrium could add a significant
computational burden. 

	As Dieckmann and Law (2000) have pointed out, there are a number of
plausible moment closure assumptions that one can make; they advocate a trial
and error approach in which one compares the various moment equations to
simulations to determine which version is the most suitable for a given
system.  My study illustrates another aspect to the problem: one must
choose the closure based not just on its accuracy, but on its ability
to answer the questions of interest.  For the SIR system, Bolker's (1999)
approach appears to predict epidemic time series more accurately than
mine.  On the other hand, my approximations yield the pseudoequilibrium
behavior needed to compute the epidemic threshold.  When particular
models are studied intensively using different closure assumptions,
probably we will find not that there is a ``best'' closure, but that
the utility of the various versions depends on the questions being asked,
with tradeoffs between accuracy, tractability, explanatory power, convergence
properties, and computational cost.  While the closure assumptions I used
may be changed to try to improve the accuracy of the predictions,
it is unclear whether any gains would justify the additional model complexity
and computational cost.  My analysis clarifies the qualitative dependence
of the epidemic threshold on spatial factors; for detailed predictions
of thresholds or dynamics in a particular system, one would want to rely
on detailed simulation models rather than the simple one presented here.

\subsection*{Biological insights}

	My analysis of a simple spatial SIR model shows that the fundamental
question of whether a disease can cause an epidemic in a sessile population
depends not only on the rate of pathogen production, recovery rate, and host 
density, but also on the interaction between pathogen dispersal and host
spatial structure.  The insight that the epidemic threshold depends on
spatial factors arose in a model that treats hosts as discrete units,
rather than a continuous quantity (as in PDE models).  The central result,
analogous to lattice model results, is that local pathogen
dispersal tends to cause local saturation of the disease; the spread
of the epidemic is checked if the local (rather than global) supply
of susceptible hosts drops below a critical level.  Clustering
of hosts increases the local supply of hosts and promotes the occurence
(although not necessarily the size) of epidemics; overdispersal of
hosts has the opposite effect.

	The MFC model I used to analyze threshold structure incorporates
the spatial structure of an emerging epidemic into the transmission 
parameter of the mass action SIR model.  This allowed explicit computation
of how the critical transmission rate needed for an epidemic depends on the 
details of pathogen
dispersal and host distribution.  The analysis yielded five qualitative
predictions:

\begin{enumerate}

	\item  When hosts are distributed randomly (Poisson) or are 
		over--dispersed, the critical
		transmission rate increases from the mean field prediction
		as the pathogen dispersal distance decreases from infinity. 

	\item When hosts are clustered, there is an intermediate dispersal
		distance at which the critical transmission rate is lowest;
		longer dispersal fails to take full advantage of locally
		high host densities, while shorter dispersal leads to
		local over--saturation of infectives.
 
	\item For a given pathogen dispersal distance, increasing the
		degree of host clustering lowers the critical transmission
		rate.

	\item For a given pathogen dispersal distance, increasing the
		degree of over--dispersal of the host raises the
		critical transmission rate.

	\item The critical transmission rate depends not only on the mean
		dispersal distance or effective area of the pathogen
		dispersal kernel, but on the kernel's shape; fat--tailed
		kernels lead to less local saturation of infectives and
		thus have a lower epidemic threshold than thin--tailed
		kernels.

\end{enumerate}
I verified the validity of these predictions by simulations of the full
stochastic model.  To the extent that the MFC model predicts that spatial
structure impedes the formation of epidemics, simulations indicate that
these predictions are conservative.  For moderate pathogen dispersal
distances, the critical transmission rate for an epidemic appears to be
significantly higher than predicted by MFC.  While simulations confirm that
host clustering promotes the occurence of epidemics, it must be remembered
that this does not necessarily mean that epidemic sizes are also always
increased by clustering.  Rather, the inability of the pathogen to travel
between distinct host clusters may limit the final size of the epidemic,
even while locally high host densities promote its early growth (Watve
and Jog, 1997).

	There is much work to do on the role
of spatial structure in determining the conditions for
epidemics in plant populations.  Within the framework of simple
SIR models like the one I analyzed, there are a number of questions that
can be addressed.  We need to generate predictions
using other moment closure assumptions in order to understand more clearly
how the spatial structure of an epidemic evolves and how it affects the
course of the epidemic.  We should study the effects of dispersal kernels
that are qualitatively different from the ones that I used; advection,
vector behavior, and spore aerodynamics can produce dispersal kernels that
are not strictly decreasing with distance from the source.  This may have
profound implications for the effect of host distribution on epidemics.
In addition, host spatial distributions will often be more complex than the
patterns produced by the simple clustering and inhibition mechanisms I used.
The net effect of positive and negative host correlations at different
distances will depend on the dispersal pattern of the pathogen; for kernels
like the ones I used, host distribution at the smallest scale dominates,
but for other types of dispersal, more complex interactions may arise.

	There is also a need to extend the simple SIR model to include more
details of both hosts and pathogens, including latent periods, severity of
infection, complex pathogen life cycles, host size structure, and exogenous 
heterogeneity.  The impact of
heterogeneity in environmental factors that affect both host and pathogen,
such as light and moisture levels, is especially important.  Whether the
conditions that favor high host density also favor pathogen growth may be
critical in determining how spatial structure affect epidemics.  As these
models increase in complexity, it will be useful to find statistical
measures of spatial structure that can be incorporated in the simple
SIR framework, as was the weighted SI correlation in this study.  Finally,
it will be important to include vital dynamics of the host in the models
in order to determine how spatial factors affect the conditions for
endemicity.

	As the theory develops, it will be crucial for the models to be
constrained by data from real systems and to have their predictions tested
experimentally.  Information on pathogen dispersal kernels and host spatial
distributions in natural systems is needed to parameterize models.  
Models tuned to specific systems will then need to have their predictions
tested by experimental manipulation of host distributions and pathogen
dispersal.  Simple theory predicts that the occurence of epidemics depends
strongly on spatial factors, but we still understand little about
the structure or importance of epidemic thresholds in natural plant 
populations.

\pagebreak

\section*{Figure Captions}

Table 1: Formulas and summary statistics for pathogen dispersal and
host clustering kernels.
\newline
\newline
Figure 1: Snapshots of simulated epidemics at $t = 1$.  Light gray points are
susceptible plants; black points are infective.  Region is 50 x 50,
so that initial plant population is 2500; initially, 1\% of
plants were made infective.  Boundaries are periodic.
Pathogen transmission rate is $\lambda = 2$ and dispersal is exponential
with $A_{d} = 20$.  (a) Poisson distributed hosts.  (b) Clustered hosts,
using exponential kernel with $A_{h} = 12$ and $n_{c} = 5$ plants per
cluster.  (c) Overdispersed hosts, with minimum distance $a = 0.3$ between
plants.
\newline
\newline
Figure 2: Comparison of simulation with approximations for an epidemic
in a Poisson host population.  Parameters are as in Figure 1a.  MF = mean
field, MFC = mean field correlation model, Moment = full correlation
equations, Simul = simulation of stochastic model. (a) Density of
infectives (starting from 0.01). (b) Density of susceptibles (starting
from 0.99).  (c) Weighted SI correlation, $\mathcal{C}_{SI}$ (starting
from 0).  MFC correlation indicates the pseudoequilibrium.
\newline
\newline
Figure 3: Comparison of simulation with approximations for an epidemic
in a Poisson host population.  Parameters are as in Figure 2, except
that dispersal scale is $A_{d} = 1$. Disease fadeout prevents computing
$\mathcal{C}_{SI}$ for full time interval.
\newline
\newline
Figure 4: Comparison of simulation with approximations for an epidemic
in a clustered host population.  Parameters are as in Figure 1b.
\newline
\newline
Figure 5: Final size of epidemic ($\lim_{t \rightarrow \infty} (1 - S)$),
as predicted by mean field, MFC, and simulations.  Host is Poisson
distributed; pathogen dispersal is exponential with $A_{d} = 7.4$.
MFC predicts a critical transmission rate of 1.14.
\newline
\newline
Figure 6: Dependence of epidemic threshold on spatial factors for
clustered hosts. Both
host clustering and pathogen dispersal use exponential kernels. (a)  Curved
surface indicates critical transmission rate predicted by MFC model,
as a function of host clustering scale ($A_{h}$) and pathogen dispersal
scale ($A_{d}$).  Plane indicates critical transmission predicted by
mean field model.  (b) Cross--section of the threshold surface with
$A_{h} = 20$.
\newline
\newline
Figure 7: Effect of dispersal distance on simulated epidemics.  Hosts are
clustered, with $A_{h} = 20$ and $n_{c} = 10$.  ``Global'' indicates
uniform dispersal across the entire region.  Finite dispersal uses
exponential kernel.  Transmission rate is $\lambda = 0.8$.  Time series
are averaged from 8 simulation runs.
\newline
\newline
Figure 8: Effect of clustering distance on simulated epidemics.  Hosts
are either Poisson distributed or clustered using an exponential kernel
and $n_{c} = 5$.  Pathogen dispersal is exponential, with $A_{d} = 12$.
Transmission rate is $\lambda = 1$.
\newline
\newline
Figure 9: Dependence of epidemic threshold on spatial factors for
overdispersed hosts.  Pathogen dispersal uses exponential kernel.
Curved surface indicates critical transmission rate predicted by MFC model,
as a function of host inhibition distance ($a$) and pathogen dispersal
scale ($A_{d}$).  Plane indicates critical transmission predicted by
mean field model.
\newline
\newline
Figure 10: Effect of host inhibition distance on simulated epidemics.
Hosts are either Poisson distributed or overdispersed using minimum
inter--plant distance $a$. Pathogen dispersal is exponential, with
$A_{d} = 12$.  Transmission rate is $\lambda = 1.5$.
\newline
\newline
Figure 11: Dependence of epidemic threshold on pathogen dispersal kernel
type in a Poisson distributed host population.
\newline
\newline
Figure 12: Dependence of epidemic threshold on spatial factors and kernel
type in clustered hosts.  (a) Pathogen dispersal and host clustering use
normal kernels, with $n_{c} = 5$.  (b) Pathogen dispersal and host clustering
use fat--tailed kernels.
\newline
\newline
Figure 13: Effect of pathogen dispersal kernel type on simulated epidemics.
Hosts are Poisson distributed.  Dispersal scale is $A_{d} = 1$;
transmission rate is $\lambda = 4$.

\pagebreak

\section*{Figures}

Table 1:
\begin{table}[h]
\begin{tabular}{|llcc|} \hline
 & & & \\
Name    & Formula       & Effective & Mean Dispersal \\
 & & Area & Distance \\
Exponential & $\frac{m^{2}}{2 \pi} e^{-m r}$
 & $\frac{8 \pi}{m^2}$ & $\frac{2}{m}$ \\
 & & & \\
Normal & $\frac{m^{2}}{\pi} e^{-m^{2} r^{2}}$
 & $\frac{2 \pi}{m^2}$ & $\frac{\sqrt{\pi}}{2 m}$ \\
 & & & \\
Fat tailed & $\frac{m^{2}}{24 \pi} e^{-\sqrt{m r}}$
 & $\frac{384 \pi}{m^2}$ & $\frac{20}{m}$ \\
 & & & \\
\hline
\end{tabular}
\end{table}

\pagebreak

Figure 1a:

\begin{figure}[h]

        \scalebox{.5}{\includegraphics{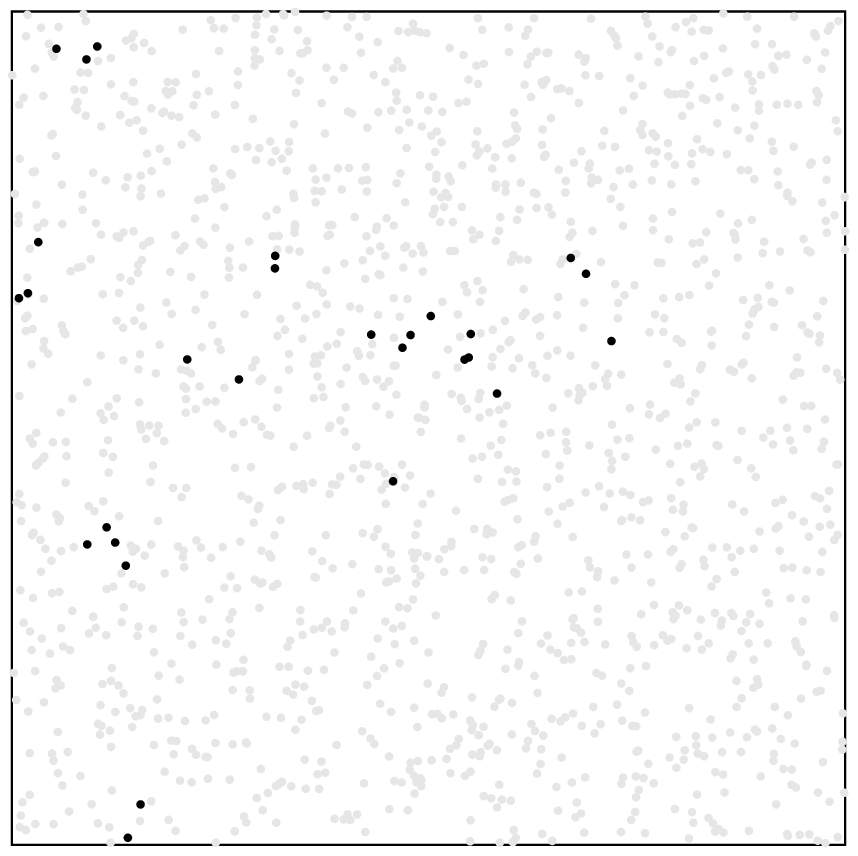}}

\end{figure}

\pagebreak

Figure 1b:

\begin{figure}[h]

        \scalebox{.5}{\includegraphics{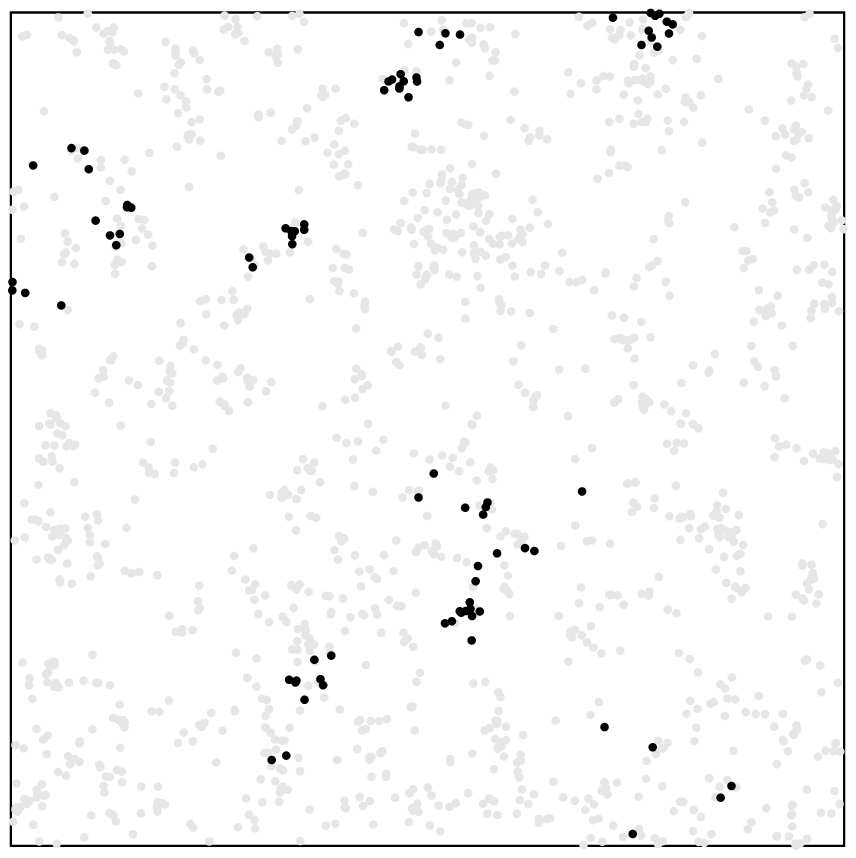}}

\end{figure}

\pagebreak

Figure 1c:
\begin{figure}[h]

        \scalebox{.5}{\includegraphics{threshfig1c.eps}}

\end{figure}

\pagebreak

Figure 2a:

\begin{figure}[h]

        \includegraphics{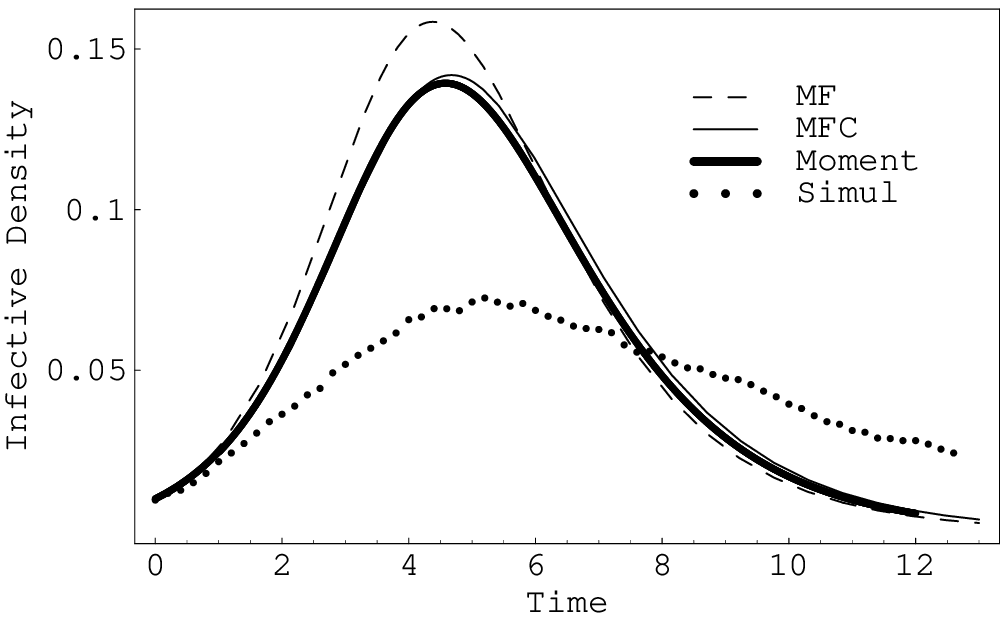}

\end{figure}

Figure 2b:

\begin{figure}[h]

        \includegraphics{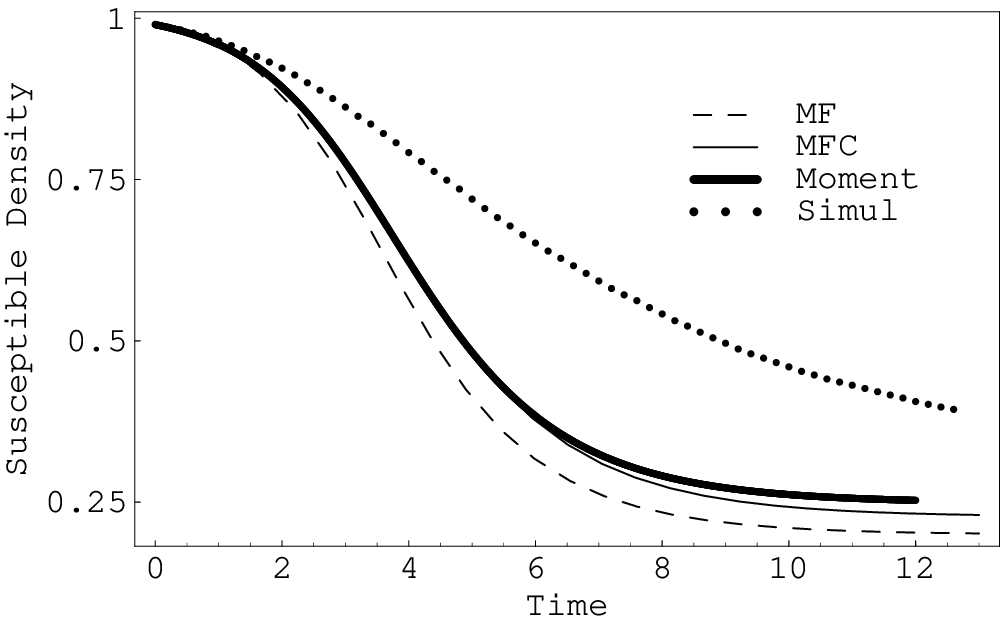}

\end{figure}

\pagebreak

Figure 2c:

\begin{figure}[h]

        \includegraphics{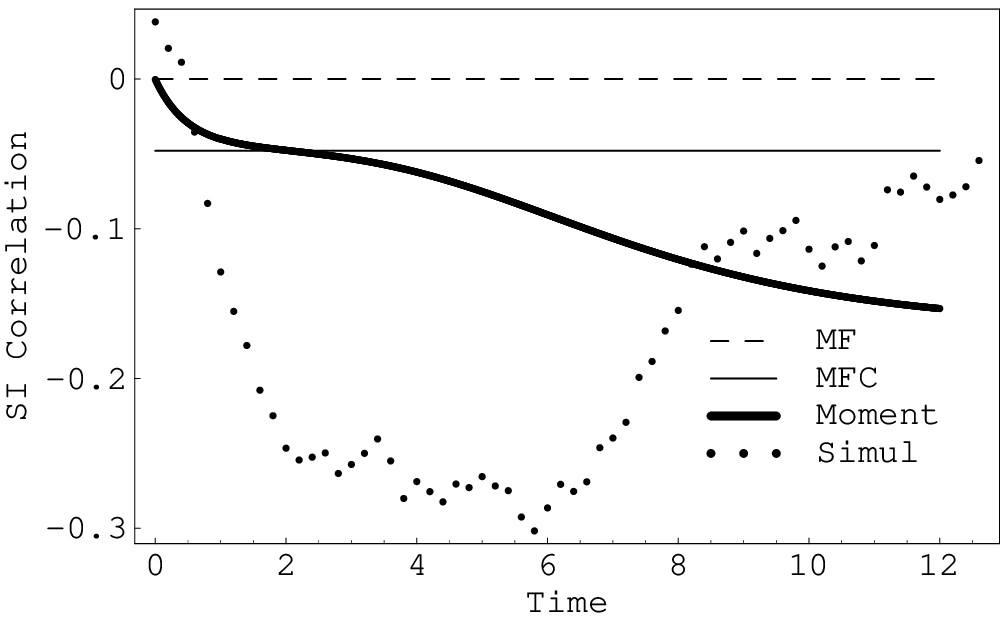}

\end{figure}

\pagebreak

Figure 3a:

\begin{figure}[h]

        \includegraphics{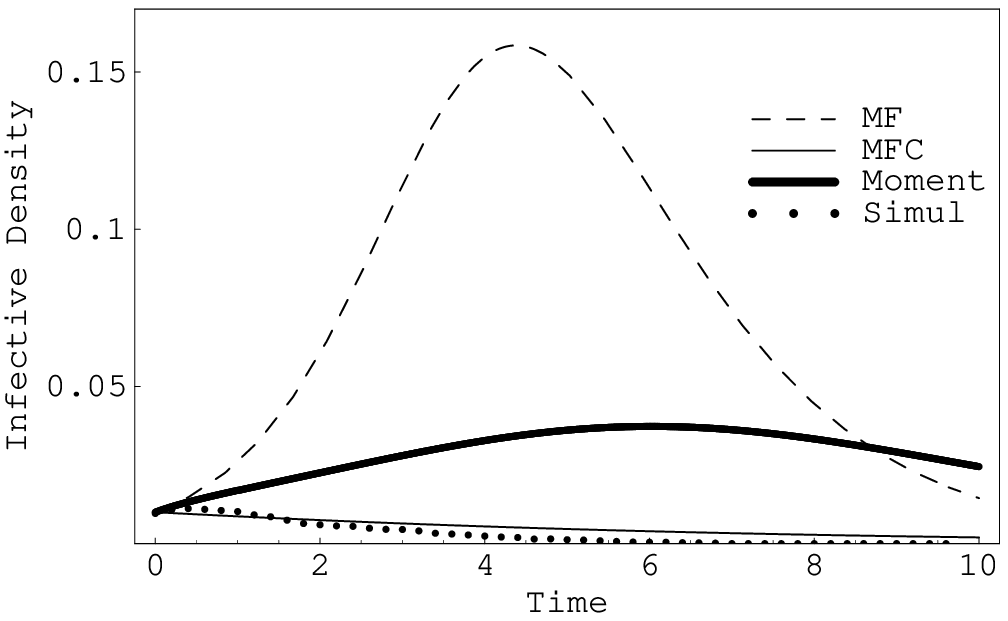}

\end{figure}

Figure 3b:

\begin{figure}[h]

        \includegraphics{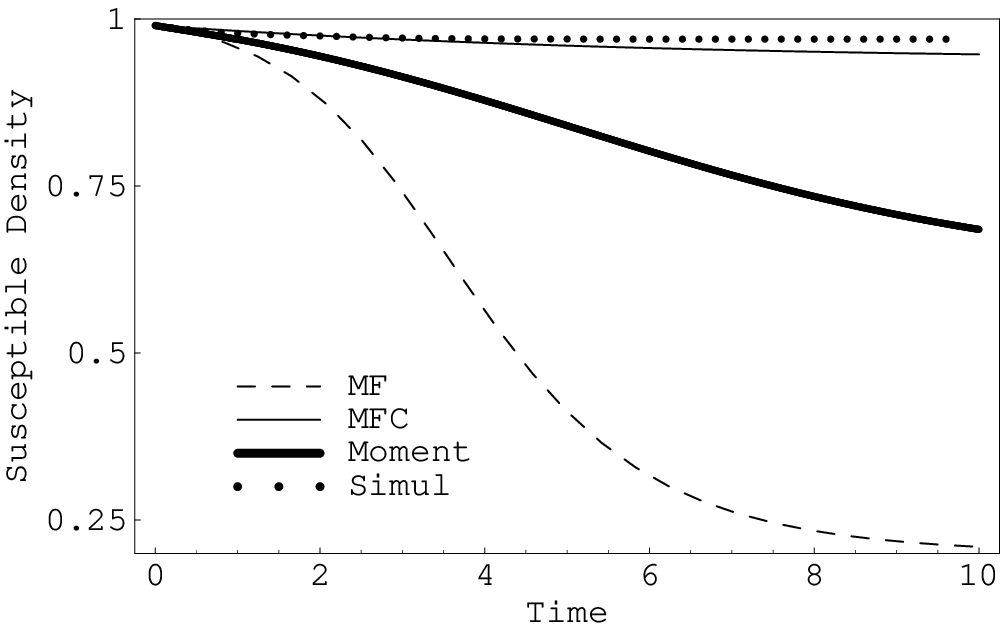}

\end{figure}

\pagebreak

Figure 3c:

\begin{figure}[h]

        \includegraphics{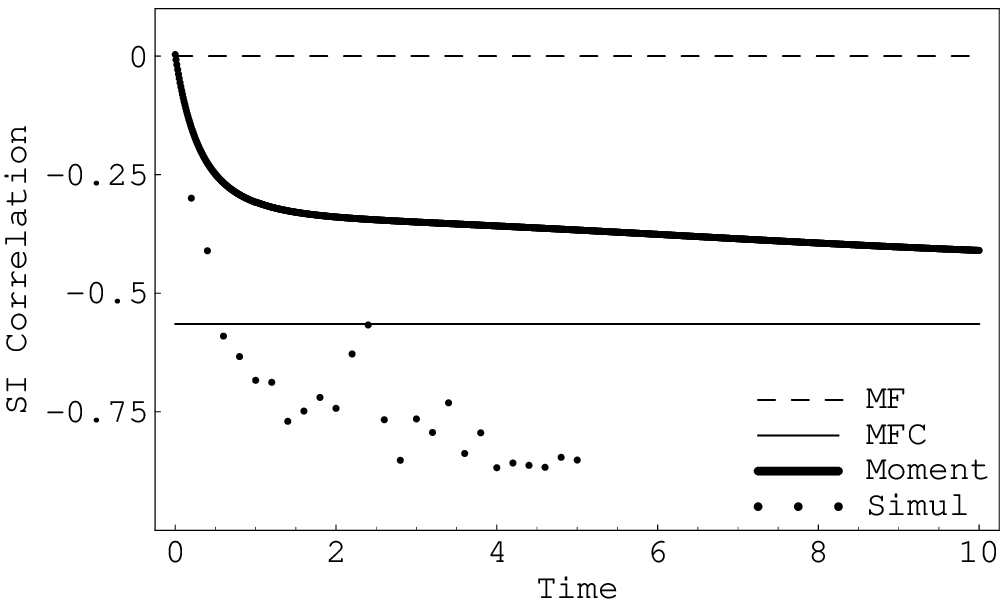}

\end{figure}

\pagebreak

Figure 4a:

\begin{figure}[h]

        \includegraphics{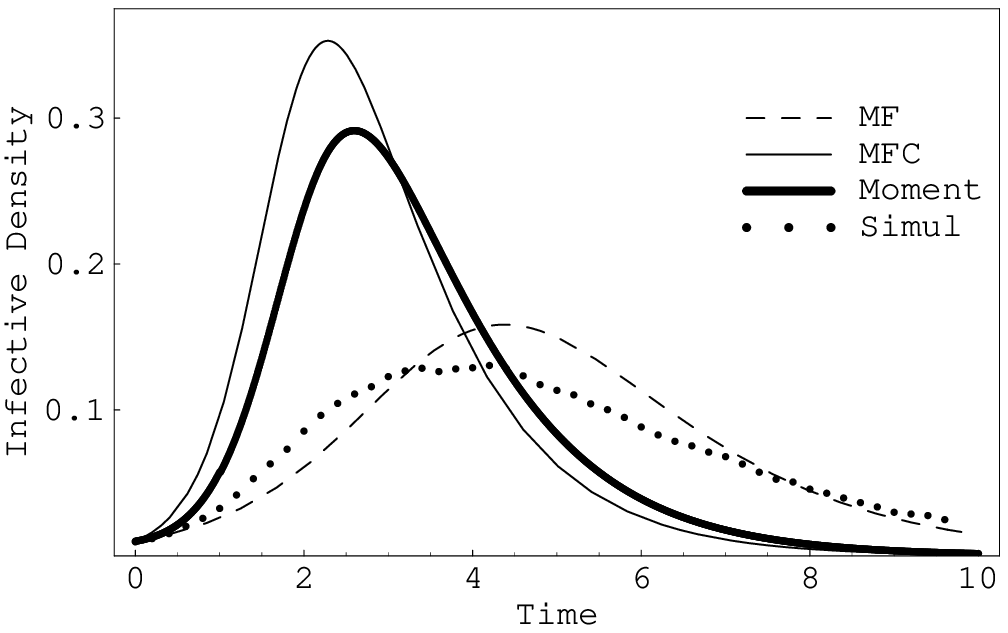}

\end{figure}

Figure 4b:

\begin{figure}[h]

        \includegraphics{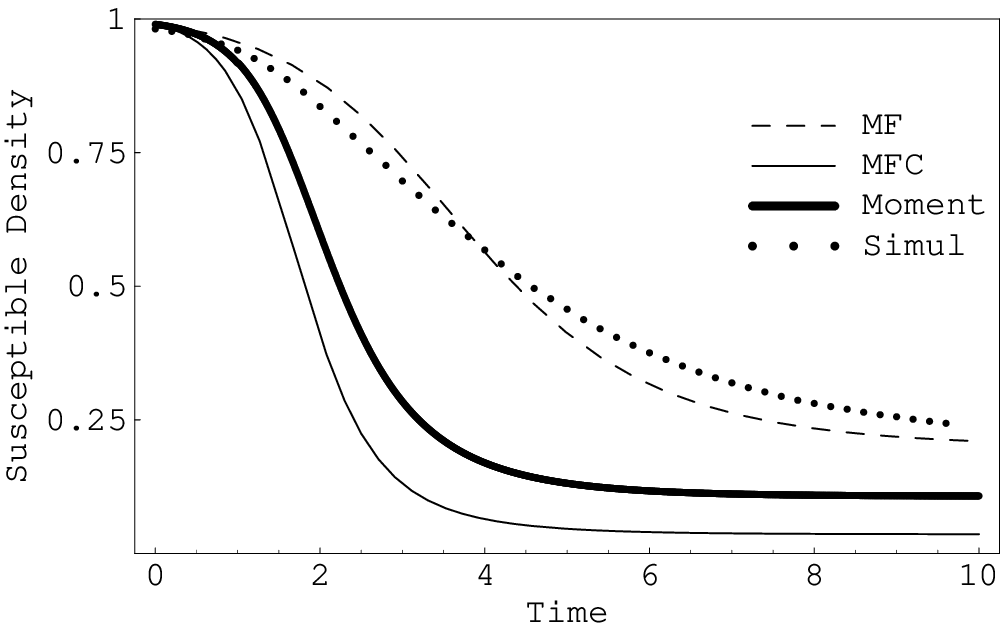}

\end{figure}

\pagebreak

Figure 4c:

\begin{figure}[h]

        \includegraphics{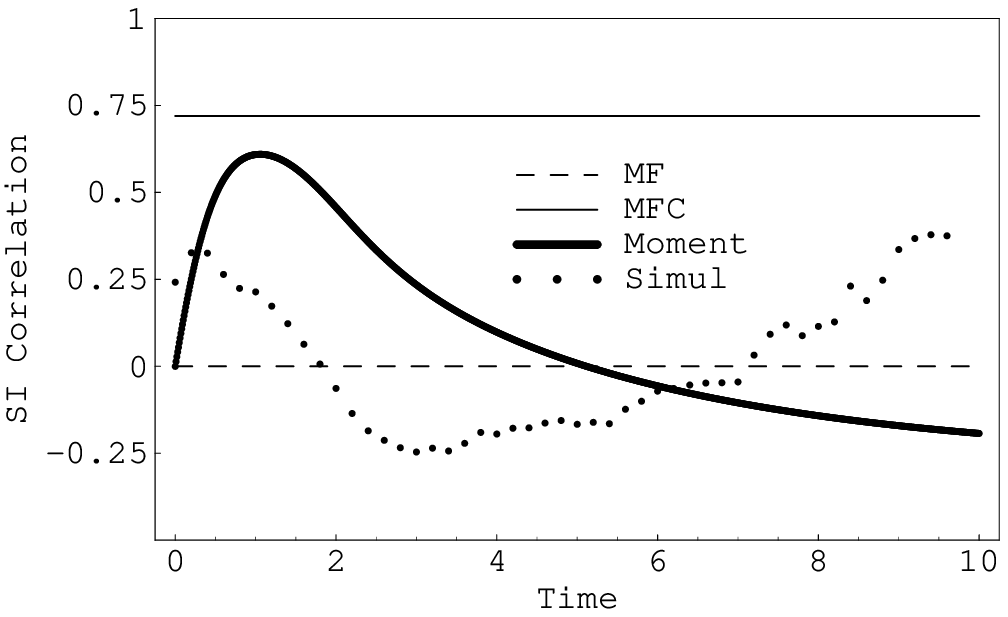}

\end{figure}

\pagebreak

Figure 5:

\begin{figure}[h]

        \includegraphics{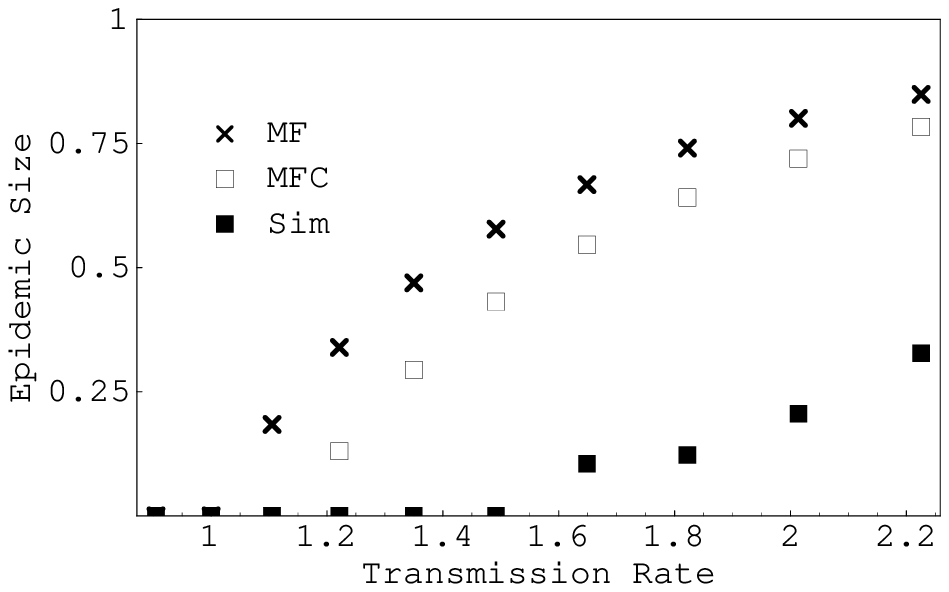}

\end{figure}

\pagebreak

Figure 6a:

\begin{figure}[h]

        \includegraphics{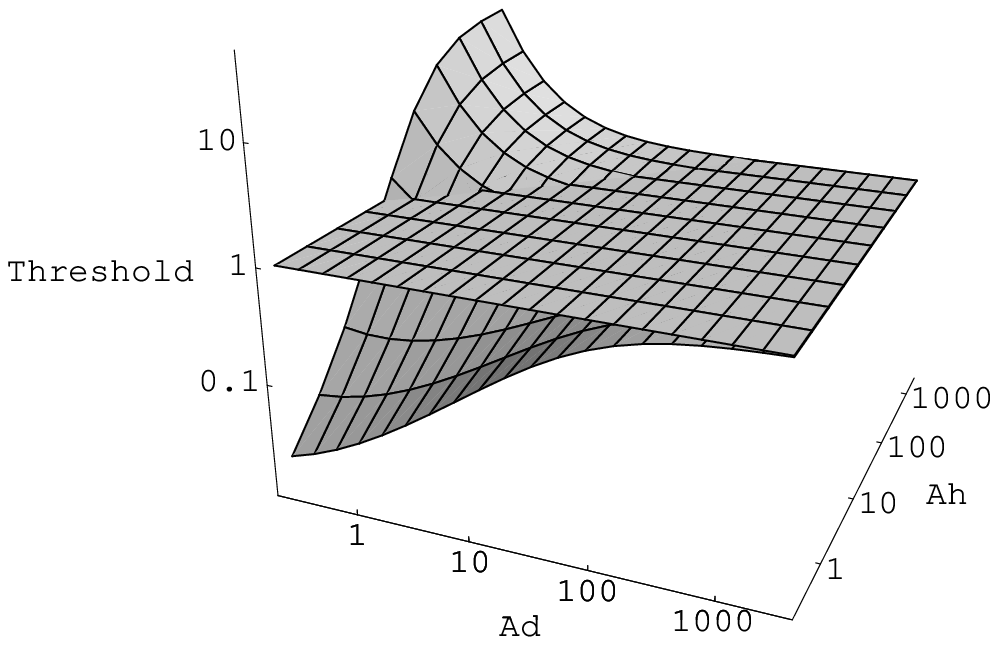}

\end{figure}

Figure 6b:

\begin{figure}[h]

        \includegraphics{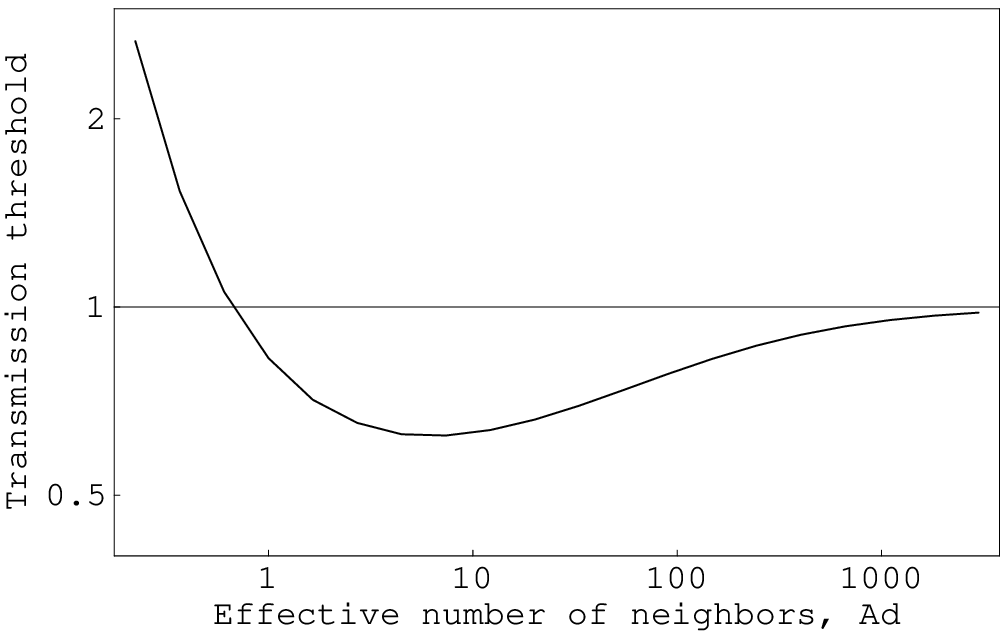}

\end{figure}

\pagebreak

Figure 7a:

\begin{figure}[h]

        \includegraphics{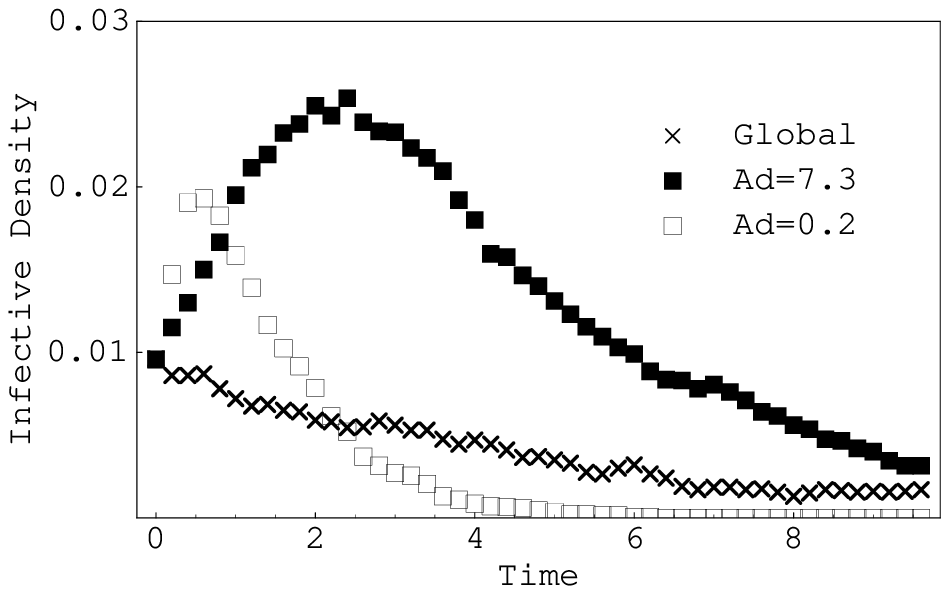}

\end{figure}

Figure 7b:

\begin{figure}[h]

        \includegraphics{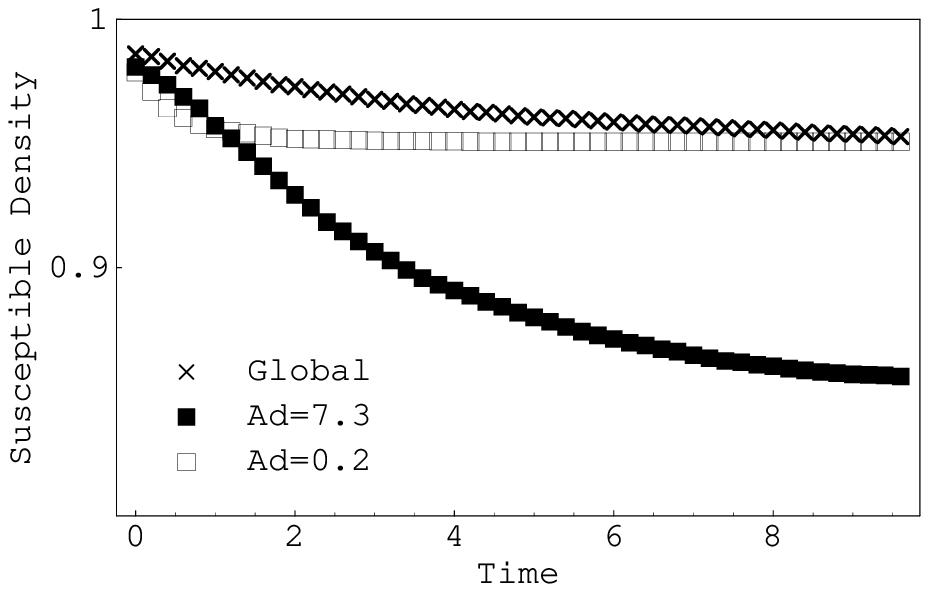}

\end{figure}

\pagebreak

Figure 8a:

\begin{figure}[h]

        \includegraphics{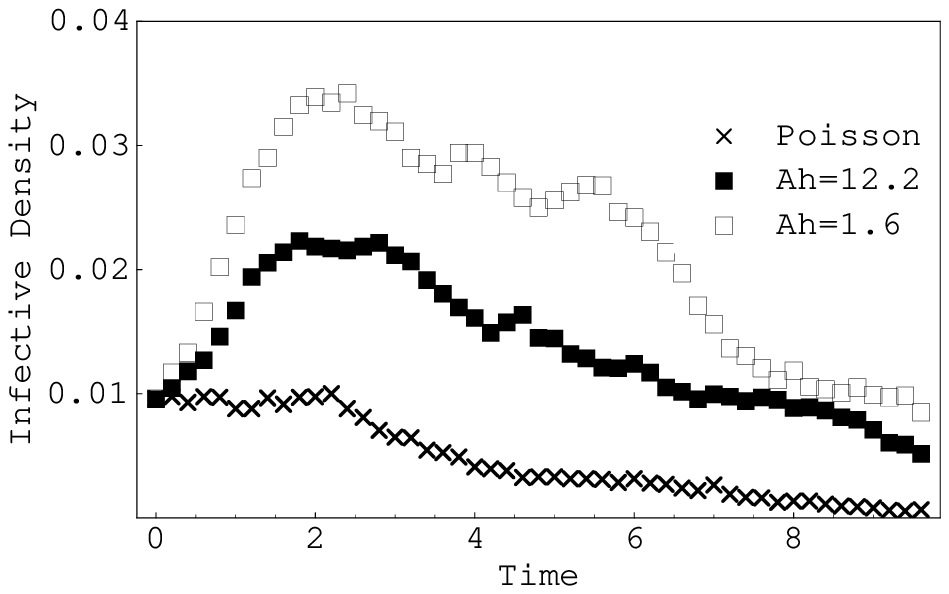}

\end{figure}

Figure 8b:

\begin{figure}[h]

        \includegraphics{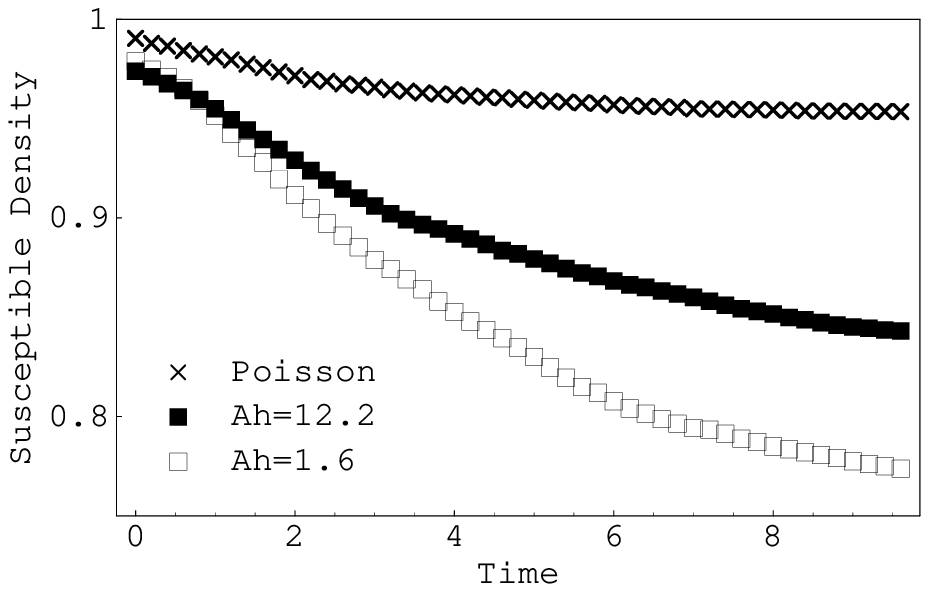}

\end{figure}

\pagebreak

Figure 9:

\begin{figure}[h]

        \includegraphics{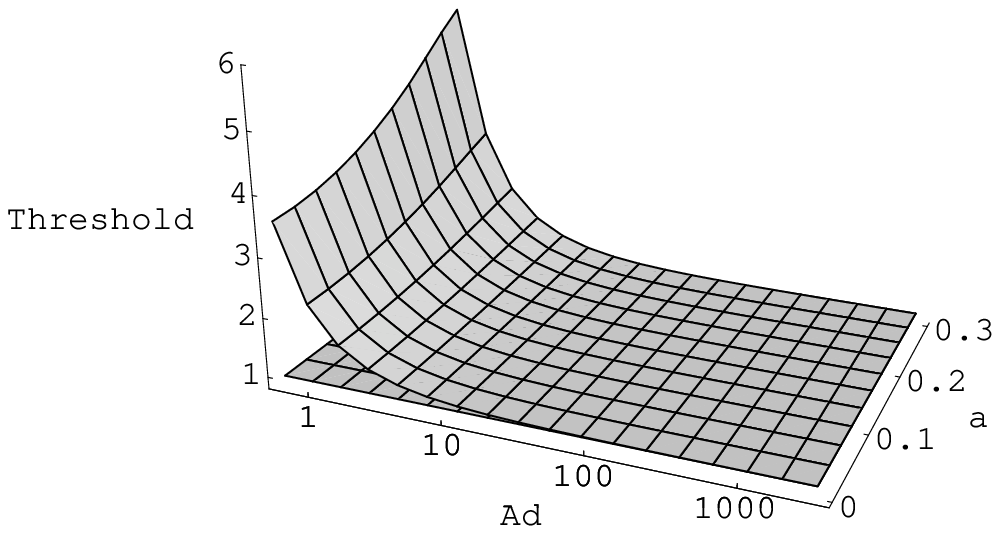}

\end{figure}

\pagebreak

Figure 10a:

\begin{figure}[h]

        \includegraphics{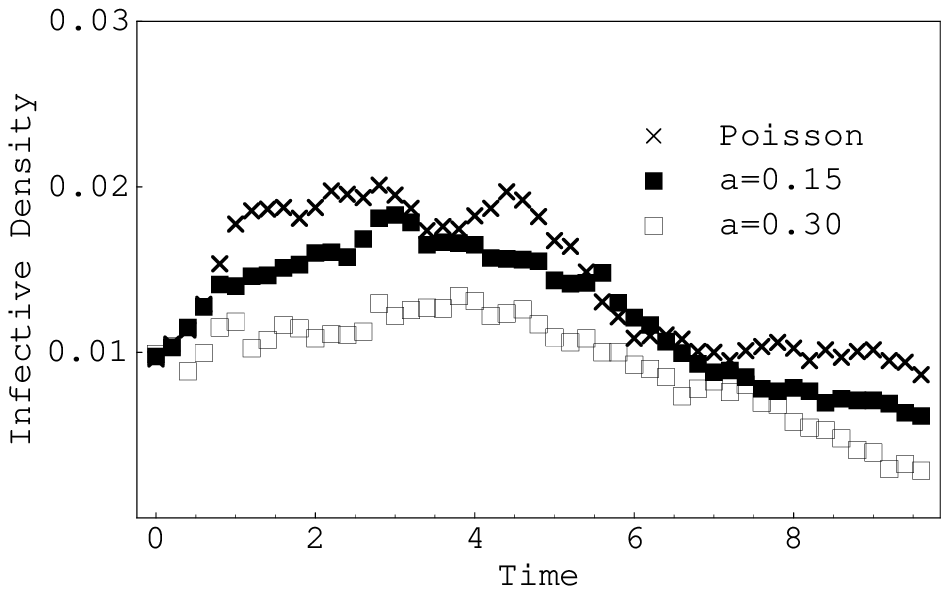}

\end{figure}

Figure 10b:

\begin{figure}[h]

        \includegraphics{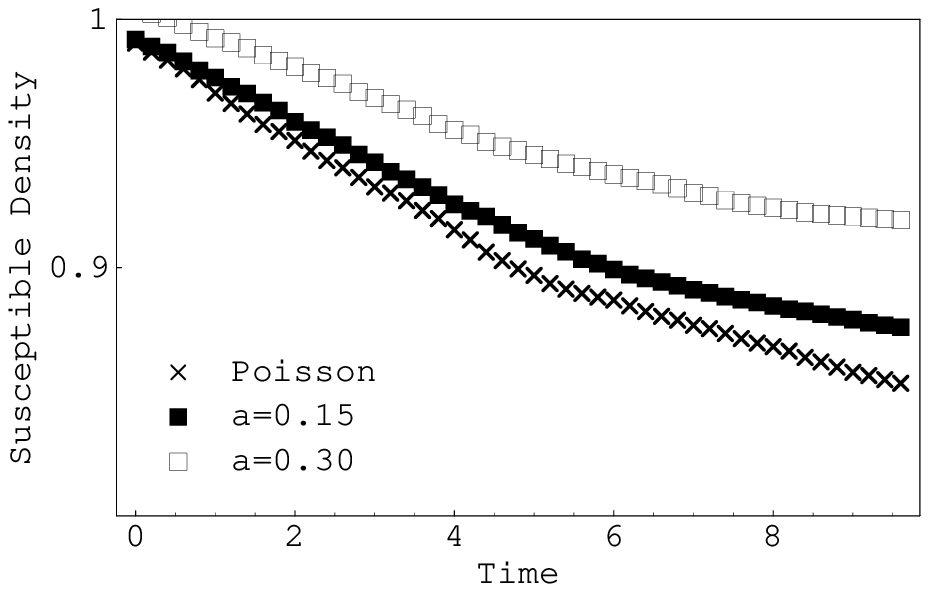}

\end{figure}

\pagebreak

Figure 11:

\begin{figure}[h]

        \includegraphics{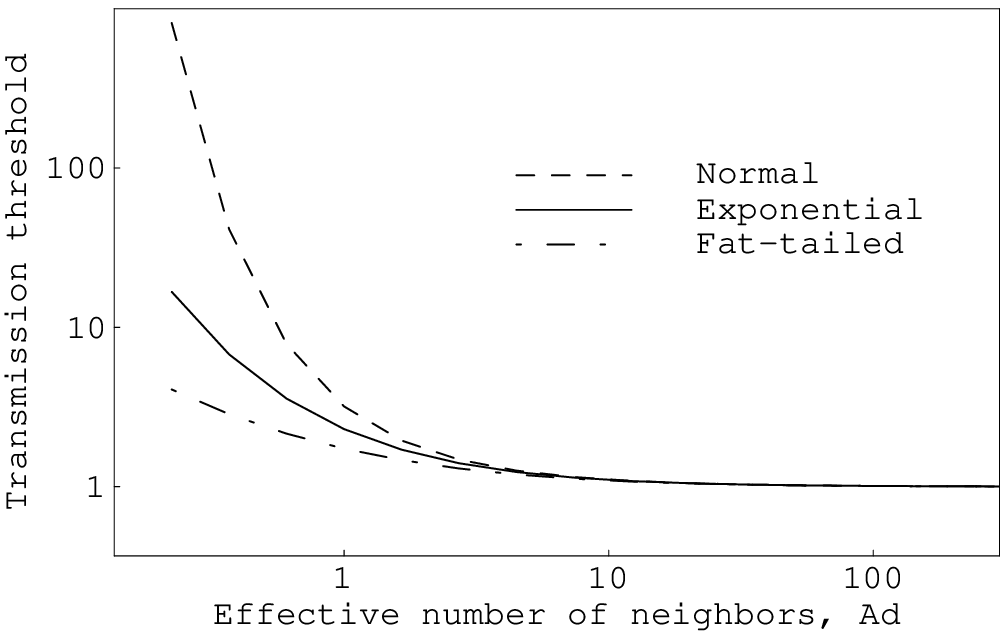}

\end{figure}

\pagebreak

Figure 12a:

\begin{figure}[h]

        \includegraphics{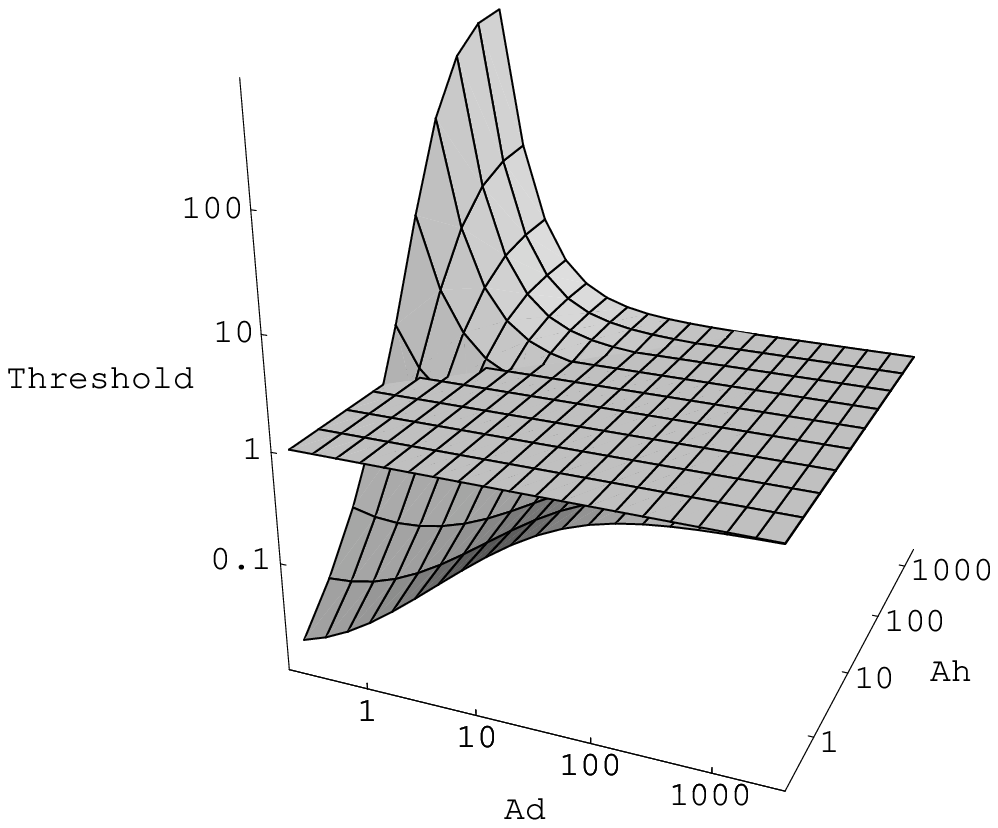}

\end{figure}

\pagebreak

Figure 12b:

\begin{figure}[h]

        \includegraphics{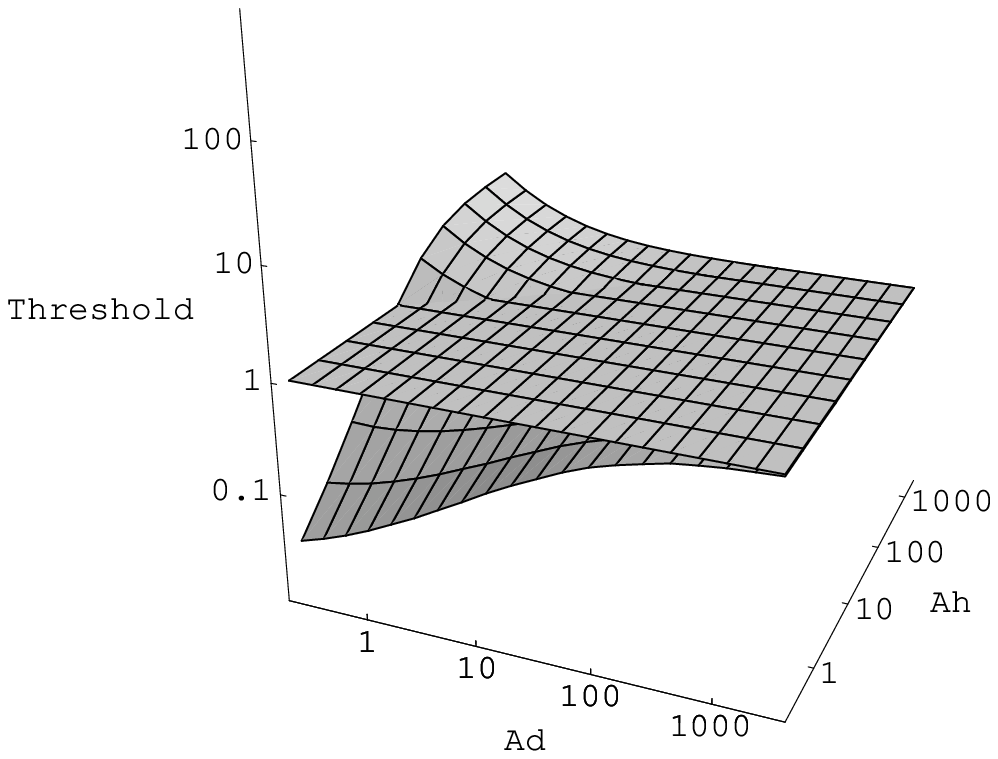}

\end{figure}

\pagebreak

Figure 13a:

\begin{figure}[h]

        \includegraphics{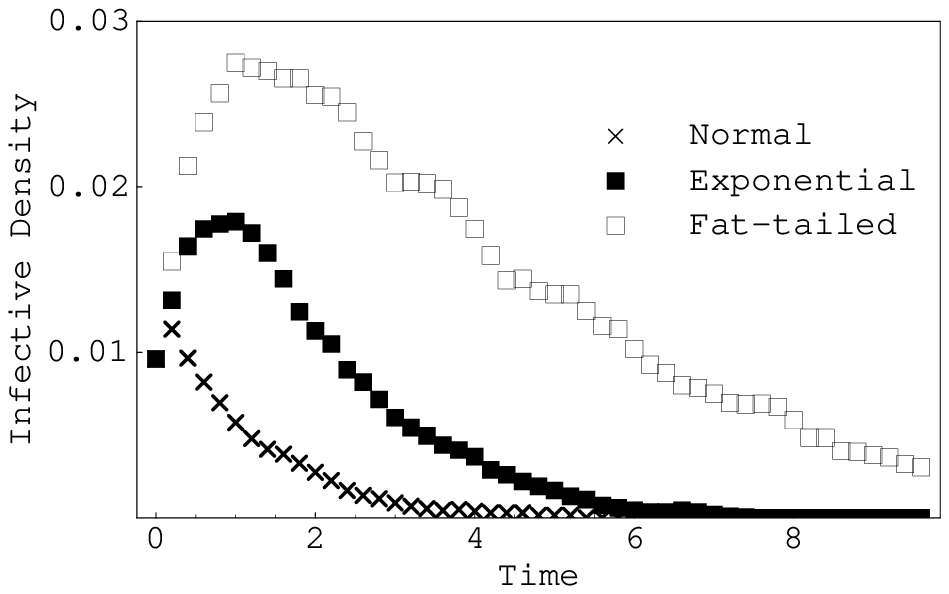}

\end{figure}

Figure 13b:

\begin{figure}[h]

        \includegraphics{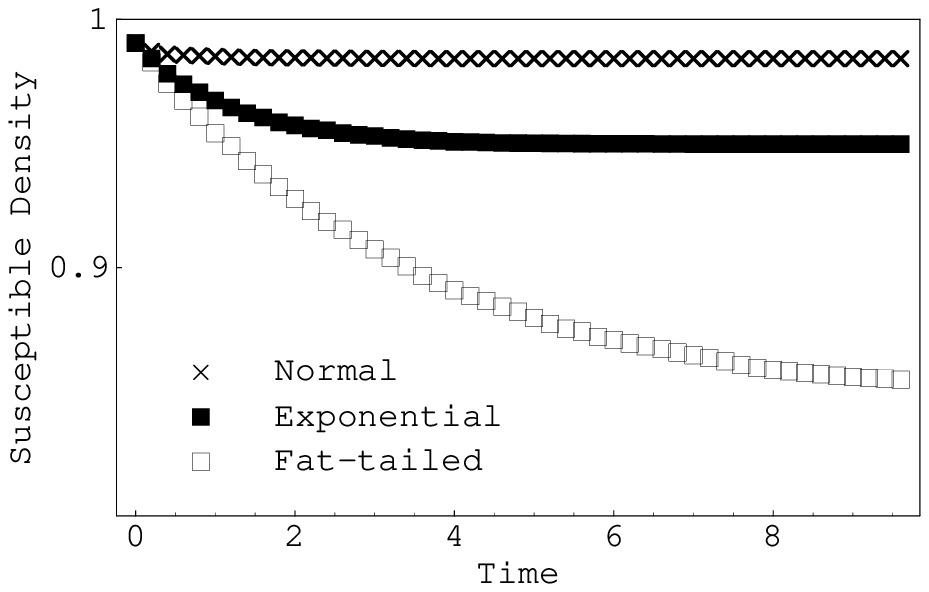}

\end{figure}

\chapter[Evolution of Resistance]{Evolution of Resistance}
\thispagestyle{myheadings}
\markright{}

\section*{Abstract}

	If a host species shares a pathogen with competing species,
the disease may provide a net benefit.  Selection for resistance will
depend on the trade--off between the damage done by the disease and
the positive effects resulting from infection of competitors.  This paper 
presents
a simple, spatially explicit model of a plant that shares a disease
with a superior competitor.  Pair approximations are used to determine
the phenotypic evolution of cost--free resistance.  Selection favors
lower resistance when transmission is spatially local and the damage
to the competitor is sufficient to outweigh the direct effects of 
infection.  This suggests that local spatial structure may be critical
in determining the coevolution of host--host--pathogen systems such
as heteroecious fungi. 

\pagebreak

\section*{Introduction}

	When competing species share a common enemy, the interaction
between direct competition for resources and apparent competition
mediated by the predator or pathogen can have important ecological
and evolutionary consequences 
(Park, 1948; Price et al., 1986, 1988; Begon and Bowers, 1995; Grosholz, 1992; 
Yan, 1996; Alexander and Holt, 1998).  
For example, it is not uncommon for directly
competing plant species to share a pathogen which can infect them both
(Rice and Westoby, 1982;  Clay, 1990).
If a pathogen infects a host's competitors, it may directly harm that host
while indirectly benefitting it by reducing competitive pressure.  
At the same time, coevolutionary forces acting on plants and their
pathogens are thought to be strong; pathogen--specific resistance to
infection is an important feature of plant evolution
(Burdon, 1987; Mitchell--Olds and Bergelson, 2000; Richter and Ronald, 2000).
The evolution of resistance may be determined by the tradeoff between
direct damage and indirect benefits that the pathogen confers to a host.
A natural question is: when will the benefits exceed the damage done
by the disease, thereby selecting for reduced resistance?
This paper addresses this issue by presenting a model which predicts
the evolution of disease resistance in a host that shares the pathogen
with a superior competitor. 

	A number of theoretical studies have addressed the evolution of
disease resistance in single host systems.  If there is no cost of
resistance, selection typically favors ever increasing resistance
(Boots and Haraguchi, 1998), although
Damgaard (1999) has shown that polymorphism for cost--free resistance may 
be maintained in a metapopulation.  Despite the apparent selective advantage 
of resistance, every species remains susceptible to infection by some
suite of pathogens.   
Two explanations for
this have been proposed.  First, resistance may come at a cost, such
as reduced fecundity or lifespan (reviewed in Purrington, 2000).  In that 
case, tradeoffs between the costs and
benefits of resistance may lead to selection for an intermediate resistance
level and/or polymorphism (Antonovics and Thrall, 1994; Bowers et al.,
1994; Boots and Haraguchi, 1999).  Second, selection for resistance may be
countered by selection for virulence in the pathogen, leading to a
coevolutionary ``arms race'' in which the host cannot achieve total resistance
(Dawkins and Krebs, 1979; Mitchell--Olds and Bergelson, 2000).
This hypothesis is supported by the gene--for--gene system in which 
resistance alleles confer protection against pathogens with specific
virulence alleles (Ellis et al., 2000; Stahl and Bishop, 2000;
Takken and Joosten, 2000).  
The interaction between these coevolutionary forces and
the costs of resistance may lead to a type of ``trench warfare''
(Stahl et al., 1999), in which 
resistance polymorphisms are maintained by spatial or temporal variation
in selection.

	Apparent competition between hosts that share a pathogen provides 
a third mechanism by which selection may not lead to complete resistance.
One host may evolve low resistance if the presence of the disease is a net
advantage because of the damage it does to a competitor.
Rice and Westoby (1982) argued that this mechanism must be invoked
to explain the existence of heteroecious rust fungi.  In order to 
complete their life cycle, heteroecious fungi must sequentially infect
two different host species (called the aecial and telial hosts, from the spore
stages that infect them).  While the fungi are highly
host specific, usually unable to infect species closely related to their
hosts, they are able to attack two different hosts which are often
taxonomically distant.  What allows them to do this?  A possible 
explanation is that the fungus has separate sets of genes responsible
for causing infections in the two host species, and that the two hosts'
resistance genes target different avirulence genes in the fungus.  
On the other hand, if
a single set of fungal genes controls infection of both hosts, and
both hosts target the same avirulence genes,
one must ask why the hosts have not developed resistance strategies
that the fungus cannot simultaneously overcome. 
In this case, the evolutionary stability of the system
can only come from the fact that one host species is not being selected
for resistance to the pathogen.  Rice and Westoby (1982) argued that
phylogenetic and ecological evidence supports this view, with the
telial host using the fungus as a weapon against a competitively
superior aecial host.  Whether or not this mechanism must be invoked
to explain the existence of heteroecious fungi depends on whether the
same genes are responsible for controlling infection and virulenc/avirulence
in both hosts.
This issue does not seem to have been resolved; however, differences
between the morphologies and infection mechanisms of the two spore stages
(Littlefield and Heath, 1979)
suggest that the two infection processes are under the control of 
separate genes
(Les Szabo, personal communication).  Moreover, there is no \emph{a priori}
reason to expect that both hosts would target the same avirulence genes.
Nevertheless,
selection for resistance in one host may be mediated by the impact
of the pathogen on the other host.  Even if apparent competition is
not a necessary explanation for the persistence of heteroecy, it may
still be an important factor in determining the evolution of resistance
in such systems. 

	This paper presents a model for the evolution of resistance
in a host plant that shares a pathogen with a superior competitor.  There
are three important features of the model.  First, it uses an evolutionarily
stable strategy (ESS) analysis of phenotypic evolution.  Although a great
deal is known about the genetics of disease resistance in plants, a phenotypic
approach allows us to focus on the structure of the ecological interactions
without incorporating specific assumptions about the genetics.  A number
of studies on the evolution of resistance have used a phenotypic approach to
generate robust predictions (Antononvics and Thrall, 1994; Bowers et al.,
1994; Boots and Haraguchi, 1998; Boots and Bowers, 1999).  Omitting population
genetics limits our ability to address the maintenance of polymorphisms
in the population (an important aspect of resistance), but it clarifies
the impact of the ecological interactions on evolutionary trends.  

	Second, the model does not incorporate costs of resistance.  
The previous theoretical studies of phenotypic evolution have
shown that the selected level of resistance can depend critically on the
costs associated with it, and there is extensive empirical evidence for 
costs of pathogen--specific resistance (Purrington, 2000).  By omitting costs
from the model, we do not imply that they do not exist.  Rather, we are
studying in isolation a different force that acts on the evolution of 
resistance: an ``ecological cost'' (Strauss et al., 1999) of resisting a 
pathogen that could serve
as an ally.  A fuller understanding of disease resistance will require
the incorporation of both direct physiological costs and indirect
ecological costs of resistance.  

	Third, the model is spatially explicit.
A number of studies have shown that limited pathogen dispersal can have 
important consequences for the evolution of resistance and virulence.
Large scale spatial structure can stabilize polymorphisms by decoupling
coevolutionary processes across the landscape 
(Burdon and Thrall, 1999; Damgaard, 1999; Stahl et al., 1999).  
At a smaller scale, spatial
structure can determine the selective pressures acting on individuals,
yielding ESS predictions that differ qualitatively from analogous
nonspatial models (Rand et al., 1995; van Baalen and Rand, 1998; Boots
and Sasaki, 2000).  It is this invidual--scale spatial structure that we
focus on here.  Rice and Westoby (1982) predicted that the dispersal scales
of the hosts and pathogen stages would be critical in determining the
evolutionary stability of the heteroecious fungi systems.  They argued
that the fungus can only serve as a useful weapon for the telial host
if spore dispersal from the telial to the aecial host is sufficiently
localized.  Hosts that accomodate the fungus with lower resistance must
reap the benefits of reduced competition; this will not occur unless
the pathogen primarily infects their neighbors.  Thus, it is expected
that small--scale spatial structure plays a key role in the evolution of 
resistance mediated by apparent competition.

	We use the model to address the following general question: what
level of resistance will be selected for in a plant that shares a pathogen
with a superior competitor?  Specifically, we examine how the ESS level
of resistance in the ``user'' host depends on the following factors:
\begin{enumerate}
	\item pathogen life history: we compare heteroecious (strictly 
		alternating) pathogens with those that are
		transmitted from either host to either host;
	\item dispersal scales of the pathogen and hosts: we compare
		local and global dispersal;
	\item increased mortality caused by infection;
	\item resistance by the other (``attacked'') host.
\end{enumerate}
The model is formulated as a stochastic, continuous--time process on a
lattice (an interacting particle system).  We analyze the model by using
pair approximations, which incorporate local, pairwise spatial structure into
a system of ordinary differential equations.	 
   
\section*{Model}
 
	Each site in a square lattice is assumed to be in one of several 
states: $E$ 
(empty), $S_{A}$ or $I_{A}$ (susceptible or infected attacked host), $S_{U}$
or $I_{U}$ (susceptible or infected user host).   In addition, we study
invasions by a new phenotype of the user host; we denote the susceptible
and infected invaders $S$ and $I$ for simplicity.  Each host and each pathogen 
stage either disperses locally (to the four nearest neighbors of a site)
or globally (uniformly across the entire system).  The attacked host
is assumed to be competitively dominant; it reproduces by placing offspring
on sites that are empty or occupied by the user host.  For simplicity,
we assume that the presence of a user host does not affect the probability
of establishment by an attacked host.  Conversely, user hosts reproduce
only onto empty sites; they cannot colonize sites already occupied.
Infection occurs when the pathogen is transmitted from an infected host
to an appropriate susceptible host, and resistance fails.  We assume
that resistance is a quantitative trait scaled to lie between $0$ 
(no resistance) and $1$ (total resistance).  We also assume that
the disease can increase the mortality rate of the hosts, but does not
affect reproduction.

	First consider the system when there is only one (``resident'') 
phenotype of the user 
host.  The probabilities that a randomly chosen site is in a particular state
satisfy a set of ordinary differential equations.  Let
$P_{\sigma}$ be the probability that a random site is in state $\sigma$,
and $P_{\sigma | \sigma'}$ be the conditional probability that, given a
site is in state $\sigma'$, a randomly chosen neighbor is in state $\sigma$.
Then by considering the possible changes to a state, and the rates at
which they occur, we can derive equations for the rates of change
of the various states.  For example, consider colonization of an
empty site by a susceptible attacked host.  Suppose each such host produces
locally dispersing offspring at rate $\beta_{A}$. Then any empty
neighbor of an $S_{A}$ site is colonized at rate $\frac{\beta_{A}}{4}$
(since there are four possible sites for the offspring to disperse to).
Since the total number of empty neighbors of $S_{A}$ sites is $4 P_{E|S_{A}}$,
this process increases the density of $S_{A}$ sites at the rate
$\beta_{A} P_{E|S_{A}} P_{S_{A}}$ (the first term in equation 1).
On the other hand, if the offspring disperse globally and are produced
at rate $B_{A}$, then colonization occurs at rate  $B_{A} P_{E} P_{S_{A}}$.
Incorporating all possible local and global interactions yields the
resident density equations:   
\begin{eqnarray}
\dot{P}_{S_{A}} & = & [\beta_{A} (P_{E|S_{A}} + P_{S_{U}|S_{A}} + 
	P_{I_{U}|S_{A}}) + B_{A} (P_{E} + P_{S_{U}} + P_{I_{U}})] P_{S_{A}}
	\nonumber \\
    & & + [\beta_{A} (P_{E|I_{A}} + P_{S_{U}|I_{A}} + P_{I_{U}|I_{A}}) 
    	+ B_{A} (P_{E} + P_{S_{U}} + P_{I_{U}})] P_{I_{A}}
	\nonumber \\
    & &	- (1 - \rho_{A}) 
	[\gamma_{AA} P_{I_{A}|S_{A}} + \Gamma_{AA} P_{I_{A}}  
    	+ \gamma_{UA} P_{I_{U}|S_{A}} + \Gamma_{UA} P_{I_{U}}] P_{S_{A}}
	- \mu_{A} P_{S_{A}} \\
\dot{P}_{I_{A}} & = & (1 - \rho_{A})
        [\gamma_{AA} P_{I_{A}|S_{A}} + \Gamma_{AA} P_{I_{A}}            
        + \gamma_{UA} P_{I_{U}|S_{A}} + \Gamma_{UA} P_{I_{U}}] P_{S_{A}}
	- \alpha_{A} \mu_{A} P_{I_{A}} \\
\dot{P}_{S_{U}} & = & [\beta_{U} P_{E|S_{U}} + B_{U} P_{E}] P_{S_{U}}
	+ [\beta_{U} P_{E|I_{U}} + B_{U} P_{E}] P_{I_{U}} 
	\nonumber \\
    & &	- [\beta_{A} (P_{S_{A}|S_{U}} + P_{I_{A}|S_{U}}) + 
	B_{A} (P_{S_{A}} + P_{I_{A}})] P_{S_{U}} 
	\nonumber \\
    & & - (1 - \rho_{U}) 
        [\gamma_{UU} P_{I_{U}|S_{U}} + \Gamma_{UU} P_{I_{U}}
        + \gamma_{AU} P_{I_{A}|S_{U}} + \Gamma_{AU} P_{I_{A}}] P_{S_{U}}
        - \mu_{U} P_{S_{U}} \\
\dot{P}_{I_{U}} & = & (1 - \rho_{U}) 
        [\gamma_{UU} P_{I_{U}|S_{U}} + \Gamma_{UU} P_{I_{U}}
        + \gamma_{AU} P_{I_{A}|S_{U}} + \Gamma_{AU} P_{I_{A}}] P_{S_{U}}
	\nonumber \\
    & & - [\beta_{A} (P_{S_{A}|I_{U}} + P_{I_{A}|I_{U}})
	+ B_{A} (P_{S_{A}} + P_{I_{A}})] P_{I_{U}} - \alpha_{U} \mu_{U} 
	P_{I_{U}} \\
P_{E} & = & 1 - P_{S_{A}} - P_{I_{A}} - P_{S_{U}} - P_{I_{U}}. \nonumber
\end{eqnarray}  
Here, $\mu_{A}$ and $\mu_{U}$ are the density--independent death rates of
uninfected attacked and user hosts, respectively; $\alpha_{A}$ and 
$\alpha_{U}$ are the factors by which infection increases host mortality -- 
we refer to this as the damage done by the pathogen.
The resistance levels of the hosts are given by $\rho_{A}$ and $\rho_{U}$.
The other parameters give rates of reproduction or transmission, and
hence describe interactions between two sites; lower case parameters
correspond to nearest--neighbor interactions, while upper case ones
corresond to global dispersal.  Thus, $\beta_{A}$ is the rate at which
the attacked host produces offspring which will disperse locally, while
$B_{A}$ is the rate at which it produces offspring which will disperse
globally. The parameters $\gamma_{\sigma \sigma'}$ and
$G_{\sigma \sigma'}$ give the transmission rates from a host of 
species $\sigma$ to one of species $\sigma'$.

	We assume throughout that the competitively dominant (attacked)
host is a perennial, while the user host is an annual.  This pattern is
often seen in the heteroecious fungi systems (Rice and Westoby, 1982).
Furthermore, we can rescale time by the lifespan of the organisms.
Thus, in all the examples presented we will set $\mu_{U} = 1$ and
arbitrarily choose
$\mu_{A} = 0.1$.  We assume that each host species or pathogen stage has
either completely local or completely global dispersal, so that for
each interaction either the upper case or lower case parameter is zero.
Finally, for heteroecious fungi we set $\gamma_{\sigma \sigma} = 
G_{\sigma \sigma} = 0$; i.e. transmission can only occur between different
species of hosts.  In the case of an arbitrarily shared pathogen,
we set $\gamma_{\sigma \sigma'} = \gamma_{\sigma' \sigma'}$, assuming
that the transmission rate depends only on the susceptible species.
In reality, the transmission rate will also depend on the rate of 
pathogen production by the infected species, but we do not include that
complication.

	If all interactions are global (the so--called ``mean field'' case),
the terms for pairs of sites drop out and the resident density equations
form a closed system.  Otherwise, the equations contain the unknown
conditional density terms $P_{\sigma|\sigma'}$.  With a little more
effort, we can write down ODEs for these terms; however, they include
terms involving triplets of sites.  We close the system at the level
of pairs by assuming that $P_{\sigma|\sigma' \sigma''} = P_{\sigma|\sigma'}$;
this is known as a pair approximation (Matsuda, 1992; Rand, 1999).
Here, $P_{\sigma|\sigma' \sigma''}$ is the probability that, given a
site is in state $\sigma'$ and a randomly chosen neighbor is in state
$\sigma''$, another randomly chosen neighbor will be in state $\sigma$.
Thus, the pair approximation assumes that the state of one neighbor of
a site does not depend on the states of the other neighbors.  The resulting
equations for the pair densities (Appendix 1) give an approximation
to the local spatial structure of the system.  

	Once we have determined the resident densities from the mean
field or pair approximation equations, we want to study the evolution
of resistance by determining the ability of a different phenotype of
the user host to invade.  We assume that the invading phenotype differs
from the resident only in its resistance, $\rho'$.  The densities of
the invading phenotype satisfy the following invasion equations:
\begin{eqnarray}
\dot{P}_{S} & = & [\beta_{U} P_{E|S} + B_{U} P_{E}] P_{S}
        + [\beta_{U} P_{E|I} + B_{U} P_{E}] P_{I}
	\nonumber \\
    & & - [\beta_{A} (P_{S_{A}|S} + P_{I_{A}|S}) +
        B_{A} (P_{S_{A}} + P_{I_{A}})] P_{S} - \mu_{U} P_{S}
        \nonumber \\
    & & - (1 - \rho')
        [\gamma_{UU} (P_{I_{U}|S} + P_{I|S}) + \Gamma_{UU} (P_{I_{U}}
	+ P_{I})
        + \gamma_{AU} P_{I_{A}|S} + \Gamma_{AU} P_{I_{A}}] P_{S} \\
\dot{P}_{I} & = & (1 - \rho')
	[\gamma_{UU} (P_{I_{U}|S} + P_{I|S}) + \Gamma_{UU} (P_{I_{U}}
	+ P_{I})
        + \gamma_{AU} P_{I_{A}|S} + \Gamma_{AU} P_{I_{A}}] P_{S}
        \nonumber \\
    & & - [\beta_{A} (P_{S_{A}|I} + P_{I_{A}|I})
        + B_{A} (P_{S_{A}} + P_{I_{A}})] P_{I} - \alpha_{U} \mu_{U}
        P_{I}. 
\end{eqnarray}
The success or failure of the invasion is determined by the dominant
eigenvalue ($\lambda^{*}(\rho')$) of the invasion equations linearized
about $P_{S} = P_{I} = 0$, and with the resident densities fixed at their
equilibrium values obtained previously.  When $\lambda^{*}(\rho') > 0$, 
the invasion succeeds
and selection favors the new phenotype.  Of course, $\lambda^{*}(\rho_{U})
 = 0$, since selection is neutral when the resident phenotype tries to
invade itself.  

	The invasion equations include terms $P_{\sigma|S}$ and 
$P_{\sigma|I}$ that give the neighborhood structure of the invading
population.  In general, we expect the local structure of the invading
phenotype to differ from that of the resident.  Thus, we need to 
determine these conditional probabilities.  Again, we can use pair
approximation to write down ODEs for these terms (Appendix 2).  As is
commonly observed in this type of model (Matsuda, 1992; Brown, 2001), we 
find that the conditional
probabilities for the invader go to equilibrium on a much faster timescale
than the overall densities.  This occurs because local interactions 
allow structure to develop at the local scale faster than at the global
scale, especially when the invasion is developing slowly because of
small phenotypic differences.  Thus, we
can solve the pair equations for equilibrium values of the conditional
probabilities with the low density assumption $P_{S} = P_{I} = 0$;
then we incorporate this spatial structure into the invasion equations
as fixed parameters.  

	To determine the direction of phenotypic evolution, we use
the following steps: compute the resident densities using equations
1--4 and Appendix 1; compute the local spatial structure of the invading
phenotype using Appendix 2; find the largest eigenvalue of the invasion
equations (5--6).  Suppose $\rho^{-} < \rho_{U} < \rho^{+}$ are 
phenotypes that differ only slightly.  Then typically we find that
the invasion eigenvalues satisfy
$\lambda^{*}(\rho^{-}) < 0 < \lambda^{*}(\rho^{+})$ or
$\lambda^{*}(\rho^{-}) > 0 > \lambda^{*}(\rho^{+})$.  In the first case,
evolution leads to higher resistance; in the second case, it leads
to lower resistance.  An evolutionarily stable state is characterized
by $\lambda^{*}(\rho^{-}) < 0$ and $\lambda^{*}(\rho^{+}) < 0$, i.e.
the resident phenotype cannot be invaded.  The direction and relative
rate of evolution can be summarized in a term called the evolutionary
flux (Rand et al., 1994; Rand et al., 1995):
\begin{equation}
\frac{\lambda^{*}(\rho^{+}) - \lambda^{*}(\rho^{-})}{\rho^{+} - \rho^{-}}.
\nonumber
\end{equation}
In the limit of arbitrarily small differences in phenotypes, the flux
can be thought of as describing the slope of the fitness landscape.
Thus, a positive flux indicates that resistance will increase, while a
negative flux indicates that resistance will decrease.  At an ESS,
the flux is zero, passing from positive to negative as we increase
resistance. 
 
\section*{Results}

	The behavior of the model depends in a complex way on the full set
of parameters.  Coexistence of both species and endemicity of the pathogen
is only possible if the birth and transmission rates are sufficiently
high.  In the heteroecious case, there exists a threshold resistance
level for each host, above which the pathogen cannot persist because
successful infection is too rare.  When the pathogen is transmitted
arbitrarily between host species, it may or may not be able to persist
in spite of total resistance by one species, depending on the other
parameters.  
	
	In order to allow comparisons between different dispersal
scales and pathogen life histories, we chose parameter values that
yielded the same resident densities.  As stated above, we chose death
rates so that the user host is shorter--lived than the attacked host.
We then chose birth rates so that, in the absence of the disease,
the equilibrium densities of the attacked and user hosts were approximately
$0.55$ and $0.2$, respectively.  With no resistance and moderate levels
of disease damage, we chose the transmission rates so that the disease was
endemic, yielding approximate equilibrium densities of $(P_{S_{A}}, 
P_{I_{A}}, P_{S_{U}}, P_{I_{U}}) = (.1,.1,.5,.1)$.  The presence of
the disease thus reversed the relative abundance of the two hosts.
For the disease damage values, we used $\alpha_{U} = 1.1$ always, while
$\alpha_{A}$ varied from $2.0$ to $2.4$ as needed to match the equilibrium
densities.  Notice that at these damage levels, the pathogen
increases the user host's mortality by $10\%$ and the attacked host's
mortality by $100\%$ or more.  However, because of the difference in
background mortality rates between the two species, this yields an
additive increase of $0.1$ to the mortality of each.  The disease
removes individuals of both species at approximately the same rate, but
the higher birth rate of the user host allows it overcome the added
mortality and exploit the space cleared by it.   
 
	For these parameter values, the presence of the disease is
advantageous to the user host.  As a result, an increase in the level
of resistance by the user host causes a decrease in its equilibrium
density.  The dependence of the equilibrium population levels on
the user's resistance is shown in Figure 1 for the mean field heteroecious
case; other cases behave similarly.  Notice that there is a threshold
level of resistance at approximately $\rho_{U} = 0.6$, above which the 
disease vanishes
and the resistance level has no effect.  Typically, an increase in
resistance by one host leads to lower infected populations of both hosts.
However, we note in passing the counterintuitive result that under
some circumstances, increasing $\rho_{U}$ from zero initially increases
the equilibrium level of $P_{I_{A}}$.  

	From the effect of user resistance on the equilibrium population
levels, we might infer that under these conditions evolution will always
lead to lower resistance by the user host.  However, an ESS approach 
indicates that this is not the case.  In the mean field case,
selection always leads to higher resistance because a user phenotype with
a higher resistance can always invade the system; a proof of this is
given in Appendix 3.  Thus, selection will lead to a user host that
maximizes resistance, even though this causes its population to shrink.
The intuitive reason for this was discussed in Rice and Westoby (1982):
less resistant user hosts bear the price of increase mortality from 
infection, but they do not reap any greater reward that more resistant
users.  When transmission from the user to the attacked host occurs
over long distances, there is nothing to prevent a higher--resistance
phenotype from ``cheating''; the damage done to the attacked hosts
benefits all user hosts equally, regardless of their level of resistance.
Because more resistant user hosts enjoy all of the benefits but
bear less of the burdens of the disease, they always displace less
resistant phenotypes.  We can only prove this result in the mean field
case, but we conjecture that it always holds when transmission
of the pathogen from the user to the attacked host is global.

	When this transmission is local, ESS analysis indicates that
selection for reduced user resistance can occur.   Figure 2a shows the
evolutionary flux in user resistance as predicted by the pair approximation
equations, in the heteroecious case with local user dispersal and
user--attacked transmission and global attacked host dispersal and
attacked--user transmission.  When the damage done to the user host 
($\alpha_{U}$) is sufficiently small, selection always favors lower
resistance.  As the damage increases, selection leads to an intermediate
ESS: below this value resistance increases, while above it resistance
decreases.  Finally, for sufficiently high damage levels, the direct
effect of infection outweighs the benefits and selection leads to
a resistance level that drives the disease extinct (here, around
$\rho_{U} = 0.45$).  Thus, the pair
approximation equations tell us that when dispersal is local, less
resistant phenotypes can sometimes displace more resistant ones.  Unlike
the mean field case, the pair approximation equations detect the fact
that the different user phenotypes encounter different competitive
pressures.  Less resistant phenotypes are more likely to infect their
neighbors, clearing space for growth.  

	A comparison between the predictions of the pair approximation
equations and simulations of the full stochastic model is shown in Figure 2b.
The simulations confirm that selection for lower resistance can occur.
In fact, to the extent that pair approximation predicts selection for
lower resistance, simulations indicate that it is conservative: lower
resistance is favored over a wider set of parameters than predicted.  
This occurs because pair approximation underestimates the clustering of
hosts.  When dispersal is local, the presence of a conspecific nearby
greatly increases the probability of finding another one, so that
$P_{\sigma|\sigma'\sigma} >> P_{\sigma|\sigma'}$ when $\sigma$ is rare.
The strategy of lowered resistance takes advantage of the local buildup of
infections, so that selection should be stronger than predicted by
pair approximations.  Corrections to the pair approximation to
deal with this phenomenon have been developed (Sato et al., 1994), 
though they increase
the complexity of the analysis.  We retain the usual pair approximation
assumption, since we believe it correctly predicts the qualitative
behavior of this evolutionary process. 

	By studying the evolutionary fluxes, we can determine how the
ESS level of resistance depends on various factors.  The damage that
the disease does to each host determines the level of resistance that
selection favors (Figure 3).   As the pathogen becomes more
lethal to the user host, the ESS resistance moves from complete
nonresistance to the threshold level at which the disease dies out.     
When the pathogen becomes more damaging to either host, this 
threshold level of resistance decreases as the infectious period
shrinks.  Still, as the disease becomes more lethal to the attacked
host, it becomes a more effective weapon, and selection favors lower 
resistance by the user host.   The more damaging the disease is to
the attacked host, the greater damage the user host can tolerate and
still gain by lowering resistance. 

	A change in the resistance level of the attacked host also 
effects the ESS level of resistance by the user host (Figure 4).  
If the attacked
host becomes more resistant, selection favors lower resistance by the
user host.  This suggests that the two host may engage in a kind of
``arms race by proxy'', with the user host lowering its own resistance
to offset increased resistance by the other species.  However, 
understanding this issue
would require a broader ESS approach, in which we examined selection for
each species' resistance and for the pathogen's virulence.  Notice that
although the user's ESS resistance decreases as the attacked host's 
resistance increases, the maximum level of damage that the user can
sustain and still favor reduced resistance decreases.  Thus, even as
the user is selected to offset the other species' resistance by lowering
its own, the range over which this strategy works is shrinking.  Although
increased damage and increased attacked host resistance both lead
to lower pathogen populations, they have different implications for the 
utility of the disease as a weapon by the user, and hence lead to 
different evolutionary trajectories for user resistance.

	As we discussed above, the dispersal scale of the pathogen
transmission from the user to the attacked host is a critical factor; 
selection only favors reduced resistance when this dispersal is local.
The other dispersal scales in the system may also be important, since
they determine the local spatial structure encountered by the user host.
Figure 5 shows the effect of varying the dispersal scales of the
attacked host and of attacked--user transmission while the other
dispersal scales are kept local.   Precise comparisons between these
cases are not meaningful, since parameters were chosen to match
resident densities approximately at an arbitrary point ($\rho_{U} = 0$, 
$\alpha_{U} = 1.1$).  However, it appears that the evolutionary advantage
of low user resistance increases when more dispersal scales are local.
Local dispersal increases spatial heterogeneity in the populations,
strengthening selective forces that depend on local spatial structure.  

	Our pair approximation approach fails when user hosts disperse
globally.  If other processes are local, this should still lead to
local spatial structure for the invader that differs from mean field
predictions.  However, the spatial invader equations (Appendix 2)
break down in this case.  These equations involve the ratio 
$\frac{P_{S|I}}{P_{I|S}}$, but when user dispersal is global, both
these quantities are zero in the low density limit.  Thus, we cannot
use the pair approximation equations to predict the evolution of resistance
when the user host disperses globally.  We conjecture that selection can
still favor lower resistance, but that it will act weakly.  Users that
disperse globally can still benefit from infecting neighbors, since these
neighbors threaten to displace them.  However, in our examples the dynamics
of the attacked host are much slower than those of the user host.  The
perennial attacked hosts affect the users primarily by acting as a barrier
to growth, rather than by displacing established users.  Thus, the
users can gain more by clearing space for colonization than by removing
neighbors that may displace them.  As a result, the amount of damage
that a user can sustain and still favor reduced resistance should be much
lower when it disperses globally than when it disperses locally.  Of course,
this would change if the attacked host had faster dynamics, threatening
the user with rapid displacement.  Analysis of this case might be feasible
by determining $\frac{P_{S|I}}{P_{I|S}} = \frac{P_{S}}{P_{I}}$ from the
eigenvector associated with the dominant eigenvalue (van Baalen and Rand,
1998), but we have not yet attempted this.

	The pair approximation approach also fails in the case that
transmission from user to attacked hosts is global, but all other processes
are local.  As we have argued, we do not expect selection to favor
reduced resistance in this case because infected users do not reap the
rewards of transmission.  However, the pair approximation equations
predict very weak selection for nonresistance when $\alpha_{U} = 0$,
with the ESS switching to complete resistance very quickly as this
damage increases.  The equations still satisfy the condition for
consistency between the invader and resident phenotypes, 
$\lambda^{*}(\rho_{U}) = 0$, so this does not result from the low density
assumptions.  Rather, the pair approximation errs in predicting that
the local spatial structure of the user depends on its resistance level
even when transmission to the attacked host is global.  This prediction
is not supported by simulations, and is an anomaly of the pair approximation
approach that we do not fully understand.  In other dispersal cases,
as we mentioned, the pair approximation's prediction of selection for
lowered resistance is conservative when compared to simulations.   

	Finally, we consider the importance of the pathogen life
cycle in determining the evolution of the user's resistance.  Thus far,
all examples have been based on a heteroecious life cycle, in which 
the pathogen alternates strictly between host species.  Figure 6
compares the heteroecious ESS resistance with that obtained when the 
disease can be transmitted between any hosts.   Recall that the parameters
were chosen so that the resident densities match.  Transmission between
arbitrary hosts allows much easier spread of the disease, so that the 
transmission levels used in this case are much lower than those used
in the corresponding heteroecious case.  As a result, infected user
hosts are less likely to infect neighboring attacked hosts when we use
this basis of comparison.  Consequently, the disease is a less potent
weapon, and selection leads to higher levels of user resistance than in 
the heteroecious case.  Of course, lowered resistance becomes a better
strategy when transmission from user to attacked hosts increases.  In
fact, for the parameter values used here, transmission between arbitrary
hosts at the same rate as used in the heteroecious case ($\lambda_{\sigma
\sigma'} = 10$) leads to the extinction of the attacked host.  
Depending on the transmission rates and relative damage, the disease
may be supported by either host at high enough levels to wipe out the
other species.  The dependence of the ESS on other parameters is
qualitatively similar to the heteroecious case: increasing the
damage done to the attacked host or that host's resistance leads to
lower resistance by the user, and the behavior is not sensitive to
the dispersal scale of the attacked host.  With transmission between
arbitrary hosts, we assume that the pathogen has a single dispersal
scale; selection for reduced user resistance is possible only if
this dispersal is local.

	In all of the cases studied, a successful invasion by a novel
phenotype leads to the extinction of the resident phenotype.  In addition,
for each set of parameters we found only one resistance phenotype stable 
against invasions by neighboring phenotypes.  Thus, the pair approximation
equations indicate that polymorphism for resistance will not be maintained
by the mechanism of apparent competition.  Simulations of the stochastic
model support this conclusion, although the time required for one 
phenotype to replace another is very long (typically, tens of thousands
of generations) because the invasion eigenvalues are so small.  Thus,
while this mechanism does not itself lead to polymorphisms, the fact
that it acts weakly and locally would allow other mechanisms to maintain
them.

\section*{Discussion}       	

	We have shown that a plant that shares an infectious disease with
a competitor may evolve less than total resistance to the pathogen.  This can
occur not because of direct physiological costs of resistance, but because
the disease may be a weapon that one host uses against the other.  Under
fairly generic conditions, lowering resistance by an inferior competitor
led to a higher population in a nonspatial model.  However, an ESS
analysis showed that in this case evolution always leads to higher 
resistance, since invading phenotypes can cheat and take advantage of 
the resident phenotype's infection.  Evolution only led to lower resistance
if the disease transmission from user to attacked host was localized in
space.   This allows less resistant phenotypes to benefit from the disease
by primarily infecting their own competitors.  

	The ESS level of resistance depended on the severity of the damage
that the disease inflicted on each host.  If infection of the user host
was sufficiently benign, selection favored nonresistance.  As the level
of damage increased, higher resistance was favored; if the damage was too
severe, selection favored total resistance.  On the other hand, the greater
the damage the disease inflicted on the attacked host, the more potent a 
weapon it was, and selection favored lower resistance by the user.  The
level of resistance by the attacked host also impacted the evolution of 
resistance by the user.  As the attacked host's resistance increased,
the user was selected to offset this by lowering its own resistance.
This has interesting implications for the coevolution of resistance and
virulence in the system, with the possibility of a multilateral arms
race.  

	The evolution of resistance was quite robust with respect
to the pathogen life cycle and the dispersal scale of the attacked host.
When we compared systems with approximately equal resident densities,
we found that lowered resistance is favored over a similar range of
disease damage.  Thus, the tradeoff between the direct damage of the disease
and the indirect benefit seems to depend in our model more on the overall
plant and pathogen
densities than on the details of the pathogen life cycle or
on the dispersal scales of the attacked host and attacked--user transmission.

	Our conclusion that evolution may lead to reduced resistance
by a plant that shares a pathogen with a superior competitor supports the
feasibility of the mechanism that Rice and Westoby (1982) invoked to explain
ecological and phylogenetic patterns of heteroecious fungi.  
We have generalized their argument
by showing that selection may favor lowered resistance whenever a user
host gains a net advantage from the disease, provided that it transmits
the pathogen to its own competitors.  It may be possible to extend the
results still further, to any system in which species compete both directly
and through a common enemy.  When both direct and indirect competition
occur, a species will undergo selection on traits which determine its
vulnerability to the predator or pathogen.  A tradeoff can then occur between
the direct damage of infection and the indirect benefits of apparent 
competition.
In our system, spatial structure was essential in determing the evolutionary
consequences of this tradeoff.  There may be other types of localization
that allow less ``resistant'' phenotypes to benefit from their strategy.
For example, social interaction networks can determine the spread of diseases
in animals (Keeling et al., 1997; Keeling and Grenfell, 2000).  
If an individual can preferentially target for infection
its own competitors, selection may favor lowered resistance.
In general, the coevolution of hosts and pathogens or
predators and prey may depend not only on the direct effects of the
pathogen or predator, but on its impact on other species in the system. 
The resistance strategy of the host or prey must then be interpreted
in terms of a tradeoff between direct and indirect effects.

\section*{Acknowledgments}

This research was conducted with support from NSF DBI-9602226, the
Research Training Grant -- Nonlinear Dynamics in Biology, awarded to the
University of California, Davis.   

\pagebreak

\section*{Appendix 1}

We present the pair approximation equations for the resident population
densities.  Let $P_{\sigma\sigma'}$ be the ``pair density'', the
probability that a randomly
chosen site is in state $\sigma$ and a randomly chosen neighbor of it
is in state $\sigma'$.  Then:
\begin{equation}
P_{\sigma|\sigma'} = \frac{P_{\sigma\sigma'}}{P_{\sigma'}}
\end{equation}
Note that $P_{\sigma\sigma'} = P_{\sigma'\sigma}$, but that
in general $P_{\sigma|\sigma'} \ne P_{\sigma'|\sigma}$.
The dynamics of the pairs can be obtained by considering interactions
between the two sites and interactions with neighbors of them.  These
last interactions involve triplet terms of the form
$P_{\sigma|\sigma'\sigma''}$.  For example,
consider again the colonization of an empty site by a susceptible
attacked host.  This event switches an $S_{A}E$ pair to $S_{A}S_{A}$.
If dispersal is local, the colonization can occur due to reproduction
by the first $S_{A}$ site or by any of the three other neighbors of
the the empty site.  Thus, local reproduction by $S_{A}$ switches
$S_{A}E$ pairs to $S_{A}S_{A}$ at the rate $\frac{\beta_{A}}{4}
(1 + 3 P_{S_{A}|ES_{A}}) P_{S_{A}E}$.  On the other hand, if dispersal is
global, the offspring could have come from anywhere, and the
rate of change due to this process is simply $B_{A} P_{S_{A}} P_{S_{A}E}$.
After replacing the triplet terms with the ordinary pair approximation
($P_{\sigma|\sigma'\sigma''} = P_{\sigma|\sigma'}$), we obtain the
pair density equations for the residents:
\begin{eqnarray}
\frac{1}{2}\dot{P}_{S_{A}S_{A}} & = & 
	[\frac{\beta_{A}}{4}(1 + 3P_{S_{A}|E} + 3P_{I_{A}|E}) 
	+ B_{A} (P_{S_{A}} + P_{I_{A}})] P_{S_{A}E}
	\nonumber \\
    & & + [\frac{\beta_{A}}{4}(1 + 3P_{S_{A}|S_{U}} + 3P_{I_{A}|S_{U}}) 
	+ B_{A} (P_{S_{A}} + P_{I_{A}})] P_{S_{A}S_{U}}
	\nonumber \\
    & &	+ [\frac{\beta_{A}}{4}(1 + 3P_{S_{A}|I_{U}} + 3P_{I_{A}|I_{U}})
	+ B_{A} (P_{S_{A}} + P_{I_{A}})] P_{S_{A}I_{U}}
	\nonumber \\
    & &	- \{\mu_{A} + (1 - \rho_{A})[\frac{3}{4}\gamma_{AA}P_{I_{A}|S_{A}}
	+ \Gamma_{AA}P_{I_{A}} 
	\nonumber \\
  & &	+ \frac{3}{4}\gamma_{UA}P_{I_{U}|S_{A}}
	+ \Gamma_{UA}P_{I_{U}}]\} P_{S_{A}S_{A}} \\
\dot{P}_{S_{A}I_{A}} & = & 
	(1 - \rho_{A})[\frac{3}{4}\gamma_{AA}P_{I_{A}|S_{A}}
        + \Gamma_{AA}P_{I_{A}} + \frac{3}{4}\gamma_{UA}P_{I_{U}|S_{A}}
        + \Gamma_{UA}P_{I_{U}}] P_{S_{A}S_{A}}
	\nonumber \\
  & &	+ [\frac{\beta_{A}}{4}(1 + 3P_{S_{A}|E} + 3P_{I_{A}|E}) 
	+ B_{A} (P_{S_{A}} + P_{I_{A}})] P_{I_{A}E} 
	\nonumber \\	  
  & &	+ [\frac{\beta_{A}}{4}(1 + 3P_{S_{A}|S_{U}} + 3P_{I_{A}|S_{U}})
        + B_{A} (P_{S_{A}} + P_{I_{A}})] P_{I_{A}S_{U}}
	\nonumber \\
  & &	+ [\frac{\beta_{A}}{4}(1 + 3P_{S_{A}|I_{U}} + 3P_{I_{A}|I_{U}})
        + B_{A} (P_{S_{A}} + P_{I_{A}})] P_{I_{A}I_{U}}
  	- \{\mu_{A} + \alpha_{A}\mu_{A} 
	\nonumber \\
  & &	+ (1 - \rho_{A})
	[\frac{\gamma_{AA}}{4}(1 + 3P_{I_{A}|S_{A}} + \Gamma_{AA}P_{I_{A}}
	\nonumber \\
  & &	+ \frac{3}{4}\gamma_{UA}P_{I_{U}|S_{A}} + \Gamma_{UA}P_{I_{U}}]\} 
	P_{S_{A}I_{A}} \\
\dot{P}_{S_{A}S_{U}} & = &
	[\frac{3}{4}\beta_{A}(P_{S_{A}|E} + P_{I_{A}|E}) 
	+ B_{A} (P_{S_{A}} + P_{I_{A}})] P_{S_{U}E}
	\nonumber \\
  & &	+ [\frac{3}{4}\beta_{A}(P_{S_{A}|S_{U}} + P_{I_{A}|S_{U}}) 
        + B_{A} (P_{S_{A}} + P_{I_{A}})] P_{S_{U}S_{U}}
	\nonumber \\
  & & 	+ [\frac{3}{4}\beta_{A}(P_{S_{A}|I_{U}} + P_{I_{A}|I_{U}})
        + B_{A} (P_{S_{A}} + P_{I_{A}})] P_{S_{U}I_{U}}
        \nonumber \\
  & &	+ [\frac{3}{4}\beta_{U}(P_{S_{U}|E} + P_{I_{U}|E}) 
        + B_{U} (P_{S_{U}} + P_{I_{U}})] P_{S_{A}E}
        \nonumber \\
  & &	- \{\mu_{A} + \mu_{U} + \frac{\beta_{A}}{4}(1 + 3P_{S_{A}|S_{U}}
	+ 3P_{I_{A}|S_{U}} + B_{A} (P_{S_{A}} + P_{I_{A}})
	\nonumber \\
  & &	+ (1 - \rho_{U})[\frac{3}{4}\gamma_{UU}P_{I_{U}|S_{U}} + 
        \Gamma_{UU}P_{I_{U}} + \frac{3}{4}\gamma_{AU}P_{I_{A}|S_{U}} + 
	\Gamma_{AU}P_{I_{A}}]
	\nonumber \\
  & &	+ (1 - \rho_{A})[\frac{3}{4}\gamma_{AA}P_{I_{A}|S_{A}} + 
        \Gamma_{AA}P_{I_{A}} 
	\nonumber \\
  & &	+ \frac{3}{4}\gamma_{UA}P_{I_{U}|S_{A}} +
        \Gamma_{UA}P_{I_{U}}]\} P_{S_{A}S_{U}} \\
\dot{P}_{S_{A}I_{U}} & = &
  	[\frac{3}{4}\beta_{A}(P_{S_{A}|E} + P_{I_{A}|E})
        + B_{A} (P_{S_{A}} + P_{I_{A}})] P_{I_{U}E}
        \nonumber \\
  & &	+ [\frac{3}{4}\beta_{A}(P_{S_{A}|S_{U}} + P_{I_{A}|S_{U}})
        + B_{A} (P_{S_{A}} + P_{I_{A}})] P_{S_{U}I_{U}}
        \nonumber \\
  & &   + [\frac{3}{4}\beta_{A}(P_{S_{A}|I_{U}} + P_{I_{A}|I_{U}})
        + B_{A} (P_{S_{A}} + P_{I_{A}})] P_{I_{U}I_{U}}
        \nonumber \\
  & &	+ (1 - \rho_{U})[\frac{3}{4}\gamma_{UU}P_{I_{U}|S_{U}} +
        \Gamma_{UU}P_{I_{U}} + \frac{3}{4}\gamma_{AU}P_{I_{A}|S_{U}} +
        \Gamma_{AU}P_{I_{A}}] P_{S_{A}S_{U}}
        \nonumber \\
  & &	- \{\mu_{A} + \alpha_{U}\mu_{U} 
	+ [\frac{\beta_{A}}{4}(1 + 3P_{S_{A}|I_{U}} + 3P_{I_{A}|I_{U}})
        + B_{A} (P_{S_{A}} + P_{I_{A}})]
        \nonumber \\
  & &	+ (1 - \rho_{A})[\frac{3}{4}\gamma_{AA}P_{I_{A}|S_{A}} +
        \Gamma_{AA}P_{I_{A}} 
	\nonumber \\
  & &	+ \frac{\gamma_{UA}}{4}(1 + 3P_{I_{U}|S_{A}})
	+ \Gamma_{UA} P_{I_{U}}]\} P_{S_{A}I_{U}} \\
\frac{1}{2}\dot{P}_{I_{A}I_{A}} & = & 
	(1 - \rho_{A})[\frac{\gamma_{AA}}{4}(1 + 3P_{I_{A}|S_{A}})
	+ \Gamma_{AA}P_{I_{A}} + \frac{3}{4}\gamma_{UA}P_{I_{U}|S_{A}}
	+ \Gamma_{UA}P_{I_{U}}] P_{S_{A}I_{A}}
	\nonumber \\
  & &	- \alpha_{A}\mu_{A}P_{I_{A}I_{A}} \\
\dot{P}_{I_{A}S_{U}} & = & 
	(1 - \rho_{A})[\frac{3}{4}\gamma_{AA}P_{I_{A}|S_{A}} +
	\Gamma_{AA}P_{I_{A}} + \frac{3}{4}\gamma_{UA}P_{I_{U}|S_{A}}
        + \Gamma_{UA}P_{I_{U}}] P_{S_{A}S_{U}}
	\nonumber \\
  & & 	+ [\frac{3}{4}\beta_{U}(P_{S_{U}|E} + P_{I_{U}|E})
        + B_{U} (P_{S_{U}} + P_{I_{U}})] P_{I_{A}E}
        \nonumber \\
  & &	- \{\alpha_{A}\mu_{A} + \mu_{U} + \frac{\beta_{A}}{4}(1 + 
	3P_{I_{A}|S_{U}} + 3P_{S_{A}|S_{U}})
	+ B_{A} (P_{S_{A}} + P_{I_{A}})
	\nonumber \\
  & &	+ (1 - \rho_{U})[\frac{\gamma_{AU}}{4}(1 + 3P_{I_{A}|S_{U}})
	+ \Gamma_{AU}P_{I_{A}} 
	\nonumber \\
  & &	+ \frac{3}{4}\gamma_{UU}P_{I_{U}|S_{U}}
	+ \Gamma_{UU}P_{I_{U}}]\} P_{I_{A}S_{U}} \\
\dot{P}_{I_{A}I_{U}} & = &
	(1 - \rho_{A})[\frac{3}{4}\gamma_{AA}P_{I_{A}|S_{A}} +
        \Gamma_{AA}P_{I_{A}} + \frac{\gamma_{UA}}{4}(1 + 3P_{I_{U}|S_{A}})
        + \Gamma_{UA} P_{I_{U}}] P_{S_{A}I_{U}} 
	\nonumber \\
  & &	+ (1 - \rho_{U})[\frac{\gamma_{AU}}{4}(1 + 3P_{I_{A}|S_{U}})
        + \Gamma_{AU}P_{I_{A}} + \frac{3}{4}\gamma_{UU}P_{I_{U}|S_{U}}
        + \Gamma_{UU}P_{I_{U}}] P_{I_{A}S_{U}} 
	\nonumber \\
  & &	- \{\alpha_{A}\mu_{A} + \alpha_{U}\mu_{U} 
	+ \frac{\beta_{A}}{4}(1 + 3P_{S_{A}|I_{U}} + 3P_{I_{A}|I_{U}})
	\nonumber \\
  & &   + B_{A} (P_{S_{A}} + P_{I_{A}})\} P_{I_{A}I_{U}} \\
\frac{1}{2}\dot{P}_{S_{U}S_{U}} & = &
	[\frac{\beta_{U}}{4}(1 + 3P_{S_{U}|E} + 3P_{I_{U}|E})
	+ B_{U} (P_{S_{U}} + P_{I_{U}})] P_{S_{U}E}
	\nonumber \\
  & &	- \{\mu_{U} + (1 - \rho_{U})[\frac{3}{4}\gamma_{UU}P_{I_{U}|S_{U}}
	+ \Gamma_{UU}P_{I_{U}} + \frac{3}{4}\gamma_{AU}P_{I_{A}|S_{U}}
	+ \Gamma_{AU}P_{I_{A}}] 
	\nonumber \\
  & &	+ \frac{3}{4}\beta_{A}(P_{S_{A}|S_{U}} + P_{I_{A}|S_{U}})
	+ B_{A} (P_{S_{A}} + P_{I_{A}})\} P_{S_{U}S_{U}} \\
\dot{P}_{S_{U}I_{U}} & = &
	[\frac{\beta_{U}}{4}(1 + 3P_{I_{U}|E} + 3P_{S_{U}|E})
	+ B_{U} (P_{S_{U}} + P_{I_{U}})] P_{I_{U}E}
	\nonumber \\
  & &	+ (1 - \rho_{U})[\frac{3}{4}\gamma_{UU}P_{I_{U}|S_{U}}
        + \Gamma_{UU}P_{I_{U}} + \frac{3}{4}\gamma_{AU}P_{I_{A}|S_{U}}
        + \Gamma_{AU}P_{I_{A}}] P_{S_{U}S_{U}}
        \nonumber \\
  & &	- \{\mu_{U} + \alpha_{U}\mu_{U}
	+ \frac{3}{4}\beta_{A}(P_{S_{A}|S_{U}} + P_{I_{A}|S_{U}}
	+ P_{S_{A}|I_{U}} + P_{I_{A}|I_{U}}) 
	\nonumber \\
  & &	+ 2B_{A} (P_{S_{A}} + P_{I_{A}})
  	+ (1 - \rho_{U})[\frac{\gamma_{UU}}{4}(1 + 3P_{I_{U}|S_{U}})
	+ \Gamma_{UU}P_{I_{U}} 
	\nonumber \\
  & &	+ \frac{3}{4}\gamma_{AU}P_{I_{A}|S_{U}}
        + \Gamma_{AU}P_{I_{A}}]\} P_{S_{U}I_{U}} \\
\frac{1}{2}\dot{P}_{I_{U}I_{U}} & = &
	(1 - \rho_{U})[\frac{\gamma_{UU}}{4}(1 + 3P_{I_{U}|S_{U}})
        + \Gamma_{UU}P_{I_{U}} + \frac{3}{4}\gamma_{AU}P_{I_{A}|S_{U}}
        + \Gamma_{AU}P_{I_{A}}] P_{S_{U}I_{U}}
	\nonumber \\
  & &	- \{\alpha_{U}\mu_{U} 
	+ \frac{3}{4}\beta_{A}(P_{S_{A}|I_{U}} + P_{I_{A}|I_{U}})
	+ B_{A} (P_{S_{A}} + P_{I_{A}})\} P_{I_{U}I_{U}}.
\end{eqnarray}
Since $P_{\sigma} = \sum_{\sigma'} P_{\sigma\sigma'}$, we can augment these
equations with equations 1--4 for the singleton densities to obtain a
closed system.  The system is too complex to solve analytically, but
numerical solution is straightforward and allows efficient exploration
of parameter space.  A more detailed explanation of the derivation
of pair density equations can be found in, for example, Matsuda (1992)
or Rand (1999).

\section*{Appendix 2}

We present the conditional pair equations for the neighborhood structure
of the invading phenotype at low density.  First we find equations for
pair densities involving $S$ or $I$, by the same procedure used above.
Now, the low density assumption means that $P_{\sigma} = 
P_{\sigma\sigma'} = 0$, when $\sigma \in {S,I}$.  However, the conditional
probabilities $P_{\sigma'|\sigma}$ are not necessarily small in this
case.  To obtain these quantities, we use the fact that
\begin{equation}
\dot{P}_{\sigma'|\sigma} = \frac{1}{P_{\sigma}}\dot{P}_{\sigma'\sigma}
	- P_{\sigma'|\sigma}\frac{\dot{P}_{\sigma}}{P_{\sigma}}
\end{equation}
to derive the dynamics of the conditional probabilities from those
of the singleton and pair densities.  Using the pair approximation
and the low density assumption, we obtain: 
\begin{eqnarray}
\dot{P}_{S_{A}|S} & = &
        [\frac{3}{4}\beta_{A}(P_{S_{A}|E} + P_{I_{A}|E})
        + B_{A} (P_{S_{A}} + P_{I_{A}})] P_{E|S}
        \nonumber \\
  & &   + [\frac{3}{4}\beta_{A}(P_{S_{A}|S} + P_{I_{A}|S})
        + B_{A} (P_{S_{A}} + P_{I_{A}})] P_{S|S}
        \nonumber \\
  & &   + [\frac{3}{4}\beta_{A}(P_{S_{A}|S_{U}} + P_{I_{A}|S_{U}})
        + B_{A} (P_{S_{A}} + P_{I_{A}})] P_{S_{U}|S}
        \nonumber \\
  & &   + [\frac{3}{4}\beta_{A}(P_{S_{A}|I} + P_{I_{A}|I})
        + B_{A} (P_{S_{A}} + P_{I_{A}})] P_{I|S}
        \nonumber \\
  & &   + [\frac{3}{4}\beta_{A}(P_{S_{A}|I_{U}} + P_{I_{A}|I_{U}})
        + B_{A} (P_{S_{A}} + P_{I_{A}})] P_{I_{U}|S}
        \nonumber \\
  & &   + [\frac{3}{4}\beta_{U}(P_{E|S} + P_{E|I}\frac{P_{I}}{P_{S}})
        + B_{U}P_{E}(1 + \frac{P_{I}}{P_{S}})] P_{S_{A}|E}
        \nonumber \\
  & &   - \{\mu_{A} + \frac{\beta_{A}}{4}(1 - P_{S_{A}|S} - P_{I_{A}|S})
	\nonumber \\
  & &	+ (1 - \rho')[-\frac{\gamma_{AU}}{4}P_{I_{A}|S}
	- \frac{\gamma_{UU}}{4}(P_{I_{U}|S} + P_{I|S})]
	\nonumber \\
  & &	+ (1 - \rho_{A})[\frac{3}{4}\gamma_{AA}P_{I_{A}|S_{A}} +
        \Gamma_{AA}P_{I_{A}} + \frac{3}{4}\gamma_{UA}P_{I_{U}|S_{A}} +
        \Gamma_{UA}P_{I_{U}}]
	\nonumber \\
  & &	+ \beta_{U}(P_{E|S} + P_{E|I}\frac{P_{I}}{P_{S}})
	+ B_{U}P_{E}(1 + \frac{P_{I}}{P_{S}})\} P_{S_{A}|S} \\
\dot{P}_{I_{A}|S} & = &
        (1 - \rho_{A})[\frac{3}{4}\gamma_{AA}P_{I_{A}|S_{A}} +
        \Gamma_{AA}P_{I_{A}} + \frac{3}{4}\gamma_{UA}P_{I_{U}|S_{A}}
        + \Gamma_{UA}P_{I_{U}}] P_{S_{A}|S}
        \nonumber \\
  & &   + [\frac{3}{4}\beta_{U}(P_{E|S} + P_{E|I}\frac{P_{I}}{P_{S}})
        + B_{U}P_{E}(1 + \frac{P_{I}}{P_{S}})] P_{I_{A}|E}
        \nonumber \\
  & &   - \{\alpha_{A}\mu_{A} + (1 - \rho')[\frac{\gamma_{AU}}{4}
	(1 - P_{I_{A}|S}) - \frac{\gamma_{UU}}{4}(P_{I_{U}|S} + P_{I|S})]
        \nonumber \\
  & &	+ \frac{\beta_{A}}{4}(1 - P_{S_{A}|S} - P_{I_{A}|S})
	+ \beta_{U}(P_{E|S} + P_{E|I}\frac{P_{I}}{P_{S}})
        \nonumber \\
  & &	+ B_{U}P_{E}(1 + \frac{P_{I}}{P_{S}})\} P_{I_{A}|S} \\
\dot{P}_{S_{U}|S} & = &
        [\frac{\beta_{U}}{4}(1 + 3P_{S_{U}|E} + 3P_{I_{U}|E})
        + B_{U} (P_{S_{U}} + P_{I_{U}})] P_{E|S}
        \nonumber \\
  & &   + [\frac{3}{4}\beta_{U}(P_{E|S} + P_{E|I}\frac{P_{I}}{P_{S}})
        + B_{U}P_{E}(1 + \frac{P_{I}}{P_{S}})] P_{S_{U}|E}
        \nonumber \\
  & &   - \{\mu_{U} + (1 - \rho_{U})[\frac{3}{4}\gamma_{UU}P_{I_{U}|S_{U}}
        + \Gamma_{UU}P_{I_{U}} + \frac{3}{4}\gamma_{AU}P_{I_{A}|S_{U}}
        + \Gamma_{AU}P_{I_{A}}]
        \nonumber \\
  & &	+ (1 - \rho')[-\frac{\gamma_{AU}}{4}P_{I_{A}|S}
        - \frac{\gamma_{UU}}{4}(P_{I_{U}|S} + P_{I|S})]
        \nonumber \\
  & &	+ \frac{3}{4}\beta_{A}(P_{S_{A}|S_{U}} + P_{I_{A}|S_{U}})
        + B_{A} (P_{S_{A}} + P_{I_{A}}) - \frac{\beta_{A}}{4}
	(P_{S_{A}|S} + P_{I_{A}|S})
	\nonumber \\
  & &	+ \beta_{U}(P_{E|S} + P_{E|I}\frac{P_{I}}{P_{S}})
        + B_{U}P_{E}(1 + \frac{P_{I}}{P_{S}})\} P_{S_{U}|S} \\
\dot{P}_{I_{U}|S} & = &
	[\frac{3}{4}\beta_{U}(P_{E|S} + P_{E|I}\frac{P_{I}}{P_{S}})
	+ B_{U}P_{E}(1 + \frac{P_{I}}{P_{S}})] P_{I_{U}|E}
	\nonumber \\
  & &	+ (1 - \rho_{U})[\frac{3}{4}\gamma_{UU}P_{I_{U}|S_{U}}
        + \Gamma_{UU}P_{I_{U}} + \frac{3}{4}\gamma_{AU}P_{I_{A}|S_{U}}
        + \Gamma_{AU}P_{I_{A}}] P_{S_{U}|S}
	\nonumber \\
  & &	- \{\alpha_{U}\mu_{U}
	+ \frac{3}{4}\beta_{A}(P_{S_{A}|I_{U}} + P_{I_{A}|I_{U}})
	+ B_{A} (P_{S_{A}} + P_{I_{A}})
	- \frac{\beta_{A}}{4}(P_{S_{A}|S} + P_{I_{A}|S})
	\nonumber \\
  & &	+ (1 - \rho')[-\frac{\gamma_{AU}}{4}P_{I_{A}|S}
        + \frac{\gamma_{UU}}{4}(1 - P_{I_{U}|S} - P_{I|S})]
        \nonumber \\
  & &   + \beta_{U}(P_{E|S} + P_{E|I}\frac{P_{I}}{P_{S}})
        + B_{U}P_{E}(1 + \frac{P_{I}}{P_{S}})\} P_{I_{U}|S} \\
\dot{P}_{S|S} & = &
	\frac{\beta_{U}}{2}P_{E|S} - \{\mu_{U}
	+ (1 - \rho')[\frac{\gamma_{UU}}{2}(P_{I_{U}|S} + P_{I|S})
	+ \Gamma_{UU} P_{I_{U}} 
	\nonumber \\
  & &	+ \frac{\gamma_{AU}}{2}P_{I_{A}|S} + \Gamma_{AU}P_{I_{A}}]
  	+ \frac{\beta_{A}}{2}(P_{S_{A}|S} + P_{I_{A}|S})
	+ B_{A}(P_{S_{A}} + P_{I_{A}})
	\nonumber \\
  & &	+ \beta_{U}(P_{E|S} + P_{E|I}\frac{P_{I}}{P_{S}})
  	+ B_{U}P_{E}(1 + \frac{P_{I}}{P_{S}})\} P_{S|S} \\
\dot{P}_{I|S} & = &
	\frac{\beta_{U}}{4}P_{E|I}\frac{P_{I}}{P_{S}}
	+ (1 - \rho')[\frac{3}{4}\gamma_{UU}P_{I_{U}|S}
        + \Gamma_{UU} P_{I_{U}} + \frac{3}{4}\gamma_{AU}P_{I_{A}|S}
        + \Gamma_{AU}P_{I_{A}}] P_{S|S}
	\nonumber \\
  & &	- \{\alpha_{U}\mu_{U}
	+ \frac{3}{4}\beta_{A}(P_{S_{A}|I} + P_{I_{A}|I})
	+ B_{A} (P_{S_{A}} + P_{I_{A}})
	- \frac{\beta_{A}}{4}(P_{S_{A}|S} + P_{I_{A}|S})
	\nonumber \\
  & &	+ (1 - \rho')[-\frac{\gamma_{AU}}{4}P_{I_{A}|S}
        + \frac{\gamma_{UU}}{4}(1 - P_{I_{U}|S} - P_{I|S})]
        \nonumber \\
  & &   + \beta_{U}(P_{E|S} + P_{E|I}\frac{P_{I}}{P_{S}})
        + B_{U}P_{E}(1 + \frac{P_{I}}{P_{S}})\} P_{I|S} \\
\dot{P}_{S_{A}|I} & = &
        [\frac{3}{4}\beta_{A}(P_{S_{A}|E} + P_{I_{A}|E})
        + B_{A} (P_{S_{A}} + P_{I_{A}})] P_{E|I}
        \nonumber \\
  & &   + [\frac{3}{4}\beta_{A}(P_{S_{A}|I_{U}} + P_{I_{A}|I_{U}})
        + B_{A} (P_{S_{A}} + P_{I_{A}})] P_{I_{U}|I}
        \nonumber \\
  & &   + [\frac{3}{4}\beta_{A}(P_{S_{A}|I} + P_{I_{A}|I})
        + B_{A} (P_{S_{A}} + P_{I_{A}})] P_{I|I}
        \nonumber \\
  & &   + [\frac{3}{4}\beta_{A}(P_{S_{A}|S_{U}} + P_{I_{A}|S_{U}})
        + B_{A} (P_{S_{A}} + P_{I_{A}})] P_{S_{U}|I}
        \nonumber \\
  & &   + [\frac{3}{4}\beta_{A}(P_{S_{A}|S} + P_{I_{A}|S})
        + B_{A} (P_{S_{A}} + P_{I_{A}})] P_{S|I}
        \nonumber \\
  & &   + (1 - \rho')[\frac{3}{4}\gamma_{UU}(P_{I_{U}|S} + P_{I|S}) +
        \Gamma_{UU}P_{I_{U}} + \frac{3}{4}\gamma_{AU}P_{I_{A}|S} +
        \Gamma_{AU}P_{I_{A}}] P_{S_{A}|S}\frac{P_{S}}{P_{I}}
        \nonumber \\
  & &	- \{\mu_{A} + \frac{\beta_{A}}{4}(1 - P_{S_{A}|I} - P_{I_{A}|I})
	\nonumber \\
  & &	+ (1 - \rho_{A})[\frac{3}{4}\gamma_{AA}P_{I_{A}|S_{A}} +
        \Gamma_{AA}P_{I_{A}} + \frac{\gamma_{UA}}{4}(1 + 3P_{I_{U}|S_{A}})
        + \Gamma_{UA} P_{I_{U}}]
	\nonumber \\
  & &	+ (1 - \rho')[\gamma_{UU}(P_{I_{U}|S} + P_{I|S}) +
        \Gamma_{UU}P_{I_{U}} 
	\nonumber \\
  & &	+ \gamma_{AU}P_{I_{A}|S} + \Gamma_{AU}P_{I_{A}}]
	\frac{P_{S}}{P_{I}}\} P_{S_{A}|I} \\
\dot{P}_{I_{A}|I} & = &
	(1 - \rho_{A})[\frac{3}{4}\gamma_{AA}P_{I_{A}|S_{A}} +
        \Gamma_{AA}P_{I_{A}} + \frac{\gamma_{UA}}{4}(1 + 3P_{I_{U}|S_{A}})
        + \Gamma_{UA} P_{I_{U}}] P_{S_{A}|I}
        \nonumber \\
  & &	+ (1 - \rho')[\frac{3}{4}\gamma_{UU}(P_{I_{U}|S} + P_{I|S}) +
        \Gamma_{UU}P_{I_{U}} 
	\nonumber \\
  & &	+ \frac{\gamma_{AU}}{4}(1 + 3P_{I_{A}|S}) +
        \Gamma_{AU}P_{I_{A}}] P_{I_{A}|S}\frac{P_{S}}{P_{I}}
  	- \{\alpha_{A}\mu_{A} 
	+ \frac{\beta_{A}}{4}(1 - P_{S_{A}|I} - P_{I_{A}|I})
	\nonumber \\
  & &   + (1 - \rho')[\gamma_{UU}(P_{I_{U}|S} + P_{I|S}) +
        \Gamma_{UU}P_{I_{U}} 
	\nonumber \\
  & &	+ \gamma_{AU}P_{I_{A}|S} + \Gamma_{AU}P_{I_{A}}]
        \frac{P_{S}}{P_{I}}\} P_{I_{A}|I} \\
\dot{P}_{S_{U}|I} & = &
        [\frac{3}{4}\beta_{U}(P_{I_{U}|E} + P_{S_{U}|E})
        + B_{U} (P_{S_{U}} + P_{I_{U}})] P_{E|I}
        \nonumber \\
  & &   + (1 - \rho')[\frac{3}{4}\gamma_{UU}(P_{I_{U}|S} + P_{I|S})
        + \Gamma_{UU}P_{I_{U}} + \frac{3}{4}\gamma_{AU}P_{I_{A}|S}
        + \Gamma_{AU}P_{I_{A}}] P_{S_{U}|S}\frac{P_{S}}{P_{I}}
        \nonumber \\
  & &   - \{\mu_{U} - \frac{\beta_{A}}{4}(P_{S_{A}|I} + P_{I_{A}|I})
	+ \frac{3}{4}\beta_{A}(P_{S_{A}|S_{U}} + P_{I_{A}|S_{U}})
	+ B_{A} (P_{S_{A}} + P_{I_{A}})
	\nonumber \\
  & &	+ (1 - \rho_{U})[\frac{\gamma_{UU}}{4}(1 + 3P_{I_{U}|S_{U}})
        + \Gamma_{UU}P_{I_{U}} + \frac{3}{4}\gamma_{AU}P_{I_{A}|S_{U}}
        + \Gamma_{AU}P_{I_{A}}]
	\nonumber \\
  & &	+ (1 - \rho')[\gamma_{UU}(P_{I_{U}|S} + P_{I|S}) +
        \Gamma_{UU}P_{I_{U}} 
	\nonumber \\
  & &	+ \gamma_{AU}P_{I_{A}|S} + \Gamma_{AU}P_{I_{A}}]
        \frac{P_{S}}{P_{I}}\} P_{S_{U}|I} \\
\dot{P}_{I_{U}|I} & = &
	(1 - \rho_{U})[\frac{\gamma_{UU}}{4}(1 + 3P_{I_{U}|S_{U}})
        + \Gamma_{UU}P_{I_{U}} + \frac{3}{4}\gamma_{AU}P_{I_{A}|S_{U}}
        + \Gamma_{AU}P_{I_{A}}] P_{S_{U}|I}
	\nonumber \\
  & &	+ (1 - \rho')[\frac{\gamma_{UU}}{4}(1 + 3P_{I_{U}|S} + 3P_{I|S})
        + \Gamma_{UU}P_{I_{U}} 
	\nonumber \\
  & &	+ \frac{3}{4}\gamma_{AU}P_{I_{A}|S}
        + \Gamma_{AU}P_{I_{A}}] P_{I_{U}|S}\frac{P_{S}}{P_{I}}
	\nonumber \\
  & &	- \{\alpha_{U}\mu_{U} - \frac{\beta_{A}}{4}(P_{S_{A}|I} + P_{I_{A}|I})
        + \frac{3}{4}\beta_{A}(P_{S_{A}|I_{U}} + P_{I_{A}|I_{U}})
        + B_{A} (P_{S_{A}} + P_{I_{A}})
        \nonumber \\
  & &   + (1 - \rho')[\gamma_{UU}(P_{I_{U}|S} + P_{I|S}) +
        \Gamma_{UU}P_{I_{U}} 
	\nonumber \\
  & &	+ \gamma_{AU}P_{I_{A}|S} + \Gamma_{AU}P_{I_{A}}]
        \frac{P_{S}}{P_{I}}\} P_{S_{U}|I} \\
\dot{P}_{I|I} & = &
	(1 - \rho')[\frac{\gamma_{UU}}{2}(1 + 3P_{I_{U}|S} + 3P_{I|S})
        + 2\Gamma_{UU}P_{I_{U}} + \frac{3}{2}\gamma_{AU}P_{I_{A}|S}
        \nonumber \\
  & &	+ 2\Gamma_{AU}P_{I_{A}}] P_{S|I}
  	- \{\alpha_{U}\mu_{U} + \frac{\beta_{A}}{2}(P_{S_{A}|I} + P_{I_{A}|I})
	+ B_{A} (P_{S_{A}} + P_{I_{A}})
	\nonumber \\
  & &   + (1 - \rho')[\gamma_{UU}(P_{I_{U}|S} + P_{I|S}) +
        \Gamma_{UU}P_{I_{U}} 
	\nonumber \\
  & &	+ \gamma_{AU}P_{I_{A}|S} + \Gamma_{AU}P_{I_{A}}]
        \frac{P_{S}}{P_{I}}\} P_{I|I} \\
\dot{P}_{S|I} & = &
	\frac{\beta_{U}}{4}P_{E|I}
	+ (1 - \rho')[\frac{3}{4}\gamma_{UU}(P_{I_{U}|S} + P_{I|S})
        + \Gamma_{UU}P_{I_{U}} 
	\nonumber \\
  & &	+ \frac{3}{4}\gamma_{AU}P_{I_{A}|S}
        + \Gamma_{AU}P_{I_{A}}] P_{S|S}\frac{P_{S}}{P_{I}}
        \nonumber \\
  & &	- \{\mu_{U} - \frac{\beta_{A}}{4}(P_{S_{A}|I} + P_{I_{A}|I})
        + \frac{3}{4}\beta_{A}(P_{S_{A}|S} + P_{I_{A}|S})
        + B_{A} (P_{S_{A}} + P_{I_{A}})
        \nonumber \\
  & &	+ (1 - \rho')[\frac{\gamma_{UU}}{4}(1 + 3P_{I_{U}|S} + 3P_{I|S})
	+ \Gamma_{UU}P_{I_{U}} + \frac{3}{4}\gamma_{AU}P_{I_{A}|S}
	+ \Gamma_{AU}P_{I_{A}}]
	\nonumber \\
  & &   + (1 - \rho')[\gamma_{UU}(P_{I_{U}|S} + P_{I|S}) +
        \Gamma_{UU}P_{I_{U}} 
	\nonumber \\
  & &	+ \gamma_{AU}P_{I_{A}|S} + \Gamma_{AU}P_{I_{A}}]
        \frac{P_{S}}{P_{I}}\} P_{S|I}
\end{eqnarray}
We obtain a closed system by including the resident equilibrium values 
and using:
\begin{equation}
\frac{P_{I}}{P_{S}} = \frac{P_{I|S}}{P_{S|I}}.
\end{equation}
By finding the equilibrium solution
numerically, we obtain the local spatial structure of the invasion
early in its development.

\section*{Appendix 3}

	We present a proof that selection cannot lead to lower user
host resistance in the mean field case.  Let $J$ be the Jacobian matrix
for the new phenotype invading at low density, i.e. equations 5--6
linearized around $P_{S} = P_{I} = 0$ and the resident populations at
the equilibrium solution of equations 1--4.  Then:
\begin{eqnarray}
J & = & \left[ \begin{array}{c}
	B_{U}P^{*}_{E} - B_{A}(P^{*}_{S_{A}}+P^{*}_{I_{A}}) - \mu_{U}
 	- (1-\rho')(\Gamma_{AU}P^{*}_{I_{A}}+\Gamma_{UU}P^{*}_{I_{U}}) \\
	(1-\rho')(\Gamma_{AU}P^{*}_{I_{A}}+\Gamma_{UU}P^{*}_{I_{U}})
	\end{array} \right.
	\nonumber \\
  & & \makebox[2.5in]{} \left. \begin{array}{c}
   	 B_{U}P^{*}_{E} \\
   	 - B_{A}(P^{*}_{S_{A}}+P^{*}_{I_{A}}) - \alpha_{U}\mu_{U}
	\end{array} \right] 
\end{eqnarray}
Now, from equations 3--4, we know that the resident equilibrium
populations satisfy:
\begin{equation}
B_{U}P^{*}_{E} - B_{A}(P^{*}_{S_{A}}+P^{*}_{I_{A}}) - \mu_{U} =
	(1-\rho_{U})(\Gamma_{AU}P^{*}_{I_{A}}+\Gamma_{UU}P^{*}_{I_{U}})
	- B_{U}P^{*}_{E}\frac{P^{*}_{I_{U}}}{P^{*}_{S_{U}}},
\label{eqn:sub}
\end{equation}
and
\begin{equation}
B_{A}(P^{*}_{S_{A}}+P^{*}_{I_{A}}) + \alpha_{U}\mu_{U} =
	(1-\rho_{U})(\Gamma_{AU}P^{*}_{I_{A}}+\Gamma_{UU}P^{*}_{I_{U}})
	\frac{P^{*}_{S_{U}}}{P^{*}_{I_{U}}}.
\end{equation}
After substituting these expressions into $J$ and rearranging, we find:
\begin{equation}
det(J) = (\rho_{U} - \rho')(\Gamma_{AU}P^{*}_{I_{A}}+\Gamma_{UU}P^{*}_{I_{U}})
	\left[(1-\rho_{U})(\Gamma_{AU}P^{*}_{I_{A}}+\Gamma_{UU}P^{*}_{I_{U}})
	\frac{P^{*}_{S_{U}}}{P^{*}_{I_{U}}} - B_{U}P^{*}_{E}\right].
\end{equation}
But from equation~\ref{eqn:sub}, we find that:
\begin{equation}
B_{U}P^{*}_{E} = 
	(1-\rho_{U})(\Gamma_{AU}P^{*}_{I_{A}}+\Gamma_{UU}P^{*}_{I_{U}})
	\frac{P^{*}_{S_{U}}}{P^{*}_{I_{U}}}
	- (\alpha_{U}-1)\mu_{U}\frac{P^{*}_{S_{U}}}{P^{*}_{S_{U}}
	+ P^{*}_{I_{U}}}.
\end{equation}
Therefore, 
\begin{equation}
det(J) = (\rho_{U} - \rho')(\Gamma_{AU}P^{*}_{I_{A}}+\Gamma_{UU}P^{*}_{I_{U}})
	\left[(\alpha_{U}-1)\mu_{U}\frac{P^{*}_{S_{U}}}{P^{*}_{S_{U}}
        + P^{*}_{I_{U}}}\right].
\end{equation}
Provided that $\alpha_{U} \ge 1$ (i.e. infection shortens the host's
lifespan), we see that $det(J)$ has the same sign as
$\rho_{U} - \rho'$.  Thus, if the invader has higher resistance than
the resident, one eigenvalue is negative and the other is positive, so
that the invasion succeeds.  If the invader's resistance is lower than
the resident's, both eigenvalues are positive or both are negative.
Now,
\begin{equation}
tr(J) = (\rho' - \rho_{U})(\Gamma_{AU}P^{*}_{I_{A}}+\Gamma_{UU}P^{*}_{I_{U}})
	- B_{U}P^{*}_{E}\frac{P^{*}_{I_{U}}}{P^{*}_{S_{U}}}
	- B_{A}(P^{*}_{S_{A}}+P^{*}_{I_{A}}) - \alpha_{U}\mu_{U},
\end{equation}
so that $tr(J) \le 0$ when $\rho' \le \rho_{U}$.  Thus, when the invader's
resistance is lower than the resident's, both eigenvalues are negative
and the invasion fails. 

\pagebreak

\section*{Figure Captions}

Figure 1: Equilibrium populations of user (a) and attacked (b) hosts as a
        function of user's resistance in the mean field heteroecious model.
        Increasing resistance lowers the user's population level,
        indicating that the disease provides a net benefit.  Note the
        resistance threshold: the disease dies out if the resistance
        exceeds $0.6$.  Parameters are: $B_{A}=0.215, B_{U}=4.5,
        \Gamma_{AU}=2.0, \Gamma_{UA}=3.0, \alpha_{A}=2.3, \alpha_{U}=1.1,
        \rho_{A}=0$.
\newline
\newline
Figure 2: Evolutionary flux in user's resistance. (a) Predictions by pair
        approximation for different levels of damage to user.  Values
        of $\alpha_{U}$ for the curves are, from bottom to top: $1.0,
        1.03, 1.06, 1.09$.  The ESS resistance occurs where the flux
        is zero.  Other parameters are: $B_{A}=.215, \beta_{U}=6.3,
        \Gamma_{AU}=1.7, \gamma_{UA}=6.0, \alpha_{A}=2.4, \rho_{A}=0$.
        (b) Comparison of pair approximation and simulations when
        $\alpha_{U}=1.03$.  Dotted curve is flux predicted by
        pair approximation, as in (a).  Box plot gives results from
        simulations, with standard error bars.  We ran 10 simulations
        on a $1000 \times 1000$ lattice for each resistance level.
        Invader's resistance was $\rho_{U} \pm 0.15$.  After resident
        system reached equilibrium, invader was introduced by switching
        $0.5\%$ of all sites to the susceptible invader state.
        Invasion eigenvalues
        were estimated by fitting an exponential curve to the total
        invader population for 500 time steps after introduction.
\newline
\newline
Figure 3: ESS user resistance as a function of damage done by the pathogen.
        Dashed lines give threshold resistance levels above which the
        disease dies out.  Curves give ESS resistance for different levels
        of damage to the attacked host.  Other parameters are as in
        Figure 2.
\newline
\newline
Figure 4: Relationship between ESS user resistance and the resistance level
        of the attacked host.  Dashed lines give threshold user resistance
        levels; curves give ESS user resistance levels for different
        levels of attacked host resistance.  Other parameters are as in
        Figure 2.
\newline
\newline
Figure 5: ESS user resistance for different dispersal scales of
        the attacked host and attacked--user transmission.

\pagebreak

\section*{Figures}

Figure 1a:

\begin{figure}[h]

        \includegraphics{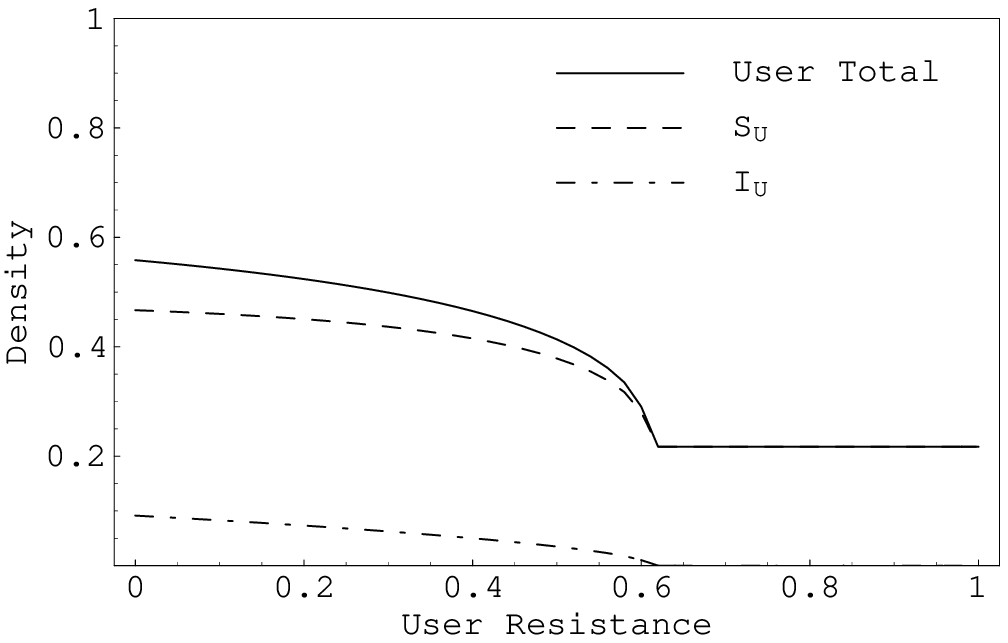}

\end{figure}

Figure 1b:

\begin{figure}[h]

        \includegraphics{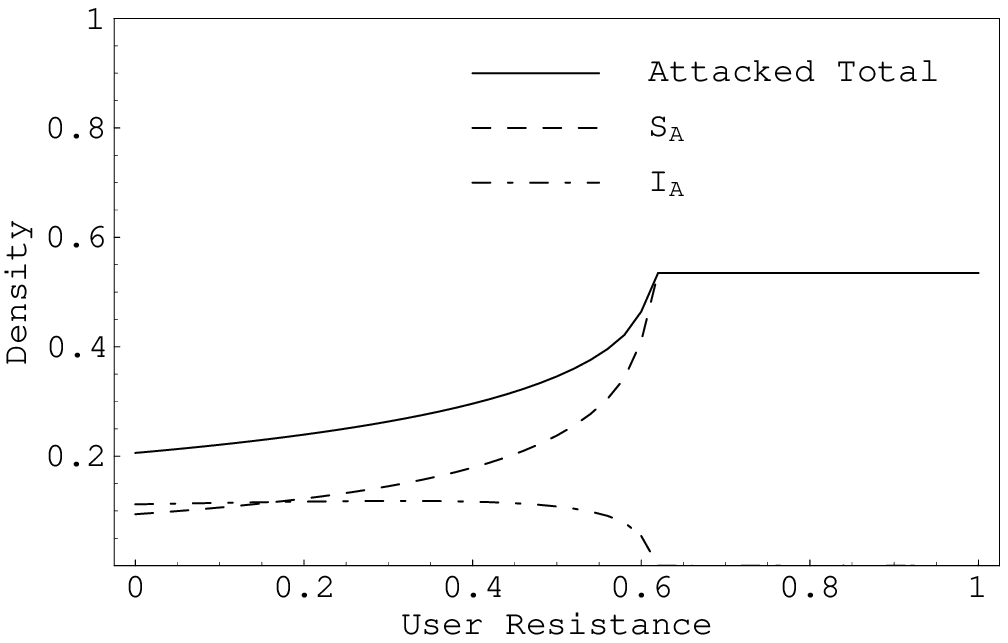}

\end{figure}

\pagebreak

Figure 2a:

\begin{figure}[h]

        \includegraphics{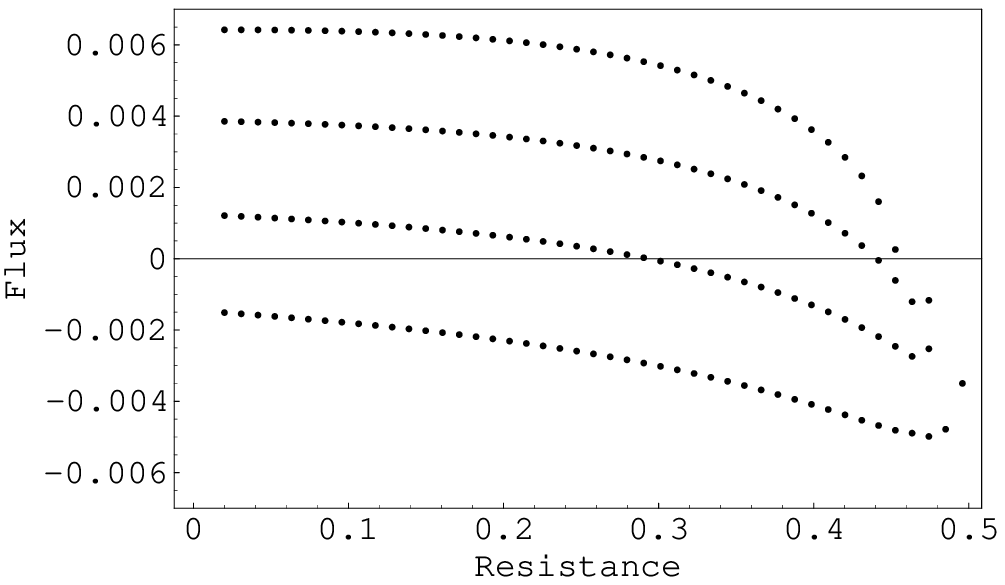}

\end{figure}

Figure 2b:

\begin{figure}[h]

        \includegraphics{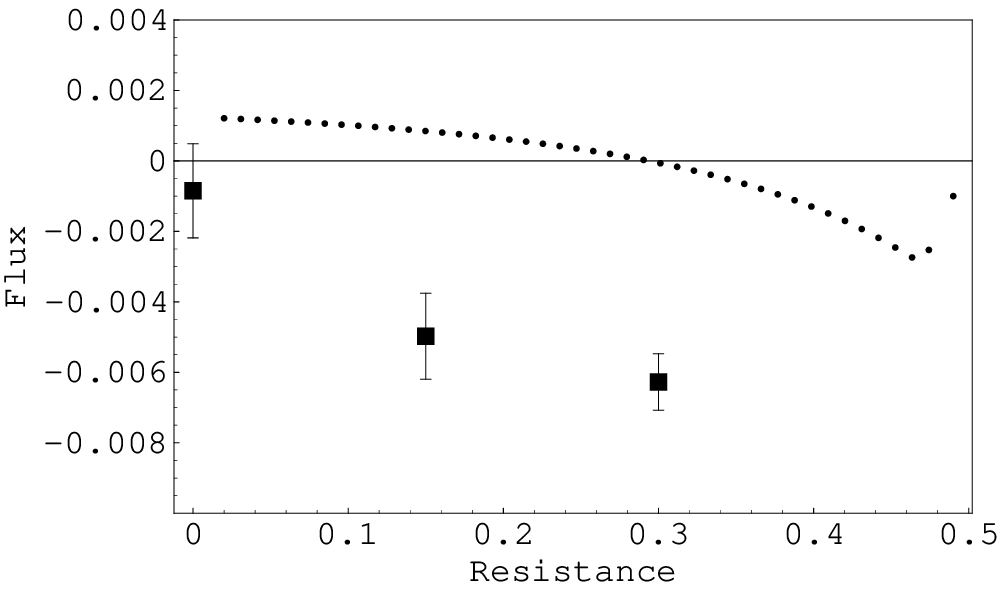}

\end{figure}

\pagebreak

Figure 3:

\begin{figure}[h]

        \includegraphics{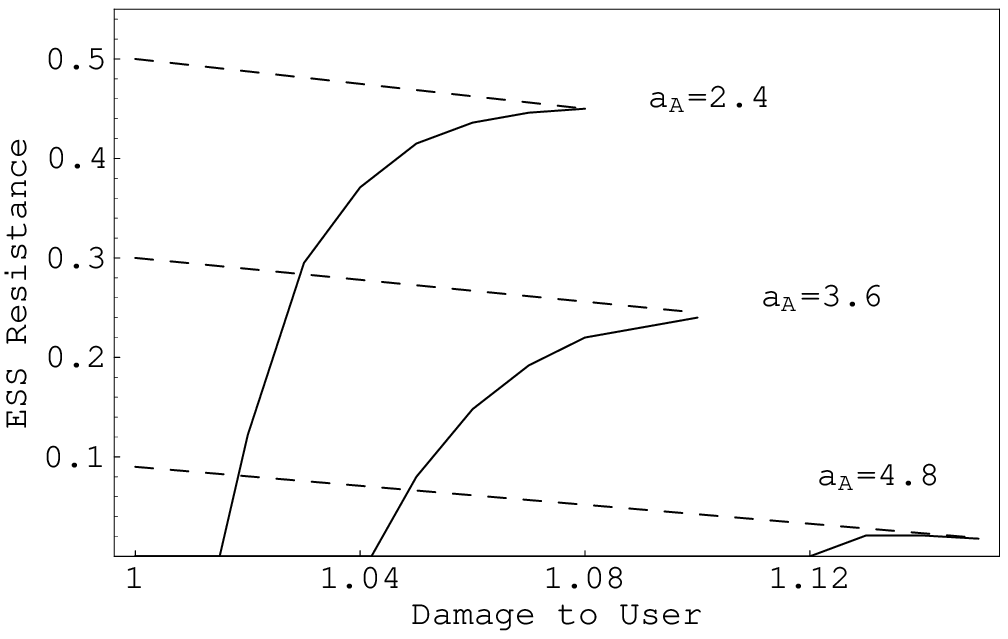}

\end{figure}

Figure 4:

\begin{figure}[h]

        \includegraphics{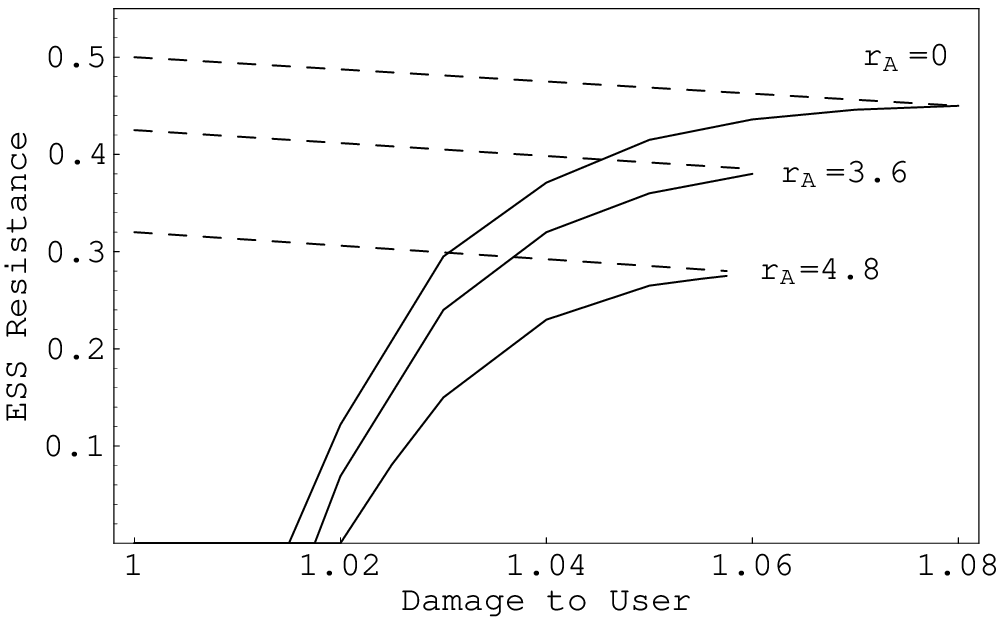}

\end{figure}

\pagebreak

Figure 5:

\begin{figure}[h]

        \includegraphics{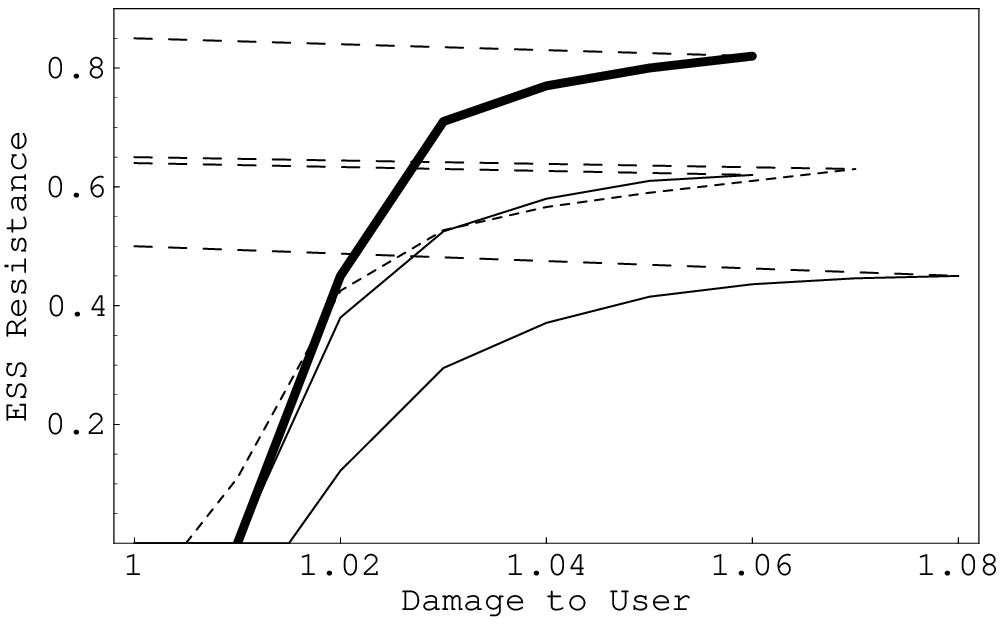}

\end{figure}

Figure 6:

\begin{figure}[h]

        \includegraphics{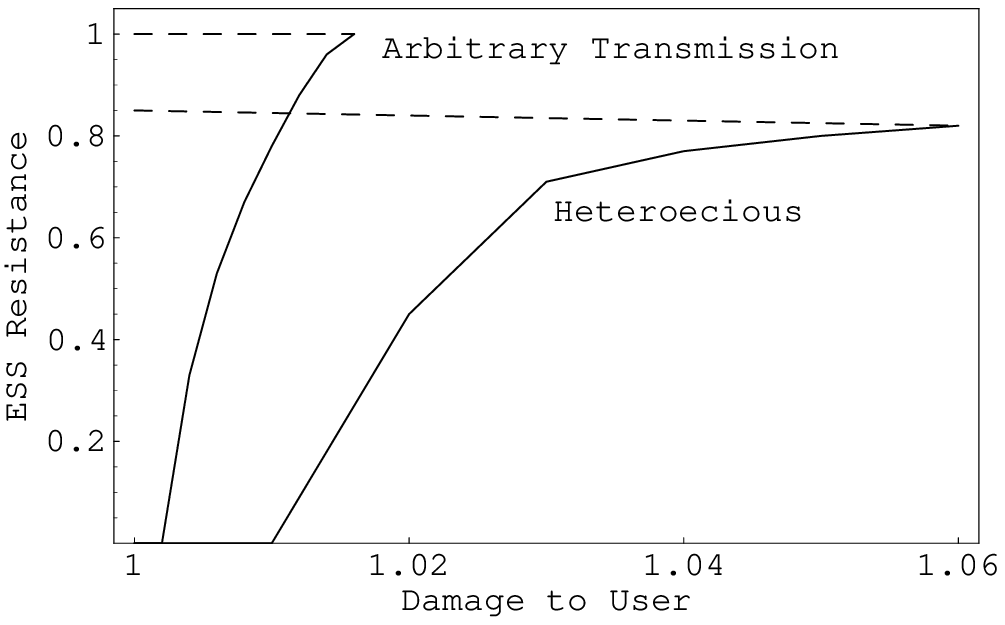}

\end{figure}

\pagebreak

\thispagestyle{myheadings}
\markright{}
\addcontentsline{toc}{chapter}{Bibliography}
\bibliography{thesis}
\thispagestyle{myheadings}
\markright{}
\bibliographystyle{plain}
\nocite{*}

\end{document}